# The Direct Shooting Method is a Complete Method

Sheng ZHANG, Jiang-Tao HUANG

(2021.02)

*Abstract:* **The direct shooting method is a classic approach for the solution of Optimal Control Problems (OCPs). It parameterizes the control variables and transforms the OCP to the Nonlinear Programming (NLP) problem to solve. This method is easy to use and it often introduces less parameters compared with all-variable parameterization method like the Pseudo-spectral (PS) method. However, it is long believed that its solution is not guaranteed to satisfy the optimality conditions of the OCP and the costates are not available in using this method. In this paper, we show that the direct shooting method may also provide the costate information, and it is proved that both the state and the costate solutions converge to the optimal as long as the control variable tends to the optimal, while the parameterized control may approach the optimal control with reasonable parameterization. This gives us the credit for the optimal control computation when employing the direct shooting method.**

## I. INTRODUCTION

Optimal control theory aims to determine the inputs to a dynamic system that optimize a specified performance index while satisfying constraints on the motion of the system. It is closely related to the engineering and has been widely studied [1]. Because of the complexity, Optimal Control Problems (OCPs) are usually solved numerically. Various numerical methods are developed and there are mainly two classes, namely, the indirect methods and the direct methods [2]. The indirect methods transform the OCP to a Boundary-value Problem (BVP) through the optimality conditions, and may hence give the accurate solution [3]. However, the indirect methods often suffer from the significant numerical difficulty due to the ill-conditioning of the Hamiltonian dynamics, that is, the stability of the costates dynamics is adverse to that of the states dynamics [4], and this makes the computation difficult without a good initial guess [5]. In a different manner, the direct methods discretize the control or/and state variables to obtain the Nonlinear Programming (NLP) problem. According to the number of variables parameterized, the direct method may be further divided into two types. The first is the all-variable parameterization method, in which all the state and control variables are parameterized. Typical methods of this type include the collation method [6] and the Legendre Pseudo-spectral (PS) method [7]. Such methods are also called the state and control parameterization method [8]. The second is the part-variable parameterization method, in which only part of the variables are parameterized. The direct shooting method [9], the differential inclusion method [10], the higher-order inclusion method [11], and the inversion-based optimization method [12] belong to this category. Generally the direct methods are easier to apply, whereas the optimality of the results is usually not guaranteed [13].

Researchers try to uncover the connection between the direct and the indirect methods. Positive studies, for example, include the direct collocation method [9], the Runge-Kutta discretization method [14], etc. In particular, the PS method [7], quite popular in

The authors are with the Computational Aerodynamics Institution, China Aerodynamics Research and Development Center, Mianyang, 621000, China. (e-mail: zszhangshengzs@hotmail.com).

recent decades, well blends the two types of methods, as it unifies the NLP and the BVP in a dualization view [15][16]. Such methods are termed as the complete method since the convergence of the primal variables and the dual variables are both guaranteed [17]. The complete method bridges the gap between the indirect methods and the direct methods, blurring their distinction. In contrast to the part-variable parameterization method, the all-variable parameterization methods have a facilitation that the costate's information may be obtained through dualization with the Lagrange multiplier rule and the Karush-Kuhn-Tucker conditions, and this enables further investigation on the costate mapping theory. For them, the feature of the complete method is that the discretization forms, from the transformed NLP and the resulting BVP respectively, are consistent when placed in the primal-dual view [18]. However, many all-variable parameterization methods are not complete methods [17].

The direct shooting method, also often called the control parameterization method, is a typical part-variable parameterization method, where only the control variables are parameterized. Due to its simplicity, it may be the most basic direct method solving OCPs and is widely utilized [2]. With the direct shooting method, the dynamic equations are satisfied by integrating the differential equations using time-marching algorithms, and the performance index is calculated in the way consistent with the numerical integrator that solves the differential equations [4]. Then the OCP is transformed to the NLP to be computed, using gradient methods or other congeneric algorithms. For example, with the dynamic method in Ref. [19], the gradient dynamic equation, which may be regarded as the continuous form of the iterative Sequential Quadratic Programming (SQP) method, is derived to seek the solution of the NLPs. Despite of the convenience and the simplicity, one adverse controversy on the direct shooting method is that this method cannot provide the costate information; therefore the optimality of the solution is not guaranteed. In Teo *et al.* [20], a set of costate system is presented in deriving the gradients on the control parameters. However, those variables therein will not converge to the costates of the optimal control solutions because they are considered for unconstrained functionals.

The Variation Evolving Method (VEM) is a recently proposed method for the optimal control computation [21]-[29]. Based on the infinite-dimensional Lyapunov principle, the Evolution Partial Differential Equation (EPDE) and Evolution Differential Equation (EDE), which describe the evolution of variables towards the optimal solution, are derived from the viewpoint of variation motion to improve the performance index. The VEM synthesizes the direct and indirect methods in that the established costate-free optimality conditions will be gradually met with theoretical guarantee. Upon the primary variables, the information on the costates and the multipliers may be reconstructed. In particular, from the principle of the VEM, the control variables may be parameterized and the resulting evolution dynamic equations for the parameters could be formulated. Like the infinite-dimensional case, it is expected that as the parameterized control approaches the optimal solution, the costates will tend to the optimal as well. An interesting question further arises that what is the relationship to the gradient dynamic equation for the NLP obtained with the direct shooting method. If they are consistent, then it suggests that the direct shooting method can provide the costate information and hence pave a way to investigate the completeness of the direct shooting method. This will be traversed in the following and the answer, to our encouragement, is positive.

Throughout the paper, our work is built upon the assumption that the solution for the optimization problem exists. We do not describe the existence conditions for the purpose of brevity. Relevant researches such as the Filippov-Cesari theorem are documented in Ref. [30]. In the following, first the preliminaries including the Lyapunov stability theory, the VEM, and the projection theory are reviewed. Then two parameterization forms are considered and the dynamic equations for the discretization parameters are derived from the view of VEM. Next, solutions using the direct shooting method are investigated, and it is shown that the gradient dynamic equations solving the resulting NLP are equivalent to the parameterization evolution equations arising from the VEM. Later, the costate mapping principle is established. After that, illustrative examples are solved to show the effectiveness and completeness of the direct shooting method.

## II. PRELIMINARIES

### A. Lyapunov dynamics stability theory

The Lyapunov stability theory investigates the dynamic behavior of states within a dynamic system, from the view of generalized energy [31]. It provides the theoretical basis for the method to be studied below. First the theory regarding the finite-dimensional dynamics is presented.

**Definition 1:** For a finite-dimensional continuous-time dynamic system like

$$\dot{\boldsymbol{x}} = \boldsymbol{f}(\boldsymbol{x}), \tag{1}$$

where $\boldsymbol{x} \in \mathbb{D} \in \mathbb{R}^n$ is the state, $\dot{\boldsymbol{x}} = \frac{\mathrm{d}\boldsymbol{x}}{\mathrm{d}t}$ is its time derivative, and $\boldsymbol{f}: \mathbb{D} \to \mathbb{R}^n$ is a vector function. $\mathbb{D}$ is a certain set. If $\hat{\boldsymbol{x}} \in \mathbb{D}$ satisfies $\boldsymbol{f}(\hat{\boldsymbol{x}}) = \boldsymbol{0}$, then $\hat{\boldsymbol{x}}$ is called an equilibrium solution.

**Definition 2:** The equilibrium solution $\hat{\boldsymbol{x}}$ is an asymptotically stable equilibrium solution in $\mathbb{D}$, if for any initial condition $\boldsymbol{x}(t)\big|_{t=0} = \boldsymbol{x}_0 \in \mathbb{D}$, there is $\lim_{t\to+\infty} \left\|\boldsymbol{x}(t) - \hat{\boldsymbol{x}}\right\|_2 = 0$. Here $\left\|\cdot\right\|_2$ denotes the 2-norm of vector.

**Lemma 1:** (see [31], with small adaptation) *For the continuous-time autonomous dynamic system* (1), *if there exists a continuously differentiable (If not, only except at* $\hat{\boldsymbol{x}}$ *) function* $V: \mathbb{D} \to \mathbb{R}$ *such that*

*i)* $V(\hat{\boldsymbol{x}}) = c$ *and* $V(\boldsymbol{x}) > c$ *in* $\mathbb{D} / \{\hat{\boldsymbol{x}}\}$,

*ii)* $\dot{V}(\boldsymbol{x}) \le 0$ *in* $\mathbb{D}$ *and* $\dot{V}(\boldsymbol{x}) < 0$ *in* $\mathbb{D} / \{\hat{\boldsymbol{x}}\}$.

*where* $c$ *is a constant. Then* $\boldsymbol{x} = \hat{\boldsymbol{x}}$ *is an asymptotically stable equilibrium solution in* $\mathbb{D}$.

For example, maybe $\boldsymbol{f}(\boldsymbol{x})$ in the dynamic system (1) satisfies $(\boldsymbol{x} - \hat{\boldsymbol{x}})^{\mathrm{T}} \boldsymbol{f}(\boldsymbol{x}) < 0$ for any $\boldsymbol{x} \ne \hat{\boldsymbol{x}}$, and then a Lyapunov function candidate can be constructed as $V = \frac{1}{2}(\boldsymbol{x} - \hat{\boldsymbol{x}})^{\mathrm{T}}(\boldsymbol{x} - \hat{\boldsymbol{x}})$, where the superscript '$\mathrm{T}$' denotes the transpose operator. The dynamics given by $\boldsymbol{f}(\boldsymbol{x})$ determines that $\dot{V} \le 0$ and $\boldsymbol{x}$ will converge to the equilibrium $\hat{\boldsymbol{x}}$. Figure 1 sketches the trajectory of some state in the stable dynamic system and the corresponding Lyapunov function value. No matter what the initial condition $\boldsymbol{x}_0$ is, as long as it falls into the stable region $\mathbb{D}$, the state $\boldsymbol{x}$ will approach the equilibrium $\hat{\boldsymbol{x}}$ gradually. Meanwhile, the 'energy' of the dynamic system, measured by the function $V$, will reach its minimum. Note that the constant $c$ in Lemma 1 is allowed to take non-zero value and this facilitates the analogy with the objective function in the optimization problem, whose value may not vanish at the optimum.

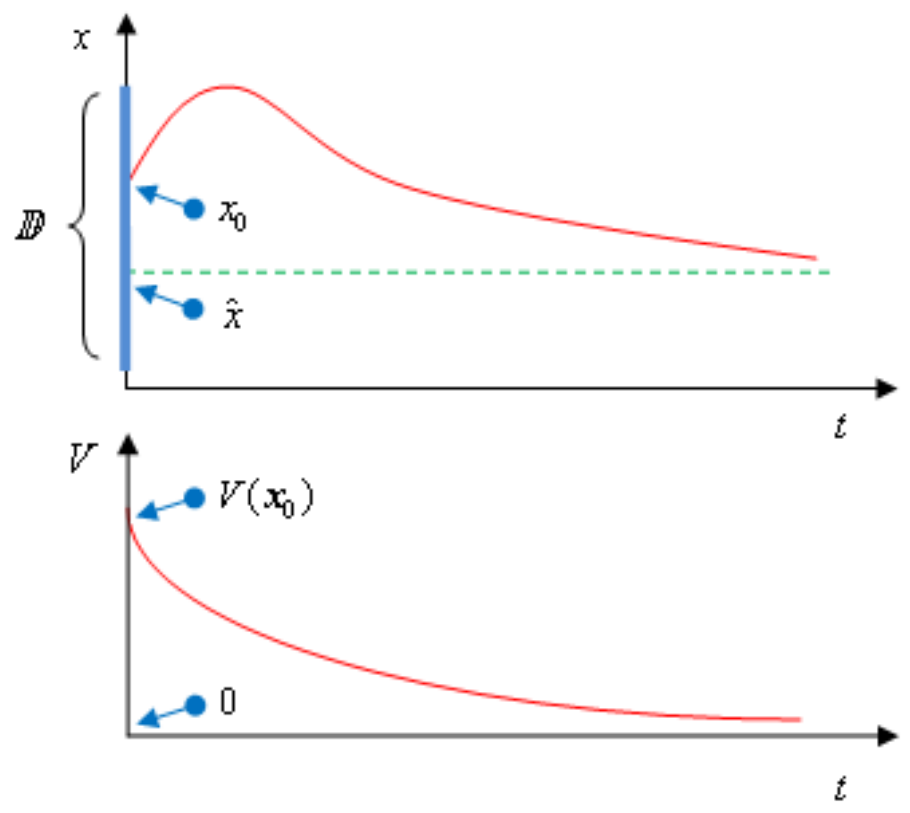

Fig. 1. Sketch for the state trajectory and the Lyapunov function value profile.

The Lyapunov principle may be generalized to the infinite-dimensional continuous-time dynamics as follows, which guarantees the effectiveness of the EPDE and the EDE formulated in the VEM.

**Definition 3:** For an infinite-dimensional dynamic system described by

$$\frac{\partial \boldsymbol{y}(x,t)}{\partial t} = \boldsymbol{f}(\boldsymbol{y}, \boldsymbol{p}, x) \tag{2}$$

$$\frac{\mathrm{d}\boldsymbol{p}(t)}{\mathrm{d}t} = \boldsymbol{\Gamma}(\boldsymbol{y}, \boldsymbol{p}) \tag{3}$$

where $t$ is the time. $x \in \mathbb{R}$ is the independent variable. $\boldsymbol{y}(x) \in \mathbb{D}(x) \in \mathbb{R}^n(x)$ is the function vector of $x$. $\boldsymbol{p} \in \Omega \in \mathbb{R}^m$ is a state vector. $\boldsymbol{f} : \mathbb{D}(x) \times \Omega \times \mathbb{R} \to \mathbb{R}^n(x)$ and $\boldsymbol{\Gamma} : \mathbb{D}(x) \times \Omega \to \mathbb{R}^m$ are vector functionals. $\mathbb{D}(x)$ is a certain function set and $\Omega$ is a certain number set. If $(\hat{\boldsymbol{y}}(x), \hat{\boldsymbol{p}}) \in (\mathbb{D}(x) \times \Omega)$ satisfies $\boldsymbol{f}\left(\hat{\boldsymbol{y}}(x), \hat{\boldsymbol{p}}, x\right) = \boldsymbol{0}$ and $\boldsymbol{\Gamma}\left(\hat{\boldsymbol{y}}(x), \hat{\boldsymbol{p}}\right) = \boldsymbol{0}$, then $(\hat{\boldsymbol{y}}(x), \hat{\boldsymbol{p}})$ is called an equilibrium solution.

**Definition 4:** The equilibrium solution $(\hat{\boldsymbol{y}}(x), \hat{\boldsymbol{p}})$ is an asymptotically stable equilibrium solution in $(\mathbb{D}(x) \times \Omega)$ if for any initial conditions $\boldsymbol{y}(x,t)\big|_{t=0} = \tilde{\boldsymbol{y}}(x) \in \mathbb{D}(x)$ and $\boldsymbol{p}(t)\big|_{t=0} = \tilde{\boldsymbol{p}} \in \Omega$, there is $\lim\limits_{t \to +\infty} \left\| (\boldsymbol{y}(x,t), \boldsymbol{p}(t)) - (\hat{\boldsymbol{y}}(x), \hat{\boldsymbol{p}}) \right\|_\infty = 0$, where $\left\| (\boldsymbol{y}(x), \boldsymbol{p}) \right\|_\infty := \max\left( \sup\limits_x \left( \sum_{i=1}^{n} y_i^2 \right)^{1/2}, \left( \sum_{i=1}^{m} p_i^2 \right)^{1/2} \right)$ is the defined supremum norm.

**Lemma 2:** *For the infinite-dimensional dynamic system* (2) *and* (3), *if there exists a continuously differentiable functional* $V : \mathbb{D}(x) \times \Omega \to \mathbb{R}$ *such that*

*i)* $V\left(\hat{\boldsymbol{y}}(x), \hat{\boldsymbol{p}}\right) = c$ *and* $V\left(\boldsymbol{y}(x), \boldsymbol{p}\right) > c$ *in* $(\mathbb{D}(x) \times \Omega) / \{(\hat{\boldsymbol{y}}(x), \hat{\boldsymbol{p}})\}$,

*ii)* $\dot{V}\left(\boldsymbol{y}(x), \boldsymbol{p}\right) \le 0$ *in* $(\mathbb{D}(x) \times \Omega)$ *and* $\dot{V}\left(\boldsymbol{y}(x), \boldsymbol{p}\right) < 0$ *in* $(\mathbb{D}(x) \times \Omega) / \{(\hat{\boldsymbol{y}}(x), \hat{\boldsymbol{p}})\}$,

*where* $c$ *is a constant. Then* $(\boldsymbol{y}(x), \boldsymbol{p}) = (\hat{\boldsymbol{y}}(x), \hat{\boldsymbol{p}})$ *is an asymptotically stable equilibrium solution in* $(\mathbb{D}(x) \times \Omega)$.

### *B. The VEM*

The VEM is a newly proposed method for the optimal solutions. Enlightened from the inverse consideration of the Lyapunov dynamics stability theory in the control field [31], the VEM analogizes the optimal solution of the OCP to the asymptotically stable equilibrium point of an infinite-dimensional dynamic system, and derives such dynamics to minimize a specific performance index that acts the Lyapunov functional. When employing Lemma 2 in the VEM, a virtual variation time, $\tau$, is introduced to describe the process that a variable evolves to the optimal solution under the dynamics governed by the EPDE and the EDE. Figure 2 illustrates the variation evolution of the control variables $\boldsymbol{u}(t)$ in the VEM to solve the OCP. Through the variation motion, the initial guess of control will evolve to the optimal solution, and the optimality conditions will be gradually achieved.

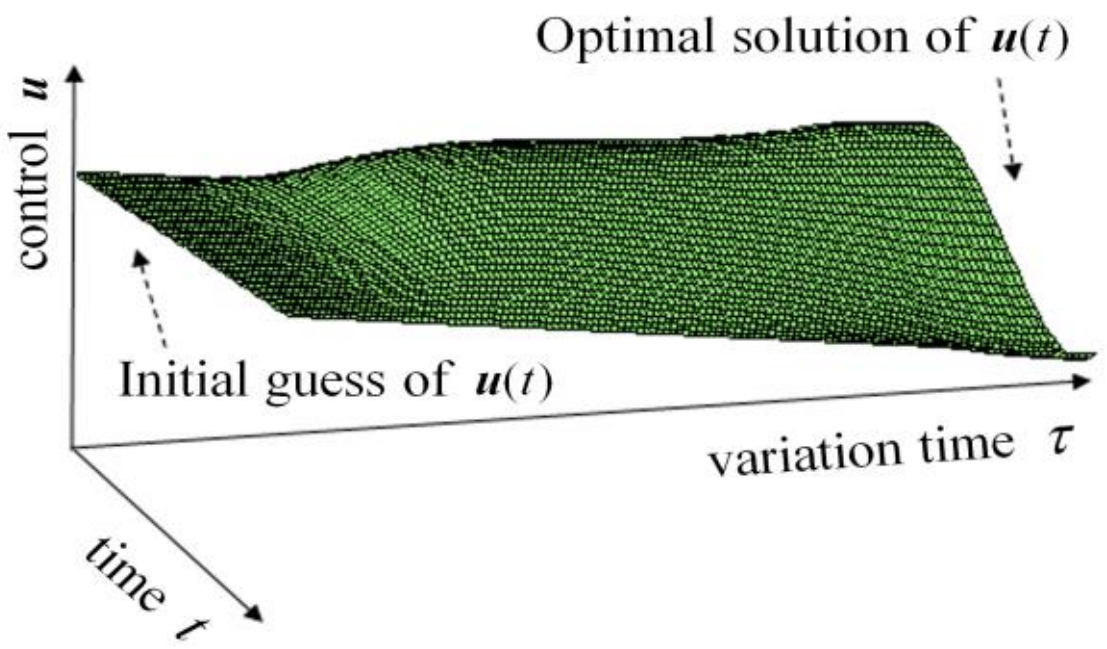

Fig. 2. The illustration of the control variable evolution along the variation time $\tau$ in the VEM.

In Refs. [21] and [22], besides the states and the controls, the costates are also employed in developing the EPDE, with arbitrary initial conditions. In Refs. [23]-[26], the VEM that uses only the state and control variables is developed. The costate-free optimality conditions are established and the evolution equations are defined within the feasible solution domain. In Refs. [27] and [28], the EPDE and the EDE are further generalized in the infeasible solution domain. In Ref. [29], a compact form EPDE that only concerns the control variables is proposed, and the state equations therein are satisfied by the numerical integration with time-marching Ordinary Differential Equation (ODE) integration methods. To seek the optimal solution, we need to solve the resulting infinite-dimensional Initial-value Problems (IVPs), defined by the EPDE and the EDEs with right definite conditions. Via the well-known semi-discrete method in the field of PDE numerical calculation [32], those equations are transformed to the finite-dimensional IVPs to be solved, with the common ODE integration methods.

An important result obtained with the VEM is the costate-free optimality conditions, which characterize the optimal solution only with the primary variables, namely the state and the control variables. The costate, the Lagrange multipliers, and the KKT multipliers may all be expressed by the primary variables, or determined by the equations related to the primary variables [26]. It is proved that the costate-free optimality conditions are equivalent to the classic optimality conditions which introduce the costates. The VEM produces the infinite-dimensional EPDE upon the Lyapunov functional. Actually, the method allows the variables to be parameterized in formulating the evolution equations. Then the Lyapunov functional degrades to the Lyapunov function, and the finite-dimensional Lyapunov stability theory will be employed to guarantee the effectiveness of the resulting equations.

### *C. Projection theory of Hilbert space*

The projection theory is involved in the following studies. Thus the main results will be reviewed in this subsection. Hilbert space is complete inner product space, and for its elements, the projection is the optimal approximation in a given subspace within based on the defined norm. For a separable Hilbert space, there is complete orthonormal system and each element in the space may be spanned by these bases. Given a subspace spanned by some bases (not necessarily orthonormal), we may investigate the least-square approximation of elements, i.e., the projection to the subspace. First, let us consider the Euclidean vector space.

**Lemma 3** [19]: *For the inner product* $<\boldsymbol{y}_1, \boldsymbol{y}_2> = \boldsymbol{y}_1^{\mathrm{T}}\boldsymbol{Q}\boldsymbol{y}_2$, *where* $\boldsymbol{y}_1, \boldsymbol{y}_2 \in \mathbb{R}^n$ *and* $\boldsymbol{Q}$ *is a positive-definite matrix, presuming there is a linear space* $\mathbb{S}$ *spanned by the row vectors of* $\boldsymbol{A} = [\boldsymbol{a}_1 \quad \boldsymbol{a}_2 \quad \dots \quad \boldsymbol{a}_m]$, *where* $\boldsymbol{A}$ *has full column rank, then the coordinates* $\boldsymbol{x}$ *in* $\mathbb{S}$ *that minimizes* $<\boldsymbol{z}-\boldsymbol{A}\boldsymbol{x}, \boldsymbol{z}-\boldsymbol{A}\boldsymbol{x}>$, *with* $\boldsymbol{z} \in \mathbb{R}^n$ *an arbitrary vector, is given by*

$$\boldsymbol{x} = (\boldsymbol{A}^{\mathrm{T}}\boldsymbol{Q}\boldsymbol{A})^{-1}\boldsymbol{A}^{\mathrm{T}}\boldsymbol{Q}\boldsymbol{z} \tag{4}$$

*and the projection of* $\boldsymbol{z}$ *into* $\mathbb{S}$ *is*

$$\mathrm{Pro}_{\mathbb{S}}(\boldsymbol{z}) = \boldsymbol{A}(\boldsymbol{A}^{\mathrm{T}}\boldsymbol{Q}\boldsymbol{A})^{-1}\boldsymbol{A}^{\mathrm{T}}\boldsymbol{Q}\boldsymbol{z} \tag{5}$$

*where* $\mathrm{Pro}_{\mathbb{S}}(\cdot)$ *denotes the projection operator to* $\mathbb{S}$ *.*

If the vectors in $\boldsymbol{A}$ are orthonormal under the given inner product, then there is $(\boldsymbol{A}^{\mathrm{T}}\boldsymbol{Q}\boldsymbol{A})^{-1}=\boldsymbol{I}$ , where $\boldsymbol{I}$ is the identity matrix. For the vector $\boldsymbol{z}\in\mathbb{R}^n$ , it may be divided into

$$\boldsymbol{z}=\mathrm{Pro}_{\mathbb{S}}(\boldsymbol{z})+\mathrm{Pro}_{\mathbb{S}^{\perp}}(\boldsymbol{z}) \tag{6}$$

where $\mathbb{S}^{\perp}$ denotes the orthogonal complement space in $\mathbb{R}^n$ . According to the projection theory, for any vector $\boldsymbol{v}\in\mathbb{S}^{\perp}$ , there is

$$<\mathrm{Pro}_{\mathbb{S}}(\boldsymbol{z}),\boldsymbol{v}>=0 \tag{7}$$

We now introduce the projection theory in the function space $L_2[t_0,t_f]$ . See

**Lemma 4** [33]: *For the inner product* $<g_1,g_2>=\int_{t_0}^{t_f}g_1(t)w(t)g_2(t)\mathrm{d}t$ *, where* $g_1$ *an* $g_2$ *are scalar functions,* $w(t)$ *is the positive weight function. Given a set of linearly independent bases* $a_1(t),a_2(t),...,a_n(t)$ *that span the function space* $\mathbb{S}(t)$ *, denote their assembling matrix by* $\boldsymbol{A}(t)=\begin{bmatrix}a_1(t) & a_2(t) & ... & a_n(t)\end{bmatrix}$ *. Then the coordinates* $\boldsymbol{x}$ *in the space* $\mathbb{S}(t)$ *, which minimizes* $<f-\boldsymbol{A}\boldsymbol{x},f-\boldsymbol{A}\boldsymbol{x}>$ *with* $f\in L_2[t_0,t_f]$ *an arbitrary vector function, is given by*

$$\boldsymbol{x}=\left(\int_{t_0}^{t_f}\boldsymbol{A}^{\mathrm{T}}(t)w(t)\boldsymbol{A}(t)\,\mathrm{d}t\right)^{-1}\left(\int_{t_0}^{t_f}\boldsymbol{A}^{\mathrm{T}}(t)w(t)\boldsymbol{f}(t)\,\mathrm{d}t\right) \tag{8}$$

*and the projection of* $f$ *into* $\mathbb{S}(t)$ *is*

$$\mathrm{Pro}_{\mathbb{S}}(f)=\boldsymbol{A}\left(\int_{t_0}^{t_f}\boldsymbol{A}^{\mathrm{T}}(t)w(t)\boldsymbol{A}(t)\,\mathrm{d}t\right)^{-1}\left(\int_{t_0}^{t_f}\boldsymbol{A}^{\mathrm{T}}(t)w(t)f(t)\,\mathrm{d}t\right) \tag{9}$$

Lemma 4 shows that the projection of the function has an analogous expression to that in the vector space. Note that Lemma 4 also aims at the general case where $a_1(t),a_2(t),...,a_n(t)$ are not required orthonormal under the defined inner product. If they are orthonormal with respect to $w(t)$, then the expressions will be simpler; that is, the matrix $\left(\int_{t_0}^{t_f}\boldsymbol{A}^{\mathrm{T}}(t)w(t)\boldsymbol{A}(t)\,\mathrm{d}t\right)$ will be an identity matrix. Extend the projection theory into the vector function space; we have

**Lemma 5**: *For the inner product* $<\boldsymbol{g}_1,\boldsymbol{g}_2>=\int_{t_0}^{t_f}\boldsymbol{g}_1^{\mathrm{T}}(t)\boldsymbol{W}(t)\boldsymbol{g}_2(t)\mathrm{d}t$ *, where* $\boldsymbol{g}_1$ *an* $\boldsymbol{g}_2$ *are vector functions,* $\boldsymbol{W}(t)$ *is the weight matrix function, being positive-definite within* $t\in[t_0,t_f]$ *. Given a set of linearly independent bases* $\boldsymbol{a}_1(t),\boldsymbol{a}_2(t),...,\boldsymbol{a}_n(t)$ *that span the vector function space* $\mathbb{S}(t)$ *, denote their assembling matrix by* $\boldsymbol{A}(t)=\begin{bmatrix}\boldsymbol{a}_1(t) & \boldsymbol{a}_2(t) & ... & \boldsymbol{a}_n(t)\end{bmatrix}$ *. Then the coordinates* $\boldsymbol{x}$ *in the space* $\mathbb{S}(t)$ *, which minimizes* $<\boldsymbol{f}-\boldsymbol{A}\boldsymbol{x},\boldsymbol{f}-\boldsymbol{A}\boldsymbol{x}>$ *, with* $\boldsymbol{f}$ *an arbitrary vector function, is given by*

$$\boldsymbol{x}=\left(\int_{t_0}^{t_f}\boldsymbol{A}^{\mathrm{T}}(t)\boldsymbol{W}(t)\boldsymbol{A}(t)\,\mathrm{d}t\right)^{-1}\left(\int_{t_0}^{t_f}\boldsymbol{A}^{\mathrm{T}}(t)\boldsymbol{W}(t)\boldsymbol{f}(t)\,\mathrm{d}t\right) \tag{10}$$

*and the projection of* $\boldsymbol{f}$ *into* $\mathbb{S}(t)$ *is*

$$\mathrm{Pro}_{\mathbb{S}}(\boldsymbol{f})=\boldsymbol{A}\left(\int_{t_0}^{t_f}\boldsymbol{A}^{\mathrm{T}}(t)\boldsymbol{W}(t)\boldsymbol{A}(t)\,\mathrm{d}t\right)^{-1}\left(\int_{t_0}^{t_f}\boldsymbol{A}^{\mathrm{T}}(t)\boldsymbol{W}(t)\boldsymbol{f}(t)\,\mathrm{d}t\right) \tag{11}$$

Lemma 5 is a natural extension of Lemma 4. Again Lemma 5 aims at the general case, which allows $\boldsymbol{a}_1(t),\boldsymbol{a}_2(t),...,\boldsymbol{a}_n(t)$ not orthonormal for the defined inner product. Similarly, a function $\boldsymbol{f}$ may be divided into

$$\boldsymbol{f} = \mathrm{Pro}_{\mathbb{S}}(\boldsymbol{f}) + \mathrm{Pro}_{\mathbb{S}^{\perp}}(\boldsymbol{f}) \tag{12}$$

and $\mathrm{Pro}_{\mathbb{S}}(\boldsymbol{f})$ is orthogonal to any functions within $\mathbb{S}^{\perp}(t)$.

## III. Parameterization in the Variation Evolution View

### *A. OCP definition*

In this paper, we consider the OCPs with terminal constraint that are defined as

**Problem 1**: Consider performance index of Bolza form

$$J = \varphi(\boldsymbol{x}(t_f), t_f) + \int_{t_0}^{t_f} L\left(\boldsymbol{x}(t), \boldsymbol{u}(t), t\right) \mathrm{d}t \tag{13}$$

subject to the dynamic equation

$$\dot{\boldsymbol{x}} = \boldsymbol{f}(\boldsymbol{x}, \boldsymbol{u}, t) \tag{14}$$

where $t \in \mathbb{R}$ is the time. $\boldsymbol{x} \in \mathbb{R}^n$ is the state vector and its elements belong to $C^2[t_0, t_f]$. $\boldsymbol{u} \in \mathbb{R}^m$ is the control vector and its elements belong to $C^1[t_0, t_f]$. The function $L: \mathbb{R}^n \times \mathbb{R}^m \times \mathbb{R} \to \mathbb{R}$ and its first-order partial derivatives are continuous with respect to $\boldsymbol{x}$, $\boldsymbol{u}$ and $t$. The function $\varphi: \mathbb{R}^n \times \mathbb{R} \to \mathbb{R}$ and its first-order and second-order partial derivatives are continuous with respect to $\boldsymbol{x}$ and $t$. The vector function $\boldsymbol{f}: \mathbb{R}^n \times \mathbb{R}^m \times \mathbb{R} \to \mathbb{R}^n$ and its first-order partial derivatives are continuous and Lipschitz in $\boldsymbol{x}$, $\boldsymbol{u}$ and $t$. The initial time $t_0$ is fixed and the terminal time $t_f$ is free. The initial and terminal boundary conditions are respectively prescribed as

$$\boldsymbol{x}(t_0) = \boldsymbol{x}_0 \tag{15}$$

$$\boldsymbol{g}\left(\boldsymbol{x}(t_f), t_f\right) = \boldsymbol{0} \tag{16}$$

where $\boldsymbol{g}: \mathbb{R}^n \times \mathbb{R} \to \mathbb{R}^q$ is a $q$ dimensional vector function with continuous first-order partial derivatives. Find the optimal solution $(\hat{\boldsymbol{x}}, \hat{\boldsymbol{u}})$ that minimizes $J$, i.e.

$$(\hat{\boldsymbol{x}}, \hat{\boldsymbol{u}}) = \arg\min(J) \tag{17}$$

### *B. Form 1 parameterization evolution equation*

In contrast to the infinite-dimensional evolution equations, we will derive the finite-dimensional parameterization evolution equation under the frame of the VEM. Consider Problem 1 within the quasi-feasible solution domain $\mathbb{D}_q$, in which the solutions of variables $\boldsymbol{x}$ and $\boldsymbol{u}$ satisfy Eqs. (14) and (15). The variation of the performance index $J$ given by Eq. (13), with respect to the virtual time $\tau$, is [23]

$$\frac{\delta J}{\delta \tau} = (\varphi_t + \varphi_{\boldsymbol{x}}{}^{\mathrm{T}} \boldsymbol{f} + L)\Big|_{t_f} \frac{\delta t_f}{\delta \tau} + \int_{t_0}^{t_f} \boldsymbol{p}_{\boldsymbol{u}}{}^{\mathrm{T}} \frac{\delta \boldsymbol{u}}{\delta \tau} \mathrm{d}t \tag{18}$$

where " $\Big|_{t_f}$ " means "evaluated at " $t_f$ ", and $\boldsymbol{p}_{\boldsymbol{u}}$ is

$$\boldsymbol{p}_{\boldsymbol{u}}(t) = L_{\boldsymbol{u}} + \boldsymbol{f}_{\boldsymbol{u}}{}^{\mathrm{T}}(t) \boldsymbol{\Phi}_o{}^{\mathrm{T}}(t_f, t) \varphi_{\boldsymbol{x}}(t_f) + \boldsymbol{f}_u{}^{\mathrm{T}} \left( \int_t^{t_f} \boldsymbol{\Phi}_o{}^{\mathrm{T}}(\sigma, t) L_{\boldsymbol{x}}(\sigma) \mathrm{d}\sigma \right) \tag{19}$$

The states are solved by

$$\boldsymbol{x}(t) = \boldsymbol{x}_0 + \int_{t_0}^{t} \boldsymbol{f}\left(\boldsymbol{x}, \boldsymbol{u}, s\right) \mathrm{d}s \tag{20}$$

which implies the variation evolution relation between $\frac{\delta \boldsymbol{x}}{\delta \tau}$ and $\frac{\delta \boldsymbol{u}}{\delta \tau}$, i.e.,

$$\frac{\delta \boldsymbol{x}}{\delta \tau} = \int_{t_0}^{t} \boldsymbol{\Phi}_o(t,s)\, \boldsymbol{f}_{\boldsymbol{u}}(s) \frac{\delta \boldsymbol{u}}{\delta \tau}(s)\,\mathrm{d}s \tag{21}$$

Like the procedure employed to formulate the generalized EPDE and EDE that are applicable in the infeasible solution domain [29], in deriving the evolution equations that ensure the achievement of the terminal boundary condition (16) from an infeasible solution, it is set that

$$\frac{\delta \boldsymbol{g}}{\delta \tau} + \boldsymbol{K}_{\boldsymbol{g}} \boldsymbol{g} = \boldsymbol{0} \tag{22}$$

where

$$\begin{aligned} \frac{\delta \boldsymbol{g}}{\delta \tau} &= \boldsymbol{g}_{\boldsymbol{x}_f} \frac{\delta \boldsymbol{x}(t_f)}{\delta \tau} + \boldsymbol{g}_{\boldsymbol{x}_f} \boldsymbol{f} \frac{\delta t_f}{\delta \tau} + \boldsymbol{g}_{t_f} \frac{\delta t_f}{\delta \tau} \\ &= \boldsymbol{g}_{\boldsymbol{x}_f} \int_{t_0}^{t_f} \boldsymbol{\Phi}_o(t_f,t)\, \boldsymbol{f}_{\boldsymbol{u}}(t) \frac{\delta \boldsymbol{u}}{\delta \tau}(t)\,\mathrm{d}t + (\boldsymbol{g}_{\boldsymbol{x}_f} \boldsymbol{f} + \boldsymbol{g}_{t_f}) \frac{\delta t_f}{\delta \tau} \end{aligned} \tag{23}$$

Then we construct the following augmented performance, which is created for the viable evolution equations (See Refs. [23] and [24] for detail).

$$J_t = J_{t1} + J_{t2} + \boldsymbol{\pi}^{\mathrm{T}} \left( \frac{\delta \boldsymbol{g}}{\delta \tau} + \boldsymbol{K}_{\boldsymbol{g}} \boldsymbol{g} \right) \tag{24}$$

where $\boldsymbol{\pi} \in \mathbb{R}^q$ is the Lagrange multipliers adjoining the variation of Eq. (22), $\boldsymbol{K}_{\boldsymbol{g}}$ is a $q \times q$ dimensional positive-definite matrix, and

$$J_{t1} = (\varphi_t + \varphi_{\boldsymbol{x}}{}^{\mathrm{T}} \boldsymbol{f} + L)\Big|_{t_f} \frac{\delta t_f}{\delta \tau} + \int_{t_0}^{t_f} \boldsymbol{p}_{\boldsymbol{u}}{}^{\mathrm{T}} \frac{\delta \boldsymbol{u}}{\delta \tau}\,\mathrm{d}t \tag{25}$$

$$J_{t2} = \frac{1}{2k_{t_f}} \left(\frac{\delta t_f}{\delta \tau}\right)^2 + \int_{t_0}^{t_f} \frac{1}{2} \left(\frac{\delta \boldsymbol{u}}{\delta \tau}\right)^{\mathrm{T}} \boldsymbol{K}^{-1} \frac{\delta \boldsymbol{u}}{\delta \tau}\,\mathrm{d}t \tag{26}$$

Note that the positive-definite matrix $\boldsymbol{K}$ may be time-varying, i.e., $\boldsymbol{K} = \boldsymbol{K}(t,\tau)$. Parameterize the control variable as

$$\boldsymbol{u}(t) = \boldsymbol{u}(t;\boldsymbol{p}) \tag{27}$$

where $\boldsymbol{p} \in \mathbb{R}^s$ is the parameters, then

$$\frac{\delta \boldsymbol{u}(t;\boldsymbol{p})}{\delta \tau} = \boldsymbol{u}_{\boldsymbol{p}}(t;\boldsymbol{p}) \frac{\mathrm{d}\boldsymbol{p}}{\mathrm{d}\tau} = \boldsymbol{u}_{p_1} \frac{\mathrm{d}p_1}{\mathrm{d}\tau} + \boldsymbol{u}_{p_2} \frac{\mathrm{d}p_2}{\mathrm{d}\tau} + \dots + \boldsymbol{u}_{p_s} \frac{\mathrm{d}p_s}{\mathrm{d}\tau} \tag{28}$$

where $\boldsymbol{u}_{\boldsymbol{p}} = \frac{\partial \boldsymbol{u}}{\partial \boldsymbol{p}} = \begin{bmatrix} \boldsymbol{u}_{p_1} & \boldsymbol{u}_{p_2} & \dots & \boldsymbol{u}_{p_s} \end{bmatrix}$ is the shorthand notation of the Jacobi partial derivative matrix and $\frac{\mathrm{d}\boldsymbol{p}}{\mathrm{d}\tau} = \begin{bmatrix} \frac{\mathrm{d}p_1}{\mathrm{d}\tau} & \frac{\mathrm{d}p_2}{\mathrm{d}\tau} & \dots & \frac{\mathrm{d}p_s}{\mathrm{d}\tau} \end{bmatrix}^{\mathrm{T}}$. Obviously, it is required that the parameters in $\boldsymbol{p}$ are irrelevant to each other, or there are redundant parameters. Therefore, the following assumption is introduced.

**Assumption 1**: The columns of matrix $\boldsymbol{u}_{\boldsymbol{p}}(t)$ are mutually independent; that is, there does not exist a non-zero constant vector of $\boldsymbol{a}$ such that $\boldsymbol{u}_{\boldsymbol{p}}(t)\boldsymbol{a} = \boldsymbol{0}$ for any $t \in [t_0, t_f]$.

Now Eq. (24) may be reformulated as (with the variation operator $\delta$ replaced)

$$J_t = (\varphi_t + \varphi_{\boldsymbol{x}}{}^{\mathrm{T}} \boldsymbol{f} + L)\Big|_{t_f} \frac{\mathrm{d}t_f}{\mathrm{d}\tau} + \int_{t_0}^{t_f} \boldsymbol{p_u}{}^{\mathrm{T}} \boldsymbol{u_p} \frac{\mathrm{d}\boldsymbol{p}}{\mathrm{d}\tau} \mathrm{d}t + \frac{1}{2k_{t_f}} (\frac{\mathrm{d}t_f}{\mathrm{d}\tau})^2 + \int_{t_0}^{t_f} \frac{1}{2} (\frac{\mathrm{d}\boldsymbol{p}}{\mathrm{d}\tau})^{\mathrm{T}} \boldsymbol{u_p}{}^{\mathrm{T}} \boldsymbol{K}^{-1} \boldsymbol{u_p} \frac{\mathrm{d}\boldsymbol{p}}{\mathrm{d}\tau} \mathrm{d}t$$
$$+ \boldsymbol{\pi}^{\mathrm{T}} \left( \boldsymbol{g}_{\boldsymbol{x}_f} \int_{t_0}^{t_f} \boldsymbol{\Phi}_o(t_f, t) \boldsymbol{f_u}(t) \boldsymbol{u_p} \frac{\mathrm{d}\boldsymbol{p}}{\mathrm{d}\tau} \mathrm{d}t + (\boldsymbol{g}_{\boldsymbol{x}_f} \boldsymbol{f} + \boldsymbol{g}_{t_f}) \frac{\mathrm{d}t_f}{\mathrm{d}\tau} + \boldsymbol{K_g g} \right) \tag{29}$$

Use the extreme conditions, i.e., $\dfrac{\partial J_t}{\partial (\frac{\mathrm{d}p_i}{\mathrm{d}\tau})} = 0, (i = 1, 2, ..., s)$ and $\dfrac{\partial J_t}{\partial \left( \frac{\mathrm{d}t_f}{\delta\tau} \right)} = 0$; we may get

$$\frac{\mathrm{d}\boldsymbol{p}}{\mathrm{d}\tau} = -\boldsymbol{M_p}{}^{-1} (\boldsymbol{r}_{1\boldsymbol{p}} + \boldsymbol{\Gamma}_{1\boldsymbol{p}} \boldsymbol{\pi}) \tag{30}$$

and

$$\frac{\mathrm{d}t_f}{\mathrm{d}\tau} = -k_{t_f} \left( \varphi_t + \varphi_{\boldsymbol{x}}{}^{\mathrm{T}} \boldsymbol{f} + L + \boldsymbol{\pi}^{\mathrm{T}} (\boldsymbol{g}_{\boldsymbol{x}_f} \boldsymbol{f} + \boldsymbol{g}_{t_f}) \right)\Big|_{t_f} \tag{31}$$

where

$$\boldsymbol{M_p} = \int_{t_0}^{t_f} \boldsymbol{u_p}{}^{\mathrm{T}}(t) \boldsymbol{K}^{-1} \boldsymbol{u_p}(t) \mathrm{d}t \tag{32}$$

$$\boldsymbol{r}_{1\boldsymbol{p}} = \int_{t_0}^{t_f} \boldsymbol{u_p}{}^{\mathrm{T}}(t) \boldsymbol{p_u} \mathrm{d}t \tag{33}$$

$$\boldsymbol{\Gamma}_{1\boldsymbol{p}} = \left( \int_{t_0}^{t_f} \boldsymbol{u_p}{}^{\mathrm{T}}(t) \boldsymbol{f_u}{}^{\mathrm{T}}(t) \boldsymbol{\Phi}_o{}^{\mathrm{T}}(t_f, t) \mathrm{d}t \right) \boldsymbol{g}_{\boldsymbol{x}_f}{}^{\mathrm{T}} \tag{34}$$

and the Lagrange multiplier $\boldsymbol{\pi}$ may be obtained by substituting Eqs. (30) and (31) into Eq. (22) as

$$\boldsymbol{\pi} = -\boldsymbol{M}_{1\boldsymbol{\pi}}{}^{-1} \boldsymbol{r}_{1\boldsymbol{\pi}} \tag{35}$$

where

$$\boldsymbol{M}_{1\boldsymbol{\pi}} = \boldsymbol{\Gamma}_{1\boldsymbol{p}}{}^{\mathrm{T}} \boldsymbol{M}_p{}^{-1} \boldsymbol{\Gamma}_{1\boldsymbol{p}} + k_{t_f} (\boldsymbol{g}_{\boldsymbol{x}_f} \boldsymbol{f} + \boldsymbol{g}_{t_f})(\boldsymbol{g}_{\boldsymbol{x}_f} \boldsymbol{f} + \boldsymbol{g}_{t_f})^{\mathrm{T}}\Big|_{t_f} \tag{36}$$

$$\boldsymbol{r}_{1\boldsymbol{\pi}} = \boldsymbol{\Gamma}_{1\boldsymbol{p}}{}^{\mathrm{T}} \boldsymbol{M}_p{}^{-1} \boldsymbol{r}_{1p} + k_{t_f} (\boldsymbol{g}_{\boldsymbol{x}_f} \boldsymbol{f} + \boldsymbol{g}_{t_f})(\varphi_t + \varphi_{\boldsymbol{x}}{}^{\mathrm{T}} \boldsymbol{f} + L)\Big|_{t_f} - \boldsymbol{K}_g \boldsymbol{g} \tag{37}$$

**Remark 1**: Assumption 1 guarantees that $\boldsymbol{u_p}{}^{\mathrm{T}}(t_i) \boldsymbol{K}^{-1} \boldsymbol{u_p}(t_i)$ has a rank of $\min(m, s)$ for any time point during $[t_0, t_f]$. If for some time point of $t_1$, $\int_{t_0}^{t_1} \boldsymbol{u_p}{}^{\mathrm{T}}(t) \boldsymbol{K}^{-1} \boldsymbol{u_p}(t) \mathrm{d}t > \boldsymbol{0}$, then there must be $\boldsymbol{M_p} > \boldsymbol{0}$; or if there exist time points $t_1, t_2, ..., t_r$ such that $\sum_{i=1}^{r} \boldsymbol{u_p}{}^{\mathrm{T}}(t_i) \boldsymbol{K}^{-1} \boldsymbol{u_p}(t_i) > \boldsymbol{0}$, then there is $\boldsymbol{M_p} > \boldsymbol{0}$. For either case, the matrix $\boldsymbol{M_p}$ is ensured positive-definite and its invertibility is guaranteed.

The matrix $\boldsymbol{M}_{1\boldsymbol{\pi}}$ given by Eq. (36) and the vector $\boldsymbol{r}_{1\boldsymbol{\pi}}$ given by Eq. (37) are consistent with those given in Ref. [24], where

$$\boldsymbol{M_\pi} = \boldsymbol{g}_{\boldsymbol{x}_f} \left( \int_{t_0}^{t_f} \boldsymbol{\Phi}_o(t_f, t) \boldsymbol{f_u K f_u}{}^{\mathrm{T}} \boldsymbol{\Phi}_o{}^{\mathrm{T}}(t_f, t) \mathrm{d}t \right) \boldsymbol{g}_{\boldsymbol{x}_f}{}^{\mathrm{T}} + k_{t_f} (\boldsymbol{g}_{\boldsymbol{x}_f} \boldsymbol{f} + \boldsymbol{g}_{t_f})(\boldsymbol{g}_{\boldsymbol{x}_f} \boldsymbol{f} + \boldsymbol{g}_{t_f})^{\mathrm{T}}\Big|_{t_f} \tag{38}$$

$$\boldsymbol{r_\pi} = \boldsymbol{g}_{\boldsymbol{x}_f} \int_{t_0}^{t_f} \boldsymbol{\Phi}_o(t_f, t) \boldsymbol{f_u K p_u} \mathrm{d}t + k_{t_f} (\boldsymbol{g}_{\boldsymbol{x}_f} \boldsymbol{f} + \boldsymbol{g}_{t_f})(\varphi_t + \varphi_{\boldsymbol{x}}{}^{\mathrm{T}} \boldsymbol{f} + L)\Big|_{t_f} - \boldsymbol{K}_g \boldsymbol{g} \tag{39}$$

and the multiplier $\boldsymbol{\pi}$ therein is determined by

$$\boldsymbol{\pi} = -\boldsymbol{M_\pi}{}^{-1} \boldsymbol{r_\pi} \tag{40}$$

To demonstrate this, let us introduce an infinite parameterization as

$$\boldsymbol{u}(t)=\int_{t_0}^{T}\boldsymbol{p}(t)\delta_F(t-\tau)\,\mathrm{d}\,\tau \tag{41}$$

where $\delta_F(t)$ is the impulse function and $T>t_f$ is fixed. Then there is

$$\frac{\partial \boldsymbol{u}}{\partial\,\boldsymbol{p}|_t}=\begin{bmatrix}\int_{t_0}^{T}\delta_F(t-\tau)\,\mathrm{d}\,\tau & 0 & 0 & 0\\ 0 & \int_{t_0}^{T}\delta_F(t-\tau)\,\mathrm{d}\,\tau & 0 & 0\\ \cdots & \cdots & \cdots & \cdots\\ 0 & 0 & 0 & \int_{t_0}^{T}\delta_F(t-\tau)\,\mathrm{d}\,\tau\end{bmatrix}_{m\times m} \tag{42}$$

where $\boldsymbol{p}|_t$ denotes the parameter at time $t$. For such special parameterization, it may be found that

$$\boldsymbol{u_p}(t)\boldsymbol{M}_p{}^{-1}\boldsymbol{r}_{2\boldsymbol{p}}=\boldsymbol{K}\boldsymbol{f_u}{}^{\mathrm{T}}\boldsymbol{\Phi}^{\mathrm{T}}(t_f,t)\boldsymbol{g}_{\boldsymbol{x}_f}{}^{\mathrm{T}} \tag{43}$$

$$\boldsymbol{u_p}(t)\boldsymbol{M}_p{}^{-1}\boldsymbol{r}_{1\boldsymbol{p}}=\boldsymbol{K}\boldsymbol{p_u} \tag{44}$$

which means

$$\boldsymbol{\Gamma}_{1\boldsymbol{p}}{}^{\mathrm{T}}\boldsymbol{M}_p{}^{-1}\boldsymbol{\Gamma}_{1\boldsymbol{p}}=\boldsymbol{g}_{\boldsymbol{x}_f}\left(\int_{t_0}^{t_f}\boldsymbol{\Phi}_o(t_f,t)\,\boldsymbol{f_u}\boldsymbol{K}\boldsymbol{f_u}{}^{\mathrm{T}}\boldsymbol{\Phi}_o{}^{\mathrm{T}}(t_f,t)\,\mathrm{d}\,t\right)\boldsymbol{g}_{\boldsymbol{x}_f}{}^{\mathrm{T}} \tag{45}$$

and

$$\boldsymbol{\Gamma}_{1\boldsymbol{p}}{}^{\mathrm{T}}\boldsymbol{M}_p{}^{-1}\boldsymbol{r}_{1p}=\boldsymbol{g}_{\boldsymbol{x}_f}\int_{t_0}^{t_f}\boldsymbol{\Phi}_o(t_f,t)\,\boldsymbol{f_u}\boldsymbol{K}\boldsymbol{p_u}\,\mathrm{d}\,t \tag{46}$$

Under the similar assumptions and with the similar proof in Ref. [28], we may establish the effectiveness of the Eqs. (30) and (31) to get an approximate optimal solution.

**Assumption 2**: The matrix $\begin{bmatrix}\boldsymbol{\Gamma}_{1\boldsymbol{p}}\\ (\boldsymbol{g}_{\boldsymbol{x}_f}\boldsymbol{f}+\boldsymbol{g}_{t_f})^{\mathrm{T}}\big|_{t_f}\end{bmatrix}$ has full column rank.

When $\boldsymbol{M_p}$ is positive-definite, Assumption 2 guarantees the invertibility of $\boldsymbol{M}_{1\boldsymbol{\pi}}$ because it is positive-definite; see

$$\boldsymbol{M}_{1\boldsymbol{\pi}}=\boldsymbol{\Gamma}_{1\boldsymbol{p}}{}^{\mathrm{T}}\boldsymbol{M}_p{}^{-1}\boldsymbol{\Gamma}_{1\boldsymbol{p}}+k_{t_f}\,(\boldsymbol{g}_{\boldsymbol{x}_f}\boldsymbol{f}+\boldsymbol{g}_{t_f})(\boldsymbol{g}_{\boldsymbol{x}_f}\boldsymbol{f}+\boldsymbol{g}_{t_f})^{\mathrm{T}}\Big|_{t_f}=\begin{bmatrix}\boldsymbol{\Gamma}_{1\boldsymbol{p}}\\ (\boldsymbol{g}_{\boldsymbol{x}_f}\boldsymbol{f}+\boldsymbol{g}_{t_f})^{\mathrm{T}}\big|_{t_f}\end{bmatrix}^{\mathrm{T}}\begin{bmatrix}\boldsymbol{M}_p{}^{-1} & \boldsymbol{0}\\ \boldsymbol{0} & k_{t_f}\end{bmatrix}\begin{bmatrix}\boldsymbol{\Gamma}_{1\boldsymbol{p}}\\ (\boldsymbol{g}_{\boldsymbol{x}_f}\boldsymbol{f}+\boldsymbol{g}_{t_f})^{\mathrm{T}}\big|_{t_f}\end{bmatrix} \tag{47}$$

Then the existence of the solution for Eq. (35) is ensured.

**Assumption 3**: During the computation process, the multiplier parameters $\boldsymbol{\pi}$ in Eq. (35) is bounded by some constant $d$, i.e., $\|\boldsymbol{\pi}\|_2\le d$.

**Theorem 1:** *Presume Assumptions 2 and 3 hold. Solving the IVP with respect to* $\tau$, *defined by the parameterization evolution equations* (30) *and* (31) *with arbitrary initial conditions of* $\boldsymbol{p}|_{\tau=0}$ *and* $t_f|_{\tau=0}$, *when* $\tau\to+\infty$, *the parameterized control* $\boldsymbol{u}(t;\boldsymbol{p})$ *and the state* $\boldsymbol{x}$ *will satisfy the terminal boundary condition* (16) *and the following optimality conditions.*

$$\boldsymbol{r}_{1\boldsymbol{p}}+\boldsymbol{\Gamma}_{1\boldsymbol{p}}\boldsymbol{\pi}=\boldsymbol{0} \tag{48}$$

$$\left(\varphi_t+\varphi_{\boldsymbol{x}}{}^{\mathrm{T}}\boldsymbol{f}+L+\boldsymbol{\pi}^{\mathrm{T}}(\boldsymbol{g}_{\boldsymbol{x}_f}\boldsymbol{f}+\boldsymbol{g}_{t_f})\right)\Big|_{t_f}=0 \tag{49}$$

The proof is based on Lemma 1, and the Lypaunov function candidate within $\mathbb{D}_q$ is

$$V=\sqrt{\boldsymbol{g}^{\mathrm{T}}\boldsymbol{g}}+c_1J \tag{50}$$

where $J$ given by Eq. (13) and $c_1 < \dfrac{\min(\mathrm{eig}(\boldsymbol{K}_{\boldsymbol{g}}))}{d\max(\mathrm{eig}(\boldsymbol{K}_{\boldsymbol{g}}))}$. $\mathrm{eig}(\cdot)$ is the function of eigenvalue. By analyzing the derivative of $V$ to $\tau$, it may be shown that the equilibrium solution of the parameterization evolution equations (30) and (31) is asymptotically stable and hence the parameterized control and state solutions will achieve the feasibility condition (16) and the optimality conditions (48), (49). In particular, Eqs. (48) and (49) may be re-presented as

$$\begin{bmatrix} \boldsymbol{r}_{1\boldsymbol{p}} \\ (\varphi_t + \varphi_{\boldsymbol{x}}{}^{\mathrm{T}}\boldsymbol{f} + L)\big|_{t_f} \end{bmatrix} + \begin{bmatrix} \boldsymbol{\Gamma}_{1\boldsymbol{p}} \\ (\boldsymbol{g}_{\boldsymbol{x}_f}\boldsymbol{f} + \boldsymbol{g}_{t_f})^{\mathrm{T}}\big|_{t_f} \end{bmatrix}\boldsymbol{\pi} = \boldsymbol{0} \tag{51}$$

This means that the vector $\begin{bmatrix} \boldsymbol{r}_{1\boldsymbol{p}} \\ (\varphi_t + \varphi_{\boldsymbol{x}}{}^{\mathrm{T}}\boldsymbol{f} + L)\big|_{t_f} \end{bmatrix}$ is within the space spanned by the column vectors of $\begin{bmatrix} \boldsymbol{\Gamma}_{1\boldsymbol{p}} \\ (\boldsymbol{g}_{\boldsymbol{x}_f}\boldsymbol{f} + \boldsymbol{g}_{t_f})^{\mathrm{T}}\big|_{t_f} \end{bmatrix}$ for the approximate optimal solution. Furthermore, we may claim that

**Proposition 1**: *Given a set of parameterization as* $\boldsymbol{u}^{(k)}(t) = \boldsymbol{u}(t;\boldsymbol{p}^{(k)})$, $k = 1, 2, ..., k_u$, *where* $k_u$ *may be infinite or finite, assuming that there is only one optimal approximation to an target optimal control solution in the function space* $\boldsymbol{u}(t;\boldsymbol{p}^{(k)}), \boldsymbol{p}^{(k)} \in \mathbb{R}^{s^{(k)}}$, *if the representation capacity of* $\boldsymbol{u}^{(k)}(t)$ *may be increased as* $k$ *increases and* $\boldsymbol{u}^{(k_u)}(t)$ *completely characterizes the optimal control* $\hat{\boldsymbol{u}}$, *then the performance index will be improved monotonically and the parameterized control* $\boldsymbol{u}^{(k)}(t)$ *will tends to the optimal control* $\hat{\boldsymbol{u}}$.

*Proof.* Since given a parameterization of $\boldsymbol{u}^{(k)}(t)$, it may be derived that the optimal results will ultimately satisfy

$$\frac{\delta J}{\delta \tau} = -k_{t_f}\left(\varphi_t + \varphi_{\boldsymbol{x}}{}^{\mathrm{T}}\boldsymbol{f} + L + \boldsymbol{\pi}^{\mathrm{T}}(\boldsymbol{g}_{\boldsymbol{x}_f}\boldsymbol{f} + \boldsymbol{g}_{t_f})\right)^2\Big|_{t_f} - \int_{t_0}^{t_f}(\boldsymbol{r}_{1\boldsymbol{p}} + \boldsymbol{\Gamma}_{1\boldsymbol{p}}\boldsymbol{\pi})^{\mathrm{T}}\boldsymbol{M}_{\boldsymbol{p}}{}^{-1}(\boldsymbol{r}_{1\boldsymbol{p}} + \boldsymbol{\Gamma}_{1\boldsymbol{p}}\boldsymbol{\pi})\,\mathrm{d}t \tag{52}$$

This term equals zero at the approximate optimal solution of $\tilde{\boldsymbol{u}}^{(k)}(t)$, with the corresponding performance index denoted by $\tilde{J}^{(k)}$. For the parameterization of $\boldsymbol{u}^{(k+1)}(t)$ that has stronger representation capacity, obviously from $\boldsymbol{u}^{(k+1)}(t) = \tilde{\boldsymbol{u}}^{(k)}(t)$, there is

$$\frac{\delta J}{\delta \tau} \le 0 \tag{53}$$

Therefore $\tilde{J}^{(k+1)} \le \tilde{J}^{(k)}$. On the premise that there is only one optimal approximation to an target optimal control solution in the function space $\boldsymbol{u}(t;\boldsymbol{p}^{(k)}), \boldsymbol{p}^{(k)} \in \mathbb{R}^{s^{(k)}}$, we have $\boldsymbol{u}^{(k_u)}(t) = \hat{\boldsymbol{u}}$ when the function space $\boldsymbol{u}(t;\boldsymbol{p}^{(k_u)}), \boldsymbol{p}^{(k_u)} \in \mathbb{R}^{s^{(k_u)}}$ includes the optimal control $\hat{\boldsymbol{u}}$. ■

It is shown in Ref. [24] that the optimality conditions for Problem 1 include Eq. (49) and the following equation.

$$\boldsymbol{p}_{\boldsymbol{u}} + \boldsymbol{f}_{\boldsymbol{u}}{}^{\mathrm{T}}\boldsymbol{\Phi}^{\mathrm{T}}(t_f, t)\boldsymbol{g}_{\boldsymbol{x}_f}{}^{\mathrm{T}}\boldsymbol{\pi} = \boldsymbol{0} \tag{54}$$

Apparently, the optimality condition regarding the terminal time is achieved for the Form 1 parameterization, while with the similar analysis with Eq. (41), Eq. (48) is consistent with the costate-free optimality condition (54).

**Remark 2**: If the OCP defined in Problem 1 has a fixed terminal time, namely, $t_f$ is no longer a parameter, then only Eq. (30) will be used to get the approximate optimal solution, and then the multiplier is calculated without considering terms involving $(\boldsymbol{g}_{\boldsymbol{x}_f}\boldsymbol{f} + \boldsymbol{g}_{t_f})\big|_{t_f}$, which may be achieved by simply setting $k_{t_f} = 0$.

According to the definition in Eq. (32), the matrix $\boldsymbol{M}_{\boldsymbol{p}}$ may be changed with $\tau$. However, it may take constant value as well, which is advantageous to the computation. Consider the case where the control is parameterized through the linear combination of base functions $\xi_i(t)$ $(i=1,2,...,s)$.

$$\boldsymbol{u} = p_1\boldsymbol{\xi}_1 + p_2\boldsymbol{\xi}_2 + ... + p_s\boldsymbol{\xi}_s \tag{55}$$

Now $\boldsymbol{u}_{\boldsymbol{p}}$ is only function of time $t$ and independent of $\boldsymbol{p}$. When the terminal time $t_f$ is fixed and $\boldsymbol{K}$ is irrelevant to $\tau$, $\boldsymbol{M}_{\boldsymbol{p}}$ is a naturally constant matrix. When the terminal time $t_f$ is free, by setting $\boldsymbol{K}$ relevant to $t_f(\tau)$, $\boldsymbol{M}_{\boldsymbol{p}}$ would always be a constant matrix. For example, for a free terminal time OCP with a scalar control input, assume $t_0 = 0$. If we use the polynomials $t^i$, $i = 0,1,2,...,(s-1)$ to parameterize the control, i.e., $u = \sum_{i=1}^{s} p_i t^{i-1}$, there is

$$\boldsymbol{M}_{\boldsymbol{p}} = \int_{t_0}^{t_f} \begin{bmatrix} 1 \\ t \\ ... \\ t^{s-1} \end{bmatrix} \frac{1}{K} \begin{bmatrix} 1 & t & ... & t^{s-1} \end{bmatrix} \mathrm{d}t = \int_{t_0}^{t_f} \begin{bmatrix} \frac{1}{K} & \frac{1}{K}t & ... & \frac{1}{K}t^{s-1} \\ \frac{1}{K}t & ... & \frac{1}{K}t^{s-1} & ... \\ ... & \frac{1}{K}t^{s-1} & ... & \frac{1}{K}t^{2s-3} \\ \frac{1}{K}t^{s-1} & ... & \frac{1}{K}t^{2s-3} & \frac{1}{K}t^{2s-2} \end{bmatrix} \mathrm{d}t = \begin{bmatrix} m_1 & m_2 & ... & m_s \\ m_2 & ... & m_s & ... \\ ... & m_s & ... & m_{2s-2} \\ m_s & ... & m_{2s-2} & m_{2s-1} \end{bmatrix} \tag{56}$$

Since $\boldsymbol{M}_{\boldsymbol{p}}$ only has $2s-1$ independent terms, let

$$\frac{1}{K} = \begin{cases} a_1 & t \in [0, \frac{t_f}{2s-1}) \\ a_2 & t \in [\frac{t_f}{2s-1}, \frac{2t_f}{2s-1}) \\ ... & \\ a_{2s-1} & t \in [\frac{(2s-2)t_f}{2s-1}, t_f] \end{cases} \tag{57}$$

and $a_j\big|_{\tau=0} > 0$ $(j = 1,2,...,2s-1)$. Then

$$\begin{bmatrix} m_1 \\ m_2 \\ ... \\ m_{2s-1} \end{bmatrix} = \boldsymbol{\Lambda}(t_f)\boldsymbol{P} \begin{bmatrix} a_1 \\ a_2 \\ ... \\ a_{2s-1} \end{bmatrix} \tag{58}$$

where $\boldsymbol{\Lambda}(t_f) = \begin{bmatrix} t_f & & & \\ & t_f^2 & & \\ & & ... & \\ & & & t_f^{2s-1} \end{bmatrix}$ and $\boldsymbol{P} = \begin{bmatrix} \frac{1}{2s-1} & \frac{1}{2s-1} & ... & \frac{1}{2s-1} \\ \frac{1^2}{2(2s-1)^2} & \frac{2^2-1^2}{2(2s-1)^2} & ... & \frac{(2s-1)^2-(2s-2)^2}{2(2s-1)^2} \\ ... & ... & ... & ... \\ \frac{1^{2s-1}}{(2s-1)(2s-1)^{2s-1}} & \frac{2^{2s-1}-1^{2s-1}}{(2s-1)(2s-1)^{2s-1}} & ... & \frac{(2s-1)^{2s-1}-(2s-2)^{2s-1}}{(2s-1)(2s-1)^{2s-1}} \end{bmatrix}$. Thus, $\boldsymbol{M}_{\boldsymbol{p}}$ may be constant by setting that

$$\begin{bmatrix} a_1(\tau) \\ a_2(\tau) \\ \dots \\ a_{2s-1}(\tau) \end{bmatrix} = \boldsymbol{P}^{-1}\boldsymbol{\Lambda}^{-1}(t_f(\tau))\boldsymbol{\Lambda}(t_f\big|_{\tau=0})\boldsymbol{P}\left.\begin{bmatrix} a_1 \\ a_2 \\ \dots \\ a_{2s-1} \end{bmatrix}\right|_{\tau=0} \tag{59}$$

and $a_j(\tau)$ $(j=1,2,...,2s-1)$ are positive around $t_f\big|_{\tau=0}$. However, $\boldsymbol{M}_{\boldsymbol{p}}$ cannot be arbitrary constant positive-definite matrix.

*C. Form 2 parameterization evolution equation*

In Sec. III.B, the parameterization of $\boldsymbol{u}(t)$ is independent of the terminal time $t_f$. This is natural when considered in the physical sense of variational motion. However, from the view of parameterization, the effect from $t_f$ should be included when the terminal time $t_f$ is free. A popular treatment often regards the control values at particular time points as the parameters and employs the interpolation technique to approximate the control profile. Such parameterization has clear physical meaning and is used widely, for example, the interpolation with Lagrange polynomials.

$$\boldsymbol{u}(t) = \sum_{i=0}^{N} \boldsymbol{u}^{\mathrm{T}}(t_i) \prod_{j=0,\, j\neq i}^{N} \frac{(t-t_j)}{(t_i-t_j)} \tag{60}$$

where $\boldsymbol{p} = \begin{bmatrix} \boldsymbol{u}^{\mathrm{T}}(t_0) & \boldsymbol{u}^{\mathrm{T}}(t_1) & \dots & \boldsymbol{u}^{\mathrm{T}}(t_N) \end{bmatrix}^{\mathrm{T}}$ is the parameter, $t_0$ denotes the initial time and $t_i$ $(i=1,2,...,N)$ are usually related to the terminal time $t_f$. In particular, $t_N = t_f$. For the parameterization in the generalized form

$$\boldsymbol{u}(t) = \boldsymbol{u}(t;\boldsymbol{p},t_f) \tag{61}$$

there is

$$\frac{\delta \boldsymbol{u}(t;\boldsymbol{p},t_f)}{\delta\tau} = \boldsymbol{u}_{\boldsymbol{p}}(t;\boldsymbol{p},t_f)\frac{\mathrm{d}\boldsymbol{p}}{\mathrm{d}\tau} + \boldsymbol{u}_{t_f}(t;\boldsymbol{p},t_f)\frac{\mathrm{d}t_f}{\mathrm{d}\tau} \tag{62}$$

We also suppose that Assumption 1 hold in the parameterization, without requiring the columns of the matrix $[\boldsymbol{u}_{\boldsymbol{p}}(t) \quad \boldsymbol{u}_{t_f}(t)]$ mutually independent.

With the variation operator $\delta$ replaced by the differential operator $\mathrm{d}$, now Eq. (24) may be expanded as

$$\begin{aligned} J_t = {} & (\varphi_t + \varphi_{\boldsymbol{x}}{}^{\mathrm{T}}\boldsymbol{f} + L)\Big|_{t_f}\frac{\mathrm{d}t_f}{\mathrm{d}\tau} + \int_{t_0}^{t_f} \boldsymbol{p}_{\boldsymbol{u}}{}^{\mathrm{T}}(\boldsymbol{u}_{\boldsymbol{p}}\frac{\mathrm{d}\boldsymbol{p}}{\mathrm{d}\tau} + \boldsymbol{u}_{t_f}\frac{\mathrm{d}t_f}{\mathrm{d}\tau})\mathrm{d}t + \frac{1}{2k_{t_f}}(\frac{\mathrm{d}t_f}{\mathrm{d}\tau})^2 + \int_{t_0}^{t_f}\frac{1}{2}\begin{bmatrix} \left(\frac{\mathrm{d}\boldsymbol{p}}{\mathrm{d}\tau}\right)^{\mathrm{T}} & \frac{\mathrm{d}t_f}{\mathrm{d}\tau} \end{bmatrix}\begin{bmatrix} \boldsymbol{u}_{\boldsymbol{p}}{}^{\mathrm{T}} \\ \boldsymbol{u}_{t_f}{}^{\mathrm{T}} \end{bmatrix}\boldsymbol{K}^{-1}\begin{bmatrix} \boldsymbol{u}_{\boldsymbol{p}} & \boldsymbol{u}_{t_f} \end{bmatrix}\begin{bmatrix} \frac{\mathrm{d}\boldsymbol{p}}{\mathrm{d}\tau} \\ \frac{\mathrm{d}t_f}{\mathrm{d}\tau} \end{bmatrix}\mathrm{d}t \\ & + \boldsymbol{\pi}^{\mathrm{T}}\left(\boldsymbol{g}_{\boldsymbol{x}_f}\int_{t_0}^{t_f}\boldsymbol{\Phi}_o(t_f,t)\boldsymbol{f}_{\boldsymbol{u}}(t)(\boldsymbol{u}_{\boldsymbol{p}}\frac{\mathrm{d}\boldsymbol{p}}{\mathrm{d}\tau} + \boldsymbol{u}_{t_f}\frac{\mathrm{d}t_f}{\mathrm{d}\tau})\mathrm{d}t + (\boldsymbol{g}_{\boldsymbol{x}_f}\boldsymbol{f} + \boldsymbol{g}_{t_f})\frac{\mathrm{d}t_f}{\mathrm{d}\tau} + \boldsymbol{K}_{\boldsymbol{g}}\boldsymbol{g}\right) \end{aligned} \tag{63}$$

Similarly using the extreme conditions, we have

$$\frac{\mathrm{d}}{\mathrm{d}\tau}\begin{bmatrix} \boldsymbol{p} \\ t_f \end{bmatrix} = -\boldsymbol{M}_{\boldsymbol{p}t_f}{}^{-1}(\boldsymbol{r}_{2\boldsymbol{p}t_f} + \boldsymbol{\Gamma}_{2\boldsymbol{p}t_f}\boldsymbol{\pi}) \tag{64}$$

where

$$\boldsymbol{M}_{\boldsymbol{p}t_f} = \begin{bmatrix} \int_{t_0}^{t_f}\boldsymbol{u}_{\boldsymbol{p}}{}^{\mathrm{T}}(t)\boldsymbol{K}^{-1}\boldsymbol{u}_{\boldsymbol{p}}(t)\,\mathrm{d}t & \int_{t_0}^{t_f}\boldsymbol{u}_{\boldsymbol{p}}{}^{\mathrm{T}}\boldsymbol{K}^{-1}\boldsymbol{u}_{t_f}\,\mathrm{d}t \\ \int_{t_0}^{t_f}\boldsymbol{u}_{t_f}{}^{\mathrm{T}}\boldsymbol{K}^{-1}\boldsymbol{u}_{\boldsymbol{p}}\mathrm{d}t & \frac{1}{k_{t_f}} + \int_{t_0}^{t_f}(\boldsymbol{u}_{t_f})^{\mathrm{T}}\boldsymbol{K}^{-1}\boldsymbol{u}_{t_f}\,\mathrm{d}t \end{bmatrix} \tag{65}$$

$$\boldsymbol{r}_{2\boldsymbol{p}t_f} = \begin{bmatrix} \int_{t_0}^{t_f} \boldsymbol{u_p}^{\mathrm{T}}(t)\boldsymbol{p_u}\,\mathrm{d}t \\ \varphi_t + \varphi_{\boldsymbol{x}}^{\mathrm{T}}\boldsymbol{f} + L\big|_{t_f} + \int_{t_0}^{t_f} \boldsymbol{u}_{t_f}^{\mathrm{T}}(t)\boldsymbol{p_u}\,\mathrm{d}t \end{bmatrix} \tag{66}$$

$$\boldsymbol{\Gamma}_{2\boldsymbol{p}t_f} = \begin{bmatrix} \left(\int_{t_0}^{t_f} \boldsymbol{u_p}^{\mathrm{T}}(t)\,\boldsymbol{f_u}^{\mathrm{T}}(t)\boldsymbol{\Phi}_o^{\mathrm{T}}(t_f,t)\mathrm{d}t\right)\boldsymbol{g_{x_f}}^{\mathrm{T}} \\ (\boldsymbol{g_{x_f}}\boldsymbol{f} + \boldsymbol{g}_{t_f})^{\mathrm{T}}\big|_{t_f} + \left(\int_{t_0}^{t_f} \boldsymbol{u}_{t_f}^{\mathrm{T}}(t)\,\boldsymbol{f_u}^{\mathrm{T}}(t)\boldsymbol{\Phi}_o^{\mathrm{T}}(t_f,t)\mathrm{d}t\right)\boldsymbol{g_{x_f}}^{\mathrm{T}} \end{bmatrix} \tag{67}$$

The Lagrange multiplier $\boldsymbol{\pi}$ may be obtained by substituting Eq. (64) into Eq. (22) as

$$\boldsymbol{\pi} = -\boldsymbol{M}_{2\boldsymbol{\pi}}^{-1}\boldsymbol{r}_{2\boldsymbol{\pi}} \tag{68}$$

where

$$\boldsymbol{M}_{2\boldsymbol{\pi}} = \boldsymbol{\Gamma}_{2pt_f}^{\mathrm{T}}\boldsymbol{M}_{pt_f}^{-1}\boldsymbol{\Gamma}_{2pt_f} \tag{69}$$

$$\boldsymbol{r}_{2\boldsymbol{\pi}} = \boldsymbol{\Gamma}_{2pt_f}^{\mathrm{T}}\boldsymbol{M}_{pt_f}^{-1}\boldsymbol{r}_{2pt_f} - \boldsymbol{K}_g\boldsymbol{g} \tag{70}$$

**Remark 3**: Due to Assumption 1 and the existence of the term $\frac{1}{k_{t_f}}$, under the similar conditions in Remark 1, the matrix $\boldsymbol{M}_{\boldsymbol{p}t_f}$ is ensured positive-definite and its invertibility is guaranteed.

In the same way, we can verify Eq. (64) in searching an approximate optimal solution, upon necessary preconditions.

**Assumption 4**: The matrix $\boldsymbol{r}_{2pt_f}$ has full column rank.

Likewise, Assumption 4 guarantees the invertibility of $\boldsymbol{M}_{2\boldsymbol{\pi}}$ when $\boldsymbol{M}_{\boldsymbol{p}t_f}$ is positive-definite.

**Assumption 5**: During the computation process, the multiplier parameters $\boldsymbol{\pi}$ in Eq. (68) is bounded by some constant $d$, i.e., $\|\boldsymbol{\pi}\|_2 \le d$.

**Theorem 2:** *Presume Assumptions 4 and 5 hold. Solving the IVP with respect to* $\tau$, *defined by the discretization evolution equation* (64) *with arbitrary initial condition, when* $\tau \to +\infty$, *the parameterized control* $\boldsymbol{u}(t;\boldsymbol{p},t_f)$ *and the state* $\boldsymbol{x}$ *will satisfy the feasibility condition* (16) *and the following optimality condition*

$$\boldsymbol{r}_{2\boldsymbol{p}t_f} + \boldsymbol{\Gamma}_{2\boldsymbol{p}t_f}\boldsymbol{\pi} = \boldsymbol{0} \tag{71}$$

Eq. (71) shows that at the approximate optimal solution, the vector $\boldsymbol{r}_{2\boldsymbol{p}t_f}$ is within the space spanned by the column vectors of $\boldsymbol{\Gamma}_{2\boldsymbol{p}t_f}$. Similarly, it may be claimed that

**Proposition 2**: *Given a set of parameterization as* $\boldsymbol{u}^{(k)}(t) = \boldsymbol{u}(t;\boldsymbol{p}^{(k)},t_f)$, $k = 1,2,\ldots,k_u$, *where* $k_u$ *may be infinite or finite, assuming that there is only one optimal approximation to an target optimal control solution in the function space* $\boldsymbol{u}(t;\boldsymbol{p}^{(k)},t_f)$, $\boldsymbol{p}^{(k)} \in \mathbb{R}^{s^{(k)}}$, *if the representation capacity of* $\boldsymbol{u}^{(k)}(t)$ *may be increased and* $\boldsymbol{u}^{(k_u)}(t)$ *completely characterizes the optimal control* $\hat{\boldsymbol{u}}$, *then the performance index may be improved monotonically and the parameterized control* $\boldsymbol{u}^{(k)}(t)$ *will tends to the optimal control* $\hat{\boldsymbol{u}}$.

The condition (71) is consistent to the optimality condition of Problem 1, namely, Eqs. (49) and (54). To show this, introduce the infinite parameterization as follows.

$$\boldsymbol{u}(t)=\int_0^1 \boldsymbol{p}(\frac{t-t_0}{t_f-t_0})\delta_F(\frac{t-t_0}{t_f-t_0}-\tau)\,\mathrm{d}\tau \tag{72}$$

Then

$$\frac{\partial \boldsymbol{u}}{\partial \boldsymbol{p}\big|_t}=\begin{bmatrix} \int_0^1 \delta_F(\frac{t-t_0}{t_f-t_0}-\tau)\,\mathrm{d}\tau & 0 & 0 & 0 \\ 0 & \int_0^1 \delta_F(\frac{t-t_0}{t_f-t_0}-\tau)\,\mathrm{d}\tau & 0 & 0 \\ \cdots & \cdots & \cdots & \cdots \\ 0 & 0 & 0 & \int_0^1 \delta_F(\frac{t-t_0}{t_f-t_0}-\tau)\,\mathrm{d}\tau \end{bmatrix}_{s\times s} \tag{73}$$

$$\boldsymbol{u}_{t_f}(t)=-\frac{t-t_0}{t_f-t_0}\frac{\mathrm{d}}{\mathrm{d}t}\boldsymbol{u}(t) \tag{74}$$

Therefore the last component in Eq. (71) determines

$$\left(\varphi_t+\varphi_{\boldsymbol{x}}{}^{\mathrm{T}}\boldsymbol{f}+L+\boldsymbol{\pi}^{\mathrm{T}}(\boldsymbol{g}_{\boldsymbol{x}_f}\boldsymbol{f}+\boldsymbol{g}_{t_f})\right)\Big|_{t_f}+\int_{t_0}^{t_f}\left(\boldsymbol{p}_{\boldsymbol{u}}+\boldsymbol{f}_{\boldsymbol{u}}{}^{\mathrm{T}}(t)\boldsymbol{\Phi}_o{}^{\mathrm{T}}(t_f,t)\boldsymbol{g}^{\mathrm{T}}{}_{\boldsymbol{x}_f}\boldsymbol{\pi}\right)^{\mathrm{T}}\boldsymbol{u}_{t_f}(t)\,\mathrm{d}t=0 \tag{75}$$

which is equivalent to Eq. (49) while Eq. (54) is established by other components in Eq. (71).

Section III.B and Section III.C consider different parameterization of control. If they have the same representation capacity, namely, given a set of parameters $\boldsymbol{p}_1$, there exist parameter vector $\boldsymbol{p}_2$ and terminal time $t_f$ such that $\boldsymbol{u}(t;\boldsymbol{p}_1)=\boldsymbol{u}(t;\boldsymbol{p}_2,t_f)$, and vice versa, then intuitively when their corresponding optimality conditions are satisfied, there should be $\tilde{\boldsymbol{u}}(t;\boldsymbol{p}_1)=\tilde{\boldsymbol{u}}(t;\boldsymbol{p}_2,t_f)$ for the optimized solutions. Consider the 1-order polynomial for a scalar control variable as follows:

Case 1: $u_a(t;\boldsymbol{p}_a)=p_{a1}+p_{a2}t$

Case 2: $u_b(t;\boldsymbol{p}_b,t_f)=p_{b1}\dfrac{(t-t_f)}{(t_0-t_f)}+p_{b2}\dfrac{(t-t_0)}{(t_f-t_0)}=\dfrac{p_{b1}t_f-p_{b2}t_0}{(t_f-t_0)}+\dfrac{(p_{b2}-p_{b1})}{(t_f-t_0)}t$

For Case 1, there is

$$\frac{\partial u_a}{\partial p_{a1}}=1,\frac{\partial u_a}{\partial p_{a2}}=t$$

For Case 2, we have

$$\frac{\partial u_b}{\partial p_{b1}}=\frac{(t-t_f)}{(t_0-t_f)},\frac{\partial u_b}{\partial p_{b2}}=\frac{(t-t_0)}{(t_f-t_0)},\frac{\partial u_b}{\partial t_f}=\frac{(p_{b1}-p_{b2})(t-t_0)}{(t_f-t_0)^2}$$

which gives

$$\begin{bmatrix}\dfrac{\partial u_b}{\partial p_{b1}}\\ \dfrac{\partial u_b}{\partial p_{b2}}\\ \dfrac{\partial u_b}{\partial t_f}\end{bmatrix}=\begin{bmatrix}-\dfrac{t_f}{t_0-t_f} & \dfrac{1}{t_0-t_f}\\ -\dfrac{t_0}{t_f-t_0} & \dfrac{1}{t_f-t_0}\\ -\dfrac{(p_{b1}-p_{b2})t_0}{(t_f-t_0)^2} & \dfrac{p_{b1}-p_{b2}}{(t_f-t_0)^2}\end{bmatrix}\begin{bmatrix}\dfrac{\partial u_a}{\partial p_{a1}}\\ \dfrac{\partial u_a}{\partial p_{a2}}\end{bmatrix}$$

The proposition given below will prove that the optimized solutions will be the same and their optimality is equivalent. This implies $p_{a1}=\dfrac{p_{b1}t_f-p_{b2}t_0}{(t_f-t_0)}$ and $p_{a2}=\dfrac{(p_{b2}-p_{b1})}{(t_f-t_0)}$. Note that in Case 2 for this simple example, the partial derivative on $t_f$ is not independent of the derivatives on $p_{b1}$ and $p_{b2}$.

**Proposition 3**: *For the control parameterizations given by Eq.* (27) *and Eq.* (61)*, respectively, if they have the same representation capacity, then their optimized solutions are the same and the optimality conditions are equivalent.*

*Proof.* If the control parameterization given by Eq. (27) and Eq. (61) have the same representation capacity; namely, there always exist $\boldsymbol{p}_1$, $\boldsymbol{p}_2$, and $t_f$ such that the following equation holds for any $t\in[t_0,t_f]$.

$$\boldsymbol{u}(t;\boldsymbol{p}_1)=\boldsymbol{u}(t;\boldsymbol{p}_2,t_f) \tag{76}$$

then the optimization problem is defined within the same parameterized function space. When the optimized solutions for the corresponding parameterization are obtained, they must be the same; or there is contradiction. Now consider their optimality conditions. Based on Eq. (76), the first-order differential relation may be established within $t\in[t_0,t_f]$, that is

$$\boldsymbol{u}_{\boldsymbol{p}_1}(t;\boldsymbol{p}_1)\mathrm{d}\boldsymbol{p}_1=\begin{bmatrix}\boldsymbol{u}_{\boldsymbol{p}_2}(t;\boldsymbol{p}_2,t_f) & \boldsymbol{u}_{t_f}(t;\boldsymbol{p}_2,t_f)\end{bmatrix}\begin{bmatrix}\mathrm{d}\boldsymbol{p}_2\\ \mathrm{d}t_f\end{bmatrix} \tag{77}$$

Since such relation always holds, there are proper dimensional constant matrixes $\boldsymbol{U}_1$, $\boldsymbol{U}_2$, and $\boldsymbol{V}_1$ such that

$$\boldsymbol{u}_{\boldsymbol{p}_1}(t;\boldsymbol{p}_1)\begin{bmatrix}\boldsymbol{U}_1 & \boldsymbol{U}_2\end{bmatrix}=\begin{bmatrix}\boldsymbol{u}_{\boldsymbol{p}_2}(t;\boldsymbol{p}_2,t_f) & \boldsymbol{u}_{t_f}(t;\boldsymbol{p}_2,t_f)\end{bmatrix} \tag{78}$$

$$\boldsymbol{u}_{\boldsymbol{p}_1}(t;\boldsymbol{p}_1)=\boldsymbol{u}_{\boldsymbol{p}_2}(t;\boldsymbol{p}_2,t_f)\boldsymbol{V}_1 \tag{79}$$

Furthermore, the same representation capacity implies the same freedom of independent parameters. Then the partial derivatives on $\boldsymbol{p}_2$ and $t_f$ are linearly dependent, and hence there exists matrix $\boldsymbol{V}_2$ such that

$$\boldsymbol{u}_{t_f}(t;\boldsymbol{p}_2,t_f)=\boldsymbol{u}_{\boldsymbol{p}_2}(t;\boldsymbol{p}_2,t_f)\boldsymbol{V}_2 \tag{80}$$

We first consider the optimality conditions of Eq. (51) for the Form 1 parameterization, with the multiplier denoted by $\boldsymbol{\pi}_1$, namely

$$\begin{bmatrix}\boldsymbol{r}_{1\boldsymbol{p}}\\ (\varphi_t+\varphi_{\boldsymbol{x}}{}^{\mathrm{T}}\boldsymbol{f}+L)\big|_{t_f}\end{bmatrix}+\begin{bmatrix}\boldsymbol{\Gamma}_{1\boldsymbol{p}}\\ (\boldsymbol{g}_{\boldsymbol{x}_f}\boldsymbol{f}+\boldsymbol{g}_{t_f})^{\mathrm{T}}\big|_{t_f}\end{bmatrix}\boldsymbol{\pi}_1=\boldsymbol{0} \tag{81}$$

Through multiplying a matrix of $\boldsymbol{C}=\begin{bmatrix}\boldsymbol{U}_1{}^{\mathrm{T}} & \boldsymbol{0} & \boldsymbol{0} & \boldsymbol{0}\\ \boldsymbol{0} & \boldsymbol{U}_1{}^{\mathrm{T}} & \boldsymbol{0} & \boldsymbol{0}\\ \cdots & \cdots & \cdots & \cdots\\ \boldsymbol{U}_2{}^{\mathrm{T}} & \boldsymbol{U}_2{}^{\mathrm{T}} & \boldsymbol{U}_2{}^{\mathrm{T}} & 1\end{bmatrix}$ by the left, there is

$$\boldsymbol{C}\begin{bmatrix}\boldsymbol{r}_{1\boldsymbol{p}}\\ (\varphi_t+\varphi_{\boldsymbol{x}}{}^{\mathrm{T}}\boldsymbol{f}+L)\big|_{t_f}\end{bmatrix}+\boldsymbol{C}\begin{bmatrix}\boldsymbol{\Gamma}_{1\boldsymbol{p}}\\ (\boldsymbol{g}_{\boldsymbol{x}_f}\boldsymbol{f}+\boldsymbol{g}_{t_f})^{\mathrm{T}}\big|_{t_f}\end{bmatrix}\boldsymbol{\pi}_1=\boldsymbol{0}\Rightarrow\boldsymbol{r}_{2\boldsymbol{p}_2t_f}+\boldsymbol{\Gamma}_{2\boldsymbol{p}_2t_f}\boldsymbol{\pi}_1=\boldsymbol{0} \tag{82}$$

On the other hand, for the optimality condition (71) of the Form 2 parameterization, denote the multiplier with $\boldsymbol{\pi}_2$. Multiplying $\boldsymbol{D}=\begin{bmatrix} \boldsymbol{V}_1^{\mathrm{T}} & \boldsymbol{0} & \boldsymbol{0} & \boldsymbol{0} \\ \boldsymbol{0} & \boldsymbol{V}_1^{\mathrm{T}} & \boldsymbol{0} & \boldsymbol{0} \\ \cdots & \cdots & \cdots & \cdots \\ -\boldsymbol{V}_2^{\mathrm{T}} & -\boldsymbol{V}_2^{\mathrm{T}} & -\boldsymbol{V}_2^{\mathrm{T}} & 1 \end{bmatrix}$ by the left, we have

$$\boldsymbol{D}\boldsymbol{r}_{2\boldsymbol{p}_2 t_f}+\boldsymbol{D}\boldsymbol{\Gamma}_{2\boldsymbol{p}_2 t_f}\boldsymbol{\pi}_2=\boldsymbol{0}\Rightarrow\begin{bmatrix} \boldsymbol{r}_{1\boldsymbol{p}} \\ (\varphi_t+\varphi_{\boldsymbol{x}}{}^{\mathrm{T}}\boldsymbol{f}+L)\big|_{t_f} \end{bmatrix}+\begin{bmatrix} \boldsymbol{\Gamma}_{1\boldsymbol{p}} \\ (\boldsymbol{g}_{\boldsymbol{x}_f}\boldsymbol{f}+\boldsymbol{g}_{t_f})^{\mathrm{T}}\big|_{t_f} \end{bmatrix}\boldsymbol{\pi}_2=\boldsymbol{0} \tag{83}$$

Eqs. (82) and (83) show that their forms may be deduced from each other. Now it only needs to show that $\boldsymbol{\pi}_1$ and $\boldsymbol{\pi}_2$ are equal. Assuming that $\boldsymbol{\pi}_1\neq\boldsymbol{\pi}_2$, then from Eq. (81) and (83), there is

$$\begin{bmatrix} \boldsymbol{\Gamma}_{1\boldsymbol{p}} \\ (\boldsymbol{g}_{\boldsymbol{x}_f}\boldsymbol{f}+\boldsymbol{g}_{t_f})^{\mathrm{T}}\big|_{t_f} \end{bmatrix}(\boldsymbol{\pi}_1-\boldsymbol{\pi}_2)=\boldsymbol{0} \tag{84}$$

Since $\begin{bmatrix} \boldsymbol{\Gamma}_{1\boldsymbol{p}} \\ (\boldsymbol{g}_{\boldsymbol{x}_f}\boldsymbol{f}+\boldsymbol{g}_{t_f})^{\mathrm{T}}\big|_{t_f} \end{bmatrix}$ has full column rank (See Assumption 2), this causes contradiction and so we have $\boldsymbol{\pi}_1=\boldsymbol{\pi}_2$, which are of the coordinate meaning in the subspace according to Lemma 3. ■

**Remark 4**: From the proof of Proposition 3, it may be found that Assumption 2 and Assumption 4 have the same implication in essence.

**Remark 5**: If the OCP defined in Problem 1 has a fixed terminal time of $t_f$, then Eq. (64) degrades to Eq. (30).

*D. The global and the local parameterization*

Compared with the global parameterization, local parameterization on a short partitioned time interval may have better approximation capacity. The discretization developed in the preceding allows both the global parameterization and the local parameterization. Between the Form 1 and Form 2 parameterizations, the latter is usually utilized when applying the local parameterization for the OCP with free terminal time, because the discretized time points are often related to the terminal time. With the local parameterization, the partial derivative of $\dfrac{\partial \boldsymbol{u}}{\partial p_j}$ $(j=1,2,...,s)$ will vanish, i.e., be a zero vector, beyond the time interval where the parameter is effective, while the form of the evolution equations is maintained. Consider two typical local parameterizations through the linear interpolation and the step approximation as follows.

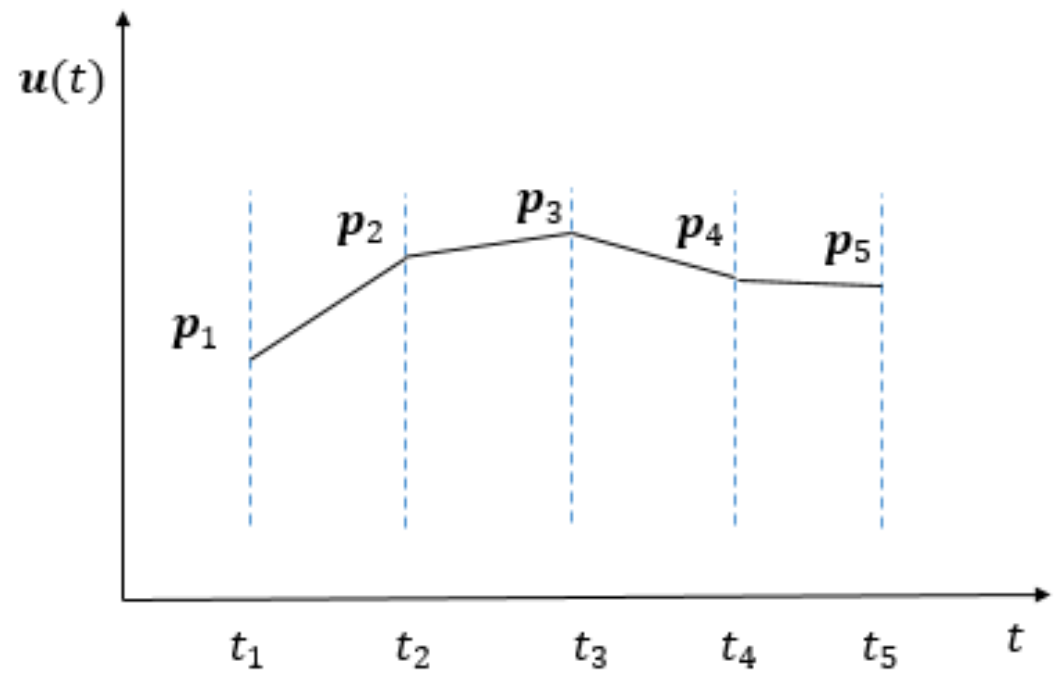

Fig. 3. The parameterization of the local linear interpolation.

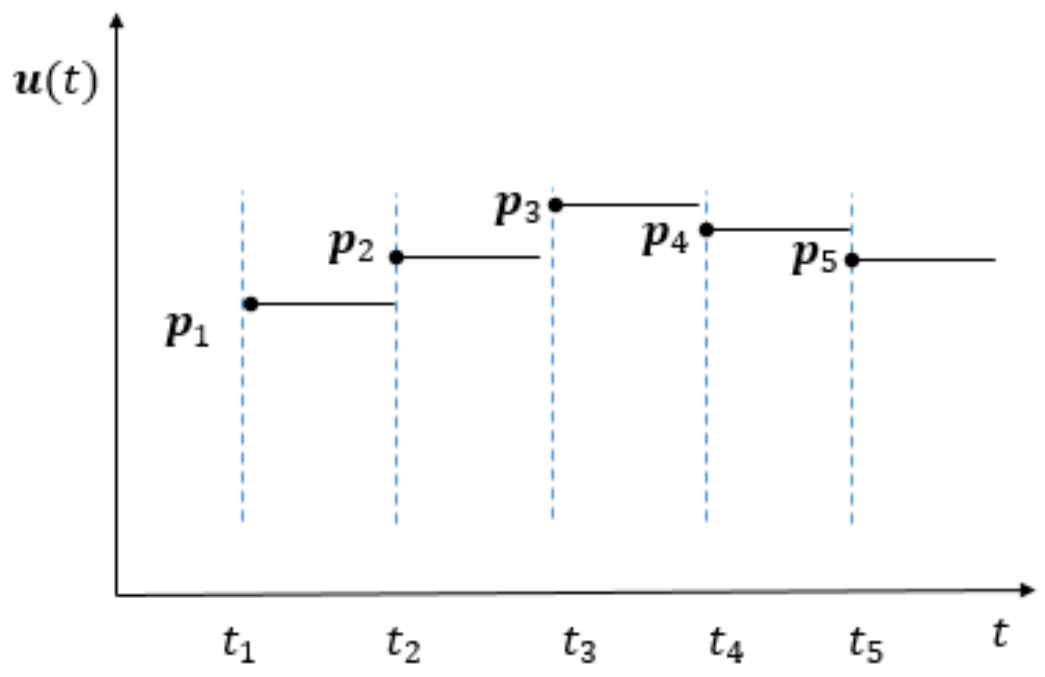

Fig. 4. The parameterization of the local step approximation.

For the linear interpolation, as shown in Fig. 3, with an equal partition for all components of $\boldsymbol{u}$ (unequal partition is also allowed), the control may be approximated as

$$\boldsymbol{u}(t) = \sum_{i=0}^{N} \boldsymbol{u}^{\mathrm{T}}(t_i)\,\chi_i(t) \tag{85}$$

where $\boldsymbol{p} = \left[\boldsymbol{u}^{\mathrm{T}}(t_0) \quad \boldsymbol{u}^{\mathrm{T}}(t_1) \quad \dots \quad \boldsymbol{u}^{\mathrm{T}}(t_N)\right]^{\mathrm{T}}$ is the parameter, $t_i = t_0 + \dfrac{t_f - t_0}{N} i \ \ (i = 0,1,\dots,N)$ are the discrete time points, and $\chi_i(t)$ is

$$\chi_0(t) = \begin{cases} \dfrac{t_1 - t}{t_1 - t_0}, & t \in [t_0, t_1] \\ 0, & \text{others} \end{cases}$$
$$\chi_i(t) = \begin{cases} \dfrac{t - t_{i-1}}{t_i - t_{i-1}}, & t \in [t_{i-1}, t_i] \\ \dfrac{t_{i+1} - t}{t_{i+1} - t_i}, & t \in [t_i, t_{i+1}] \\ 0, & \text{others} \end{cases} \quad (i = 1,2,\dots,N-1) \tag{86}$$
$$\chi_N(t) = \begin{cases} \dfrac{t_N - t}{t_N - t_{N-1}}, & t \in [t_{N-1}, t_N] \\ 0, & \text{others} \end{cases}$$

The Jacobi matrix $\boldsymbol{u}_{\boldsymbol{p}}$ is

$$\boldsymbol{u}_{\boldsymbol{p}} = \left[\chi_0(t)\boldsymbol{I}_{m\times m} \quad \chi_1(t)\boldsymbol{I}_{m\times m} \quad \dots \quad \chi_N(t)\boldsymbol{I}_{m\times m}\right] \tag{87}$$

where $\boldsymbol{I}_{m\times m}$ is the $m\times m$ dimensional identity matrix.

For the step approximation, as shown in Fig. 4, with an equal partition for all components (not necessarily), the control may be approximated as

$$\boldsymbol{u}(t) = \sum_{i=0}^{N-1} \boldsymbol{u}^{\mathrm{T}}(t_i)\,\chi_i(t) \tag{88}$$

where $\boldsymbol{p} = \left[\boldsymbol{u}^{\mathrm{T}}(t_0) \quad \boldsymbol{u}^{\mathrm{T}}(t_1) \quad \dots \quad \boldsymbol{u}^{\mathrm{T}}(t_{N-1})\right]^{\mathrm{T}}$, $t_i = t_0 + \dfrac{t_f - t_0}{N} i$, and $\chi_i(t)$ is

$$\chi_i(t) = \begin{cases} 1, & t \in [t_i, t_{i+1}) \\ 0, & \text{others} \end{cases} \quad (i = 0,1,\dots,N-1) \tag{89}$$

Obviously

$$\boldsymbol{u}_{\boldsymbol{p}} = \left[\chi_0(t)\boldsymbol{I}_{m\times m} \quad \chi_1(t)\boldsymbol{I}_{m\times m} \quad \dots \quad \chi_{N-1}(t)\boldsymbol{I}_{m\times m}\right] \tag{90}$$

## IV. Solution with the direct shooting method

### A. The transformed NLP

In solving the OCPs defined in Problem 1, the direct shooting method parameterizes the control variables and then transforms the OCPs to NLPs to solve. In this section, we will derive the computation formula from the view of NLP. Between the two types of parameterization, the Form 1 parameterization may apply to the problem with free or fixed terminal time, while the Form 2 parameterization usually aims at the case of free terminal time. Moreover, since the equivalence of the Form 1 and Form 2 parameterization, when they have the same representation capacity, has been established in the last section, in the following we use the Form 1 parameterization. Consider the transformation of Problem 1 upon the parameterization by Eq. (27). With the states expressed through the integral form of Eq.(20), it then takes the form

**Problem 2**: Consider the objective

$$J(\boldsymbol{\theta})=\varphi\left(\boldsymbol{x}_0+\int_{t_0}^{t_f}\boldsymbol{f}\left(\boldsymbol{x},\boldsymbol{u}(t;\boldsymbol{p}),t\right)\mathrm{d}t,t_f\right)+\int_{t_0}^{t_f}L\left(\boldsymbol{x}_0+\int_{t_0}^{t}\boldsymbol{f}\left(\boldsymbol{x},\boldsymbol{u}(s;\boldsymbol{p}),s\right)\mathrm{d}s,\boldsymbol{u}(t),t\right)\mathrm{d}t \tag{91}$$

subject to the equality constraint

$$\boldsymbol{g}\left(\int_{t_0}^{t_f}\boldsymbol{f}\left(\boldsymbol{x},\boldsymbol{u}(t;\boldsymbol{p}),t\right)\mathrm{d}t,t_f\right)=\boldsymbol{0} \tag{92}$$

where the functions $L:\mathbb{R}^n\times\mathbb{R}^m\times\mathbb{R}\to\mathbb{R}$, $\varphi:\mathbb{R}^n\times\mathbb{R}\to\mathbb{R}$, $\boldsymbol{f}:\mathbb{R}^n\times\mathbb{R}^m\times\mathbb{R}\to\mathbb{R}^n$, and $\boldsymbol{g}:\mathbb{R}^n\times\mathbb{R}\to\mathbb{R}^q$ conform the similar requirement as in Problem 1. Find the optimal parameter $\hat{\boldsymbol{\theta}}=\begin{bmatrix}\hat{\boldsymbol{p}}\\ \hat{t}_f\end{bmatrix}$ that minimizes $J$, i.e.

$$\hat{\boldsymbol{\theta}}=\arg\min(J) \tag{93}$$

### B. Gradient dynamic equation for NLP

The transformed NLP are usually solved with classic numerical iterative methods, namely

$$\boldsymbol{\theta}_{k+1}=\boldsymbol{M}_{\boldsymbol{\theta}}\left[\boldsymbol{\theta}_k\right] \tag{94}$$

where $\boldsymbol{M}_{\boldsymbol{\theta}}$ represents the mapping operator associated with some iterative optimization method. It has been pointed out that the gradient dynamic equation is a continuous version of the mechanism given by Eq. (94). A typical example is that Schropp proved that SQP method for the NLP with the equality constraints can be regarded as a variable step size Euler-Cauchy integration method applied to the gradient dynamic equation [34]. Therefore, the results will not be altered regardless of a continuous dynamic method or a discrete iterative method. Since the dynamic method also addresses the NLP in a continuous manner, before we move forward, the gradient dynamic equation that solves the general equality constrained NLP is presented.

**Lemma 6** [19]**:** *For the NLP defined with the objective function*

$$\min\quad f(\boldsymbol{\theta}) \tag{95}$$

*subject to the equality constraint*

$$\boldsymbol{h}(\boldsymbol{\theta})=\boldsymbol{0} \tag{96}$$

*where* $\boldsymbol{\theta}\in\mathbb{R}^n$ *is the optimization parameter vector.* $f:\mathbb{R}^n\to\mathbb{R}$ *is a scalar function with continuous first-order partial derivatives with respect to* $\boldsymbol{\theta}$. $\boldsymbol{h}:\mathbb{R}^n\to\mathbb{R}^s$ *is s-dimensional vector function with continuous first-order partial derivatives. Assume that* $\boldsymbol{h}_{\boldsymbol{\theta}}$ *is row full-rank; solve the IVP defined by the following gradient dynamic equation*

$$\frac{\mathrm{d}\boldsymbol{\theta}}{\mathrm{d}\tau} = -\boldsymbol{K}_{\boldsymbol{\theta}}(f_{\boldsymbol{\theta}} + \boldsymbol{h}_{\boldsymbol{\theta}}{}^{\mathrm{T}}\boldsymbol{\pi}_E) \tag{97}$$

*with arbitrary initial condition* $\boldsymbol{\theta}|_{\tau=0}$*, where the parameters* $\boldsymbol{\pi}_E \in \mathbb{R}^s$ *is determined by*

$$\boldsymbol{\pi}_E = -(\boldsymbol{h}_{\boldsymbol{\theta}}\boldsymbol{K}_{\boldsymbol{\theta}}\boldsymbol{h}_{\boldsymbol{\theta}}{}^{\mathrm{T}})^{-1}(\boldsymbol{h}_{\boldsymbol{\theta}}\boldsymbol{K}_{\boldsymbol{\theta}}\, \boldsymbol{f}_{\boldsymbol{\theta}} - \boldsymbol{K}_{\boldsymbol{h}}\boldsymbol{h}) \tag{98}$$

$\boldsymbol{K}_{\boldsymbol{\theta}}$ *is an* $n \times n$ *dimensional positive-definite matrix. Then under the assumption that* $\boldsymbol{\pi}_E$ *is bounded, when* $\tau \to +\infty$*,* $\boldsymbol{\theta}$ *will satisfy the feasibility conditions given by Eq.* (96) *and the optimality condition given by*

$$f_{\boldsymbol{\theta}} + \boldsymbol{h}_{\boldsymbol{\theta}}{}^{\mathrm{T}}\boldsymbol{\pi}_E = \boldsymbol{0} \tag{99}$$

*and then* $\boldsymbol{\pi}_E$ *also takes the value of* $-(\boldsymbol{h}_{\boldsymbol{\theta}}{}^{\mathrm{T}})^{+} f_{\boldsymbol{\theta}}$*, which is independent of* $\boldsymbol{K}_{\boldsymbol{\theta}}$*.*

Let $f(\boldsymbol{\theta}) = J(\boldsymbol{\theta})$ defined by Eq. (91), we first derive the gradient of $f_{\boldsymbol{\theta}}$

$$\mathrm{d}f = \varphi_{\boldsymbol{x}}{}^{\mathrm{T}}(t_f)\delta\boldsymbol{x}(t_f) + (\varphi_{\boldsymbol{x}}{}^{\mathrm{T}}\boldsymbol{f} + \varphi_t + L)\big|_{t_f}\,\mathrm{d}t_f + \int_{t_0}^{t_f}(L_{\boldsymbol{x}}{}^{\mathrm{T}}\delta\boldsymbol{x} + L_{\boldsymbol{u}}{}^{\mathrm{T}}\delta\boldsymbol{u})\,\mathrm{d}t \tag{100}$$

Because

$$\delta\boldsymbol{u}(t) = \boldsymbol{u}_{\boldsymbol{p}}(t)\mathrm{d}\boldsymbol{p} \tag{101}$$

and

$$\delta\boldsymbol{x}(t) = \int_{t_0}^{t}\boldsymbol{\Phi}_o(t,s)\,\boldsymbol{f}_{\boldsymbol{u}}(s)\delta\boldsymbol{u}(s)\,\mathrm{d}s = \left(\int_{t_0}^{t}\boldsymbol{\Phi}_o(t,s)\,\boldsymbol{f}_{\boldsymbol{u}}(s)\boldsymbol{u}_{\boldsymbol{p}}(s)\,\mathrm{d}s\right)\mathrm{d}\boldsymbol{p} \tag{102}$$

then there is

$$\mathrm{d}f = \varphi_{\boldsymbol{x}}{}^{\mathrm{T}}(t_f)\left(\int_{t_0}^{t_f}\boldsymbol{\Phi}_o(t_f,t)\,\boldsymbol{f}_{\boldsymbol{u}}(t)\boldsymbol{u}_{\boldsymbol{p}}(t)\,\mathrm{d}t\right)\mathrm{d}\boldsymbol{p} + (\varphi_{\boldsymbol{x}}{}^{\mathrm{T}}\boldsymbol{f} + \varphi_t + L)\big|_{t_f}\,\mathrm{d}t_f + \left(\int_{t_0}^{t_f}(L_{\boldsymbol{x}}{}^{\mathrm{T}}\left(\int_{t_0}^{t}\boldsymbol{\Phi}_o(t,s)\,\boldsymbol{f}_{\boldsymbol{u}}(s)\boldsymbol{u}_{\boldsymbol{p}}(s)\,\mathrm{d}s\right) + L_{\boldsymbol{u}}{}^{\mathrm{T}}\boldsymbol{u}_{\boldsymbol{p}})\,\mathrm{d}t\right)\mathrm{d}\boldsymbol{p} \tag{103}$$

Since

$$\int_{t_0}^{t_f}L_{\boldsymbol{x}}{}^{\mathrm{T}}\left(\int_{t_0}^{t}\boldsymbol{\Phi}_o(t,s)\,\boldsymbol{f}_{\boldsymbol{u}}(s)\boldsymbol{u}_{\boldsymbol{p}}(s)\,\mathrm{d}s\right)\mathrm{d}t = \int_{t_0}^{t_f}\left(\int_{t}^{t_f}L_{\boldsymbol{x}}{}^{\mathrm{T}}(\sigma)\boldsymbol{\Phi}_o(\sigma,t)\,\mathrm{d}\sigma\right)\boldsymbol{f}_{\boldsymbol{u}}(t)\boldsymbol{u}_{\boldsymbol{p}}(t)\,\mathrm{d}t \tag{104}$$

Eq. (103) may be re-written as

$$\mathrm{d}f = \varphi_{\boldsymbol{x}}{}^{\mathrm{T}}(t_f)\left(\int_{t_0}^{t_f}\boldsymbol{\Phi}_o(t_f,t)\,\boldsymbol{f}_{\boldsymbol{u}}(t)\boldsymbol{u}_{\boldsymbol{p}}(t)\,\mathrm{d}t\right)\mathrm{d}\boldsymbol{p} + (\varphi_{\boldsymbol{x}}{}^{\mathrm{T}}\boldsymbol{f} + \varphi_t + L)\big|_{t_f}\,\mathrm{d}t_f + \left(\int_{t_0}^{t_f}(\left(\int_{t}^{t_f}L_{\boldsymbol{x}}{}^{\mathrm{T}}(\sigma)\boldsymbol{\Phi}_o(\sigma,t)\,\mathrm{d}\sigma\right)\boldsymbol{f}_{\boldsymbol{u}}(t) + L_{\boldsymbol{u}}{}^{\mathrm{T}})\boldsymbol{u}_{\boldsymbol{p}}\,\mathrm{d}t\right)\mathrm{d}\boldsymbol{p} \tag{105}$$

Then we have

$$\begin{aligned} f_{\boldsymbol{\theta}} &= \begin{bmatrix} f_{\boldsymbol{p}} \\ f_{t_f} \end{bmatrix} = \begin{bmatrix} \left(\int_{t_0}^{t_f}\boldsymbol{u}_{\boldsymbol{p}}{}^{\mathrm{T}}(t)\,\boldsymbol{f}_{\boldsymbol{u}}{}^{\mathrm{T}}(t)\boldsymbol{\Phi}_o{}^{\mathrm{T}}(t_f,t)\mathrm{d}t\right)\varphi_{\boldsymbol{x}} + \left(\int_{t_0}^{t_f}\boldsymbol{u}_{\boldsymbol{p}}{}^{\mathrm{T}}(\,\boldsymbol{f}_{\boldsymbol{u}}{}^{\mathrm{T}}(t)\left(\int_{t}^{t_f}\boldsymbol{\Phi}_o{}^{\mathrm{T}}(\sigma,t)L_{\boldsymbol{x}}(\sigma)\,\mathrm{d}\sigma\right) + L_{\boldsymbol{u}})\,\mathrm{d}t\right) \\ (\varphi_{\boldsymbol{x}}{}^{\mathrm{T}}\boldsymbol{f} + \varphi_t + L)\big|_{t_f} \end{bmatrix} \\ &= \begin{bmatrix} \int_{t_0}^{t_f}\boldsymbol{u}_{\boldsymbol{p}}{}^{\mathrm{T}}(t)\,\boldsymbol{p}_{\boldsymbol{u}}(t)\mathrm{d}t \\ (\varphi_{\boldsymbol{x}}{}^{\mathrm{T}}\boldsymbol{f} + \varphi_t + L)\big|_{t_f} \end{bmatrix} \end{aligned} \tag{106}$$

where $\boldsymbol{p}_{\boldsymbol{u}}(t)$ is given by Eq. (19). On the other hand, for the equality constraint (92), there is

$$\begin{aligned} \mathrm{d}\boldsymbol{g} &= \boldsymbol{g}_{\boldsymbol{x}_f}\,\delta\boldsymbol{x}(t_f) + (\boldsymbol{g}_{\boldsymbol{x}_f}\,\boldsymbol{f} + \boldsymbol{g}_{t_f})\mathrm{d}t_f \\ &= \boldsymbol{g}_{\boldsymbol{x}_f}\left(\int_{t_0}^{t_f}\boldsymbol{\Phi}_o(t_f,t)\,\boldsymbol{f}_{\boldsymbol{u}}(t)\boldsymbol{u}_{\boldsymbol{p}}(t)\,\mathrm{d}t\right)\mathrm{d}\boldsymbol{p} + (\boldsymbol{g}_{\boldsymbol{x}_f}\,\boldsymbol{f} + \boldsymbol{g}_{t_f})\mathrm{d}t_f \end{aligned} \tag{107}$$

Thus the Jacobi matrix of $\boldsymbol{g}_{\boldsymbol{\theta}}$ is

$$\boldsymbol{g}_{\boldsymbol{\theta}} = \begin{bmatrix} \boldsymbol{g}_{\boldsymbol{p}} & \boldsymbol{g}_{t_f} \end{bmatrix} = \begin{bmatrix} \boldsymbol{g}_{\boldsymbol{x}_f}\left(\int_{t_0}^{t_f}\boldsymbol{\Phi}_o(t_f,t)\,\boldsymbol{f}_{\boldsymbol{u}}(t)\boldsymbol{u}_{\boldsymbol{p}}(t)\,\mathrm{d}t\right) & (\boldsymbol{g}_{\boldsymbol{x}_f}\,\boldsymbol{f} + \boldsymbol{g}_{t_f})\big|_{t_f} \end{bmatrix} \tag{108}$$

According to Eqs. (97), we may obtain the gradient dynamic equations as

$$\frac{\mathrm{d}}{\mathrm{d}\tau}\begin{bmatrix}\boldsymbol{p}\\ t_f\end{bmatrix}=-\boldsymbol{K}_{\boldsymbol{\theta}}\left(\boldsymbol{r}_{1\boldsymbol{p}t_f}+\boldsymbol{\Gamma}_{1\boldsymbol{p}t_f}\boldsymbol{\pi}_E\right) \tag{109}$$

with

$$\boldsymbol{\pi}_E=-(\boldsymbol{\Gamma}_{1\boldsymbol{p}t_f}{}^{\mathrm{T}}\boldsymbol{K}_{\boldsymbol{\theta}}\boldsymbol{\Gamma}_{1\boldsymbol{p}t_f})^{-1}(\boldsymbol{\Gamma}_{1\boldsymbol{p}t_f}{}^{\mathrm{T}}\boldsymbol{K}_{\boldsymbol{\theta}}\boldsymbol{r}_{1\boldsymbol{p}t_f}-\boldsymbol{K}_{\boldsymbol{g}}\boldsymbol{g}) \tag{110}$$

where $\boldsymbol{r}_{1\boldsymbol{p}t_f}=\begin{bmatrix}\int_{t_0}^{t_f}\boldsymbol{u}_{\boldsymbol{p}}{}^{\mathrm{T}}(t)\boldsymbol{p}_{\boldsymbol{u}}(t)\mathrm{d}t\\ (\varphi_{\boldsymbol{x}}{}^{\mathrm{T}}\boldsymbol{f}+\varphi_t+L)\big|_{t_f}\end{bmatrix}$ and $\boldsymbol{\Gamma}_{1\boldsymbol{p}t_f}=\begin{bmatrix}\left(\int_{t_0}^{t_f}\boldsymbol{u}_{\boldsymbol{p}}{}^{\mathrm{T}}\boldsymbol{f}_{\boldsymbol{u}}{}^{\mathrm{T}}\boldsymbol{\Phi}_o{}^{\mathrm{T}}(t_f,t)\mathrm{d}t\right)\boldsymbol{g}_{\boldsymbol{x}}{}^{\mathrm{T}}\\ (\boldsymbol{g}_{\boldsymbol{x}_f}\boldsymbol{f}+\boldsymbol{g}_{t_f})^{\mathrm{T}}\big|_{t_f}\end{bmatrix}$.

*C. Consistency between the gradient dynamic equation from NLP and parameterization evolution equation from VEM*

We now investigate the relation between the gradient dynamic equation, derived from the NLP, and the parameterization evolution equation, stemming from the VEM. Compare Eqs. (109) and (110) with Eqs. (30), (31) and (35). If they are equivalent, then the consistency is established.

First consider a special case. Set $\boldsymbol{K}_{\boldsymbol{\theta}}=\begin{bmatrix}\boldsymbol{K}_{\boldsymbol{p}} & \\ & k_{t_f}\end{bmatrix}$ in Eq. (109) and let $\boldsymbol{K}_{\boldsymbol{p}}=\left(\int_{t_0}^{t_f}\boldsymbol{u}_{\boldsymbol{p}}{}^{\mathrm{T}}(t)\boldsymbol{K}^{-1}\boldsymbol{u}_{\boldsymbol{p}}(t)\mathrm{d}t\right)^{-1}$; we have Eq. (110) and Eq. (35) are the same. It is also easy to find that the equation regarding $t_f$ in Eq. (109) is same to Eq. (31), and the equation regarding $\boldsymbol{p}$ in Eq. (109) is same to Eq. (30).

For the general form of $\boldsymbol{K}_{\boldsymbol{\theta}}$, if we set the penalty performance index in Eq. (26) as

$$J_{t2}=\int_{t_0}^{t_f}\frac{1}{2}\begin{bmatrix}\dfrac{\delta\boldsymbol{u}^{\mathrm{T}}}{\delta\tau} & \dfrac{\delta t_f}{\delta\tau}\end{bmatrix}\boldsymbol{K}^{-1}\begin{bmatrix}\dfrac{\delta\boldsymbol{u}}{\delta\tau}\\ \dfrac{\delta t_f}{\delta\tau}\end{bmatrix}\mathrm{d}t \tag{111}$$

we may derive that

$$\frac{\mathrm{d}}{\mathrm{d}\tau}\begin{bmatrix}\boldsymbol{p}\\ t_f\end{bmatrix}=-\boldsymbol{M}_{\boldsymbol{p}}{}^{-1}(\boldsymbol{r}_{1\boldsymbol{p}t_f}+\boldsymbol{\Gamma}_{1\boldsymbol{p}t_f}\boldsymbol{\pi}) \tag{112}$$

where

$$\boldsymbol{M}_{\boldsymbol{p}}=\int_{t_0}^{t_f}\frac{1}{2}\begin{bmatrix}\boldsymbol{u}_p{}^{\mathrm{T}} & \\ & 1\end{bmatrix}\boldsymbol{K}^{-1}\begin{bmatrix}\boldsymbol{u}_p & \\ & 1\end{bmatrix}\mathrm{d}t \tag{113}$$

and

$$\boldsymbol{\pi}=-\left(\boldsymbol{\Gamma}_{1\boldsymbol{p}t_f}{}^{\mathrm{T}}\boldsymbol{M}_{\boldsymbol{p}}{}^{-1}\boldsymbol{\Gamma}_{1\boldsymbol{p}t_f}\right)^{-1}\left(\boldsymbol{\Gamma}_{1\boldsymbol{p}t_f}{}^{\mathrm{T}}\boldsymbol{M}_{\boldsymbol{p}}{}^{-1}\boldsymbol{r}_{1\boldsymbol{p}t_f}-\boldsymbol{K}_g\boldsymbol{g}\right) \tag{114}$$

Let $\boldsymbol{K}_{\boldsymbol{\theta}}=\boldsymbol{M}_{\boldsymbol{p}}{}^{-1}$ in Eq. (109). Then the gradient dynamic equation from the NLP and the generalized parameterization evolution equation from the VEM are exactly the same in form.

It has been shown that in the parameterization evolution equation in Sec. III.B, the matrix of $\boldsymbol{M}_{\boldsymbol{p}}$ may be constant. Here the matrix $\boldsymbol{M}_{\boldsymbol{p}}$, for a more general case given in Eq. (113), may also take a constant value, independent of the virtual time $\tau$. Again consider the free terminal time OCP with a scalar control input in Sec. III.B. Let $\boldsymbol{K}^{-1}=\begin{bmatrix}a_{11} & a_{12}\\ a_{21} & a_{22}\end{bmatrix}$, which is positive definite. Then

$$\boldsymbol{M}_{\boldsymbol{p}}=\int_{t_0}^{t_f}\frac{1}{2}\begin{bmatrix}1 & 0\\ ... & ...\\ t^{s-1} & 0\\ 0 & 1\end{bmatrix}\begin{bmatrix}a_{11} & a_{12}\\ a_{21} & a_{22}\end{bmatrix}\begin{bmatrix}1 & ... & t^{s-1} & 0\\ 0 & ... & 0 & 1\end{bmatrix}\mathrm{d}t=\frac{1}{2}\begin{bmatrix}\int_{t_0}^{t_f}\begin{bmatrix}1\\ ...\\ t^{s-1}\end{bmatrix}a_{11}\begin{bmatrix}1 & ... & t^{s-1}\end{bmatrix}\mathrm{d}t & \int_{t_0}^{t_f}a_{12}\begin{bmatrix}1\\ ...\\ t^{s-1}\end{bmatrix}\mathrm{d}t\\ \int_{t_0}^{t_f}a_{21}\begin{bmatrix}1 & ... & t^{s-1}\end{bmatrix}\mathrm{d}t & \int_{t_0}^{t_f}a_{22}\mathrm{d}t\end{bmatrix} \tag{115}$$

We may set $a_{12}(t,\tau)$ and $a_{22}(t,\tau)$ such that $\int_{t_0}^{t_f}a_{12}\begin{bmatrix}1\\ ...\\ t^{s-1}\end{bmatrix}\mathrm{d}t$ and $\int_{t_0}^{t_f}a_{22}\mathrm{d}t$ are independent of $t_f$. However, $\boldsymbol{M}_{\boldsymbol{p}}$ cannot be an arbitrary positive-definite matrix. For example, $\boldsymbol{M}_{\boldsymbol{p}}$ cannot be diagonal. To the contrary, in the gradient dynamic equation (109), $\boldsymbol{K}_{\boldsymbol{\theta}}$ may be an arbitrary positive-definite matrix with constant value. This difference arises from the penalty index $J_{t2}$ in Eq. (111), which restricts the control variable variation. Actually, if the penalty index $J_{t2}$ is set as

$$J_{t2}=\frac{1}{2}\begin{bmatrix}\frac{\mathrm{d}\boldsymbol{p}^{\mathrm{T}}}{\mathrm{d}\tau} & \frac{\mathrm{d}t_f}{\mathrm{d}\tau}\end{bmatrix}\boldsymbol{K}^{-1}\begin{bmatrix}\frac{\mathrm{d}\boldsymbol{p}}{\mathrm{d}\tau}\\ \frac{\mathrm{d}t_f}{\mathrm{d}\tau}\end{bmatrix} \tag{116}$$

Then the equation same to Eq. (109) will be produced. For the parameterized function space, restricting the parameters through Eq. (116) allows a more flexible weight matrix, while their effects are the same in that they both guarantee a finite control solution which improves the performance index.

This section considers the Form 1 parameterization. Actually, with the Form 2 parameterization, it may also be shown that the parameterization evolution equation form the VEM and the gradient dynamic equation from the NLP are consistent as well.

## V. THE COSTATE MAPPING PRINCIPLE

### *A. The implicit costate mapping*

Originally the costate mapping principle shows the correspondence between the costates under the discretization of the BVP and the multipliers from the dualization of the NLP. Here we will show the mapping from the primary variables to the implicit costates and their convergence to the optimal counterpart, when employing the direct shooting method. For Problem 1, the classic optimality conditions contain the information about the costates and the Lagrange multipliers. Through the previous analysis on the direct shooting method, the Lagrange multiplier may be computed as Eq. (35) or Eq. (68) gives. Moreover, according to the VEM, the costates may be expressed as

$$\boldsymbol{\lambda}(t)=\boldsymbol{\Phi}_o{}^{\mathrm{T}}(t_f,t)\left(\varphi_{\boldsymbol{x}}(t_f)+\boldsymbol{g}_{\boldsymbol{x}_f}{}^{\mathrm{T}}\boldsymbol{\pi}\right)+\int_t^{t_f}\boldsymbol{\Phi}_o{}^{\mathrm{T}}(\sigma,t)L_{\boldsymbol{x}}(\sigma)\,\mathrm{d}\sigma \tag{117}$$

which satisfies the costate dynamic equation

$$\dot{\boldsymbol{\lambda}}(t)=-\boldsymbol{f}_{\boldsymbol{x}}{}^{\mathrm{T}}\boldsymbol{\lambda}(t)-L_{\boldsymbol{x}} \tag{118}$$

and the transversality conditions [35]

$$\boldsymbol{\lambda}(t_f)=\varphi_{\boldsymbol{x}}(t_f)+\boldsymbol{g}_{\boldsymbol{x}_f}{}^{\mathrm{T}}(t_f)\boldsymbol{\pi} \tag{119}$$

Note that the state solutions $\boldsymbol{x}$ and the multiplier $\boldsymbol{\pi}$ are determined by the parameterized control. Thus, Eq. (117) gives an implicit mapping to the costate from the discrete parameters.

In parameterizing the control profile, the Form 1 and Form 2 parameterizations are equivalent when they have the same representation capacity, while the Form 1 parameterization brings simpler results in that the evolution equations, on the terminal time and on the control parameters respectively, are decoupled. Actually, this parameterization has a more natural and direct relation to the infinite-dimensional evolution equations and the costate-free optimality conditions established in Ref. [24]. Thus, we again consider the Form 1 parameterization given by Eq. (27) in the following. We will reinvestigate the optimality conditions for the parameters, and further show that as the parameterized control profile approaches the optimal, the costate given by Eq. (117) will converge to the optimal, which also implies the convergence of the multiplier given by Eq. (35).

*B. Reinvestigation of the optimality conditions for parameters*

The function projection theory will be utilized to carry out the studies. First, the inner product involved is defined as

$$<\boldsymbol{g}_1, \boldsymbol{g}_2> = \int_{t_0}^{t_f} \boldsymbol{g}_1^{\mathrm{T}}(t)\boldsymbol{K}\boldsymbol{g}_2(t)\mathrm{d}t \tag{120}$$

where $\boldsymbol{g}_1$ and $\boldsymbol{g}_2$ are vector functions, and $\boldsymbol{K}$ is weight matrix used in Eq. (26). The upper limit and the lower limit of the integral are the initial time $t_0$ and the terminal time $t_f$ of Problem 1, respectively. Then we introduce the column vectors of $\boldsymbol{K}^{-1}\boldsymbol{u}_{\boldsymbol{p}}(t)$ as the bases and let $\mathbb{S}(t)$ denote the function subspace spanned. Note that the bases are linearly independent while they are not necessarily orthonormal.

According to Lemma 5, for the vector function $\boldsymbol{p}_{\boldsymbol{u}}$, its projection into $\mathbb{S}(t)$ is

$$\mathrm{Pro}_{\mathbb{S}}(\boldsymbol{p}_{\boldsymbol{u}}) = \boldsymbol{K}^{-1}\boldsymbol{u}_{\boldsymbol{p}}(t)\boldsymbol{M}_{\boldsymbol{p}}^{-1}\boldsymbol{r}_{1p} \tag{121}$$

and the coordinate relative to the given base is $\boldsymbol{M}_{\boldsymbol{p}}^{-1}\boldsymbol{r}_{1p}$, where $\boldsymbol{M}_p$ is given by Eq. (32) and $\boldsymbol{r}_{1p}$ is given by Eq. (33). For the column vectors of $\boldsymbol{f}_{\boldsymbol{u}}^{\mathrm{T}}\boldsymbol{\Phi}_o^{\mathrm{T}}(t_f,t)\boldsymbol{g}_{\boldsymbol{x}_f}^{\mathrm{T}}$, their projections into $\mathbb{S}(t)$ are

$$\mathrm{Pro}_{\mathbb{S}}\left(\boldsymbol{f}_{\boldsymbol{u}}^{\mathrm{T}}\boldsymbol{\Phi}_o^{\mathrm{T}}(t_f,t)\boldsymbol{g}_{\boldsymbol{x}_f}^{\mathrm{T}}\right) = \boldsymbol{K}^{-1}\boldsymbol{u}_{\boldsymbol{p}}(t)\boldsymbol{M}_{\boldsymbol{p}}^{-1}\boldsymbol{\Gamma}_{1p} \tag{122}$$

Their coordinates to the given base are the column vectors of $\boldsymbol{M}_{\boldsymbol{p}}^{-1}\boldsymbol{\Gamma}_{1p}$, where $\boldsymbol{\Gamma}_{1p}$ is given by Eq. (34).

Consider the optimality condition (54) for Problem 1 in the view of function space. It means that $\boldsymbol{p}_{\boldsymbol{u}}$ falls into the space spanned by the column vector functions of $\boldsymbol{f}_{\boldsymbol{u}}^{\mathrm{T}}(t)\boldsymbol{\Phi}_o^{\mathrm{T}}(t_f,t)\boldsymbol{g}_{\boldsymbol{x}_f}^{\mathrm{T}}$, with the coefficients given by $\boldsymbol{\pi}$. Through the projection theory, we may reinterpret the optimality condition on parameters given by Eq. (48), that is

**Theorem 3**: *The optimality condition* (48) *is equivalent to the following expression, namely*

$$\mathrm{Pro}_{\mathbb{S}}(\boldsymbol{p}_{\boldsymbol{u}}) + \mathrm{Pro}_{\mathbb{S}}(\boldsymbol{f}_{\boldsymbol{u}}^{\mathrm{T}}(t)\boldsymbol{\Phi}_o^{\mathrm{T}}(t_f,t)\boldsymbol{g}_{\boldsymbol{x}_f}^{\mathrm{T}})\boldsymbol{\pi} = \boldsymbol{0} \tag{123}$$

*which means that the projection of* $\boldsymbol{p}_{\boldsymbol{u}}$ *to* $\mathbb{S}(t)$, *i.e.,* $\mathrm{Pro}_{\mathbb{S}}(\boldsymbol{p}_{\boldsymbol{u}})$, *falls into the space spanned by the column vector functions in* $\mathrm{Pro}_{\mathbb{S}}\left(\boldsymbol{f}_{\boldsymbol{u}}^{\mathrm{T}}(t)\boldsymbol{\Phi}_o^{\mathrm{T}}(t_f,t)\boldsymbol{g}_{\boldsymbol{x}_f}^{\mathrm{T}}\right)$.

*Proof*. From the optimality condition (48), Eq. (123) is easy to be established by multiplying $\boldsymbol{K}^{-1}\boldsymbol{u}_{\boldsymbol{p}}(t)\boldsymbol{M}_{\boldsymbol{p}}^{-1}$. On the other hand, the identity of two elements in a space means that their coordinates are the same. Thus Eq. (123) implies for the given bases, there is

$$\boldsymbol{M}_{\boldsymbol{p}}^{-1}\boldsymbol{r}_{1p} + \boldsymbol{M}_{\boldsymbol{p}}^{-1}\boldsymbol{\Gamma}_{1p}\boldsymbol{\pi} = \boldsymbol{0} \tag{124}$$

which leads to Eq. (48) when removing the invertible matrix $\boldsymbol{M}_{\boldsymbol{p}}$. ■

**Remark 6**: Theorem 3 means that the parameter optimality condition of the transformed NLP is compatible with the optimality condition of the OCP, in that they both show the function dependence at the optimal solutions.

In the optimality condition, the multiplier $\boldsymbol{\pi}$ gives the coordinates for the linear combination. Let us re-investigate Eq. (35), which solves the multiplier.

**Proposition 4**: *Consider the matrixes* $\boldsymbol{M}_1 = \boldsymbol{g}_{\boldsymbol{x}_f}\left(\int_{t_0}^{t_f} \boldsymbol{\Phi}_o(t_f,t)\,\boldsymbol{f}_{\boldsymbol{u}}\boldsymbol{K}\boldsymbol{f}_{\boldsymbol{u}}^{\mathrm{T}}\boldsymbol{\Phi}_o^{\mathrm{T}}(t_f,t)\,\mathrm{d}t\right)\boldsymbol{g}_{\boldsymbol{x}_f}^{\mathrm{T}}$ *and* $\boldsymbol{M}_2 = \boldsymbol{\Gamma}_{1p}^{\mathrm{T}}\boldsymbol{M}_p^{-1}\boldsymbol{\Gamma}_{1p}$, *where* $\boldsymbol{M}_p$ *is given by Eq.* (32) *and* $\boldsymbol{\Gamma}_{1p}$ *is given by Eq.* (34). $\boldsymbol{M}_2$ *is an approximation of* $\boldsymbol{M}_1$ *in the sense of the projection into the vector function space spanned by the column vectors of* $\boldsymbol{K}^{-1}\boldsymbol{u}_{\boldsymbol{p}}(t)$ . *For the vectors* $\boldsymbol{a}_1 = \boldsymbol{g}_{\boldsymbol{x}_f}\int_{t_0}^{t_f} \boldsymbol{\Phi}_o(t_f,t)\,\boldsymbol{f}_{\boldsymbol{u}}\boldsymbol{K}\boldsymbol{p}_{\boldsymbol{u}}\,\mathrm{d}t$ *and* $\boldsymbol{a}_2 = \boldsymbol{\Gamma}_{1p}^{\mathrm{T}}\boldsymbol{M}_p^{-1}\boldsymbol{r}_{1p}$, *where* $\boldsymbol{r}_{1p}$ *is given by Eq.* (33), $\boldsymbol{a}_2$ *is an approximation of* $\boldsymbol{a}_1$ *in the similar sense. In particular,* $\boldsymbol{M}_1 = \boldsymbol{M}_2$ *if i) all the column vectors of* $\boldsymbol{f}_{\boldsymbol{u}}^{\mathrm{T}}\boldsymbol{\Phi}_o^{\mathrm{T}}(t_f,t)\boldsymbol{g}_{\boldsymbol{x}_f}^{\mathrm{T}}$ *are within the function space spanned by the column vectors of* $\boldsymbol{K}^{-1}\boldsymbol{u}_{\boldsymbol{p}}(t)$ ; $\boldsymbol{a}_1 = \boldsymbol{a}_2$ *if i) holds or ii)* $\boldsymbol{p}_{\boldsymbol{u}}$ *is within the function space spanned by the column vectors of* $\boldsymbol{K}^{-1}\boldsymbol{u}_{\boldsymbol{p}}(t)$.

*Proof*. With the projection relations given by Eqs. (121) and (122), there are

$$\boldsymbol{M}_2 = \int_{t_0}^{t_f} \mathrm{Pro}_{\mathbb{S}}(\boldsymbol{g}_{\boldsymbol{x}_f}\boldsymbol{\Phi}_o(t_f,t)\,\boldsymbol{f}_{\boldsymbol{u}})\boldsymbol{K}\mathrm{Pro}_{\mathbb{S}}(\boldsymbol{f}_{\boldsymbol{u}}^{\mathrm{T}}\boldsymbol{\Phi}_o^{\mathrm{T}}(t_f,t)\boldsymbol{g}_{\boldsymbol{x}_f}^{\mathrm{T}})\,\mathrm{d}t = \boldsymbol{\Gamma}_{1p}^{\mathrm{T}}\boldsymbol{M}_p^{-1}\boldsymbol{\Gamma}_{1p} \tag{125}$$

$$\boldsymbol{a}_2 = \int_{t_0}^{t_f} \mathrm{Pro}_{\mathbb{S}}\left(\boldsymbol{f}_{\boldsymbol{u}}^{\mathrm{T}}\boldsymbol{\Phi}_o^{\mathrm{T}}(t_f,t)\boldsymbol{g}_{\boldsymbol{x}_f}^{\mathrm{T}}\right)^{\mathrm{T}}\boldsymbol{K}\mathrm{Pro}_{\mathbb{S}}(\boldsymbol{p}_{\boldsymbol{u}})\,\mathrm{d}t = \boldsymbol{\Gamma}_{1p}^{\mathrm{T}}\boldsymbol{M}_p^{-1}\boldsymbol{r}_{1p} \tag{126}$$

Thus, $\boldsymbol{M}_2$ is an approximation of $\boldsymbol{M}_1$ and $\boldsymbol{a}_2$ is an approximation of $\boldsymbol{a}_1$, in the sense of considering the function components within the space $\mathbb{S}(t)$ only.

If the column vectors of $\boldsymbol{f}_{\boldsymbol{u}}^{\mathrm{T}}\boldsymbol{\Phi}_o^{\mathrm{T}}(t_f,t)\boldsymbol{g}_{\boldsymbol{x}_f}^{\mathrm{T}}$ may be expressed by the linear combination of $\boldsymbol{K}^{-1}\boldsymbol{u}_{\boldsymbol{p}}(t)$ , there is $\mathrm{Pro}_{\mathbb{S}}\left(\boldsymbol{f}_{\boldsymbol{u}}^{\mathrm{T}}\boldsymbol{\Phi}_o^{\mathrm{T}}(t_f,t)\boldsymbol{g}_{\boldsymbol{x}_f}^{\mathrm{T}}\right) = \boldsymbol{f}_{\boldsymbol{u}}^{\mathrm{T}}\boldsymbol{\Phi}_o^{\mathrm{T}}(t_f,t)\boldsymbol{g}_{\boldsymbol{x}_f}^{\mathrm{T}}$. Thus $\boldsymbol{M}_1 = \boldsymbol{M}_2$ are established. Express $\boldsymbol{p}_{\boldsymbol{u}} = \mathrm{Pro}_{\mathbb{S}}(\boldsymbol{p}_{\boldsymbol{u}}) + \mathrm{Pro}_{\mathbb{S}^{\perp}}(\boldsymbol{p}_{\boldsymbol{u}})$, we have

$$\begin{aligned}
\boldsymbol{a}_1 &= \int_{t_0}^{t_f} \mathrm{Pro}_{\mathbb{S}}\left(\boldsymbol{f}_{\boldsymbol{u}}^{\mathrm{T}}\boldsymbol{\Phi}_o^{\mathrm{T}}(t_f,t)\boldsymbol{g}_{\boldsymbol{x}_f}^{\mathrm{T}}\right)^{\mathrm{T}}\boldsymbol{K}\boldsymbol{p}_{\boldsymbol{u}}\,\mathrm{d}t = \int_{t_0}^{t_f} \mathrm{Pro}_{\mathbb{S}}\left(\boldsymbol{f}_{\boldsymbol{u}}^{\mathrm{T}}\boldsymbol{\Phi}_o^{\mathrm{T}}(t_f,t)\boldsymbol{g}_{\boldsymbol{x}_f}^{\mathrm{T}}\right)^{\mathrm{T}}\boldsymbol{K}\left(\mathrm{Pro}_{\mathbb{S}}(\boldsymbol{p}_{\boldsymbol{u}}) + \mathrm{Pro}_{\mathbb{S}^{\perp}}(\boldsymbol{p}_{\boldsymbol{u}})\right)\mathrm{d}t \\
&= \int_{t_0}^{t_f} \mathrm{Pro}_{\mathbb{S}}\left(\boldsymbol{f}_{\boldsymbol{u}}^{\mathrm{T}}\boldsymbol{\Phi}_o^{\mathrm{T}}(t_f,t)\boldsymbol{g}_{\boldsymbol{x}_f}^{\mathrm{T}}\right)^{\mathrm{T}}\boldsymbol{K}\mathrm{Pro}_{\mathbb{S}}(\boldsymbol{p}_{\boldsymbol{u}})\,\mathrm{d}t \\
&= \boldsymbol{a}_2
\end{aligned} \tag{127}$$

Similarly, if $\boldsymbol{p}_{\boldsymbol{u}} = \mathrm{Pro}_{\mathbb{S}}(\boldsymbol{p}_{\boldsymbol{u}})$, there is

$$\begin{aligned}
\boldsymbol{a}_1 &= \int_{t_0}^{t_f} (\boldsymbol{f}_{\boldsymbol{u}}^{\mathrm{T}}\boldsymbol{\Phi}_o^{\mathrm{T}}(t_f,t)\boldsymbol{g}_{\boldsymbol{x}_f}^{\mathrm{T}})^{\mathrm{T}}\boldsymbol{K}\boldsymbol{p}_{\boldsymbol{u}}\,\mathrm{d}t = \int_{t_0}^{t_f}\left(\mathrm{Pro}_{\mathbb{S}}(\boldsymbol{f}_{\boldsymbol{u}}^{\mathrm{T}}\boldsymbol{\Phi}_o^{\mathrm{T}}(t_f,t)\boldsymbol{g}_{\boldsymbol{x}_f}^{\mathrm{T}}) + \mathrm{Pro}_{\mathbb{S}^{\perp}}(\boldsymbol{f}_{\boldsymbol{u}}^{\mathrm{T}}\boldsymbol{\Phi}_o^{\mathrm{T}}(t_f,t)\boldsymbol{g}_{\boldsymbol{x}_f}^{\mathrm{T}})\right)^{\mathrm{T}}\boldsymbol{K}\left(\mathrm{Pro}_{\mathbb{S}}(\boldsymbol{p}_{\boldsymbol{u}})\right)\mathrm{d}t \\
&= \int_{t_0}^{t_f} \mathrm{Pro}_{\mathbb{S}}\left(\boldsymbol{f}_{\boldsymbol{u}}^{\mathrm{T}}\boldsymbol{\Phi}_o^{\mathrm{T}}(t_f,t)\boldsymbol{g}_{\boldsymbol{x}_f}^{\mathrm{T}}\right)^{\mathrm{T}}\boldsymbol{K}\mathrm{Pro}_{\mathbb{S}}(\boldsymbol{p}_{\boldsymbol{u}})\,\mathrm{d}t \\
&= \boldsymbol{a}_2
\end{aligned} \tag{128}$$

Then the proposition is proved. ■

**Remark 7**: According to Proposition 4, it is found that Assumption 2 (and Assumption 4) has the same implication as the controllability requirement in Refs. [22] and [28], which both guarantee the existence for the solution of the multiplier.

Furthermore, we have

**Proposition 5**: *When* $\boldsymbol{p}_{\boldsymbol{u}}$ *is within the function space spanned by the column vectors of* $\boldsymbol{f}_{\boldsymbol{u}}^{\mathrm{T}}\boldsymbol{\Phi}_o^{\mathrm{T}}(t_f,t)\boldsymbol{g}_{\boldsymbol{x}_f}^{\mathrm{T}}$, *if there exists constant vector* $\boldsymbol{y}$ *such that* $\boldsymbol{M}_1\boldsymbol{y} = \boldsymbol{a}_1$, *then there is* $\boldsymbol{M}_2\boldsymbol{y} = \boldsymbol{a}_2$, *where* $\boldsymbol{M}_1$, $\boldsymbol{M}_2$, $\boldsymbol{a}_1$ *and* $\boldsymbol{a}_2$ *are given in Proposition 4.*

*Proof*. When $\boldsymbol{p}_{\boldsymbol{u}}$ is within the function space spanned by the column vectors of $\boldsymbol{f}_{\boldsymbol{u}}{}^{\mathrm{T}}\boldsymbol{\Phi}_o{}^{\mathrm{T}}(t_f,t)\boldsymbol{g}_{\boldsymbol{x}_f}{}^{\mathrm{T}}$, then there exists constant vector $\boldsymbol{y}$ such that

$$\boldsymbol{p}_{\boldsymbol{u}} = \left(\boldsymbol{f}_{\boldsymbol{u}}{}^{\mathrm{T}}\boldsymbol{\Phi}_o{}^{\mathrm{T}}(t_f,t)\boldsymbol{g}_{\boldsymbol{x}_f}{}^{\mathrm{T}}\right)\boldsymbol{y} \tag{129}$$

which is equivalent to

$$\boldsymbol{a}_1 = \boldsymbol{g}_{\boldsymbol{x}_f}\int_{t_0}^{t_f}\boldsymbol{\Phi}_o(t_f,t)\,\boldsymbol{f}_{\boldsymbol{u}}\boldsymbol{K}\boldsymbol{p}_{\boldsymbol{u}}\,\mathrm{d}t = \boldsymbol{g}_{\boldsymbol{x}_f}\left(\int_{t_0}^{t_f}\boldsymbol{\Phi}_o(t_f,t)\,\boldsymbol{f}_{\boldsymbol{u}}\boldsymbol{K}\boldsymbol{f}_{\boldsymbol{u}}{}^{\mathrm{T}}\boldsymbol{\Phi}_o{}^{\mathrm{T}}(t_f,t)\,\mathrm{d}t\right)\boldsymbol{g}_{\boldsymbol{x}_f}{}^{\mathrm{T}}\boldsymbol{y} = \boldsymbol{M}_1\boldsymbol{y} \tag{130}$$

Obviously the projection of Eq. (129) will satisfy

$$\mathrm{Pro}_{\mathbb{S}}(\boldsymbol{p}_{\boldsymbol{u}}) = \mathrm{Pro}_{\mathbb{S}}\left(\boldsymbol{f}_{\boldsymbol{u}}{}^{\mathrm{T}}\boldsymbol{\Phi}_o{}^{\mathrm{T}}(t_f,t)\boldsymbol{g}_{\boldsymbol{x}_f}{}^{\mathrm{T}}\right)\boldsymbol{y} \tag{131}$$

Then

$$\boldsymbol{a}_2 = \int_{t_0}^{t_f}\mathrm{Pro}_{\mathbb{S}}\left(\boldsymbol{f}_{\boldsymbol{u}}{}^{\mathrm{T}}\boldsymbol{\Phi}_o{}^{\mathrm{T}}(t_f,t)\boldsymbol{g}_{\boldsymbol{x}_f}{}^{\mathrm{T}}\right)^{\mathrm{T}}\boldsymbol{K}\mathrm{Pro}_{\mathbb{S}}(\boldsymbol{p}_{\boldsymbol{u}})\,\mathrm{d}t = \left(\int_{t_0}^{t_f}\mathrm{Pro}_{\mathbb{S}}\left(\boldsymbol{f}_{\boldsymbol{u}}{}^{\mathrm{T}}\boldsymbol{\Phi}_o{}^{\mathrm{T}}(t_f,t)\boldsymbol{g}_{\boldsymbol{x}_f}{}^{\mathrm{T}}\right)^{\mathrm{T}}\boldsymbol{K}\mathrm{Pro}_{\mathbb{S}}(\boldsymbol{f}_{\boldsymbol{u}}{}^{\mathrm{T}}\boldsymbol{\Phi}_o{}^{\mathrm{T}}(t_f,t)\boldsymbol{g}_{\boldsymbol{x}_f}{}^{\mathrm{T}})\,\mathrm{d}t\right)\boldsymbol{y} = \boldsymbol{M}_2\boldsymbol{y} \tag{132}$$

Q.E.D. ■

**Remark 8**: Proposition 5 means that the multipliers $\boldsymbol{\pi}$, calculated by Eq. (35) for the parameterization results and by Eq. (40) for the non-parameterization results respectively, are the same when the parameterized control achieves the optimal control, even if the bases given by $\boldsymbol{K}^{-1}\boldsymbol{u}_{\boldsymbol{p}}(t)$ are not complete.

*C. The convergences theorem*

To show the convergence of the costates in utilizing the direct shooting method, the Lipschitz condition of functions, which is important for the existence and uniqueness of the solution for ODEs, will be extended to the functional case. For a vector function $\boldsymbol{f}(\boldsymbol{x},t)$, if the matrix norm for the partial derivative $\boldsymbol{f}_{\boldsymbol{x}}$ is bounded, then we have the Lipschitz property of

$$\|\boldsymbol{f}(\boldsymbol{x}_1,t) - \boldsymbol{f}(\boldsymbol{x}_2,t)\| \le L_{\boldsymbol{x}}\|\boldsymbol{x}_1 - \boldsymbol{x}_2\| \tag{133}$$

where $L$ denotes the Lipschitz constant and the subscript "$\boldsymbol{x}$" indicates the variable. For a vector $\boldsymbol{z} = \begin{bmatrix}\boldsymbol{x}\\ \boldsymbol{y}\end{bmatrix} = \begin{bmatrix}\boldsymbol{x}\\ \boldsymbol{0}\end{bmatrix} + \begin{bmatrix}\boldsymbol{0}\\ \boldsymbol{y}\end{bmatrix}$, through the Minkowski inequality, there is $\|\boldsymbol{z}\| \le \|\boldsymbol{x}\| + \|\boldsymbol{y}\|$. Therefore, given a function of $\boldsymbol{f}(\boldsymbol{x},\boldsymbol{u},t)$ with control input $\boldsymbol{u}$, the corresponding Lipschitz condition may be expressed as

$$\|\boldsymbol{f}(\boldsymbol{x}_1,\boldsymbol{u}_1,t) - \boldsymbol{f}(\boldsymbol{x}_2,\boldsymbol{u}_2,t)\| \le L_{\boldsymbol{x}}\|\boldsymbol{x}_1 - \boldsymbol{x}_2\| + L_{\boldsymbol{u}}\|\boldsymbol{u}_1 - \boldsymbol{u}_2\| \tag{134}$$

Consider this concept for the functional; that is, given a vector functional as follows

$$\boldsymbol{F}(\boldsymbol{x}(t),\boldsymbol{u}(t),t_f) = \int_{t_0}^{t_f}\boldsymbol{f}\left(\boldsymbol{x}(t),\boldsymbol{u}(t),t\right)\mathrm{d}t \tag{135}$$

Its Lipschitz condition may be analogously expressed as

$$\|\boldsymbol{F}(\boldsymbol{x}_1(t),\boldsymbol{u}_1(t),t_{f1}) - \boldsymbol{F}(\boldsymbol{x}_2(t),\boldsymbol{u}_2(t),t_{f2})\| \le L_{\boldsymbol{x}}\|\boldsymbol{x}_1(t) - \boldsymbol{x}_2(t)\|_{\infty} + L_{\boldsymbol{u}}\|\boldsymbol{u}_1(t) - \boldsymbol{u}_2(t)\|_{\infty} + L_{t_f}\left|t_{f1} - t_{f2}\right| \tag{136}$$

Note that here the norm $\|\cdot\|_{\infty}$ denotes the vector function supremum norm.

Moreover, the following lemma that is based on the Gronwall-Bellman inequality will be used in estimating the error development on ODEs.

**Lemma 7** [31]: *Assume that* $\boldsymbol{f}(\boldsymbol{x},\boldsymbol{u},\boldsymbol{v},t)$ *is piecewise continuous in* $t$ *, and Lipschitz in* $\boldsymbol{x}$ *,* $\boldsymbol{u}$ *, and* $\boldsymbol{v}$ *with the Lipschitz constants of* $L_{\boldsymbol{x}}$ *,* $L_{\boldsymbol{u}}$ *, and* $L_{\boldsymbol{v}}$ *. Let* $\boldsymbol{x}_1(t)$ *and* $\boldsymbol{x}_2(t)$ *the solutions of*

$$\dot{\boldsymbol{x}}_1 = \boldsymbol{f}(\boldsymbol{x}_1,\boldsymbol{u}_1,\boldsymbol{v}_1,t),\ \boldsymbol{x}_1(t_0)=\boldsymbol{x}_{10} \tag{137}$$

*and*

$$\dot{\boldsymbol{x}}_2 = \boldsymbol{f}(\boldsymbol{x}_2,\boldsymbol{u}_2,\boldsymbol{v}_2,t),\ \boldsymbol{x}_2(t_0)=\boldsymbol{x}_{20} \tag{138}$$

*with* $\boldsymbol{u}_1$ *,* $\boldsymbol{v}_1$ *and* $\boldsymbol{u}_2$ *,* $\boldsymbol{v}_2$ *the control inputs. Then*

$$\|\boldsymbol{x}_1-\boldsymbol{x}_2\| \le e^{L_x(t-t_0)}\|\boldsymbol{x}_{10}-\boldsymbol{x}_{20}\| + \frac{L_{\boldsymbol{u}}}{L_{\boldsymbol{x}}}\left(e^{L_x(t-t_0)}-1\right)\|\boldsymbol{u}_1-\boldsymbol{u}_2\|_\infty + \frac{L_{\boldsymbol{v}}}{L_{\boldsymbol{x}}}\left(e^{L_x(t-t_0)}-1\right)\|\boldsymbol{v}_1-\boldsymbol{v}_2\|_\infty \tag{139}$$

*where* $\|\cdot\|$ *denotes the vector norm and* $\|\cdot\|_\infty$ *denotes the supremum norm of vector function. Furthermore, it is easy to see that for a given time span* $[t_0,t_f]$ *, there is*

$$\|\boldsymbol{x}_1-\boldsymbol{x}_2\|_\infty \le e^{L_x(t_f-t_0)}\|\boldsymbol{x}_{10}-\boldsymbol{x}_{20}\| + \frac{L_{\boldsymbol{u}}}{L_{\boldsymbol{x}}}\left(e^{L_x(t_f-t_0)}-1\right)\|\boldsymbol{u}_1-\boldsymbol{u}_2\|_\infty + \frac{L_{\boldsymbol{v}}}{L_{\boldsymbol{x}}}\left(e^{L_x(t_f-t_0)}-1\right)\|\boldsymbol{v}_1-\boldsymbol{v}_2\|_\infty \tag{140}$$

**Theorem 4**: *Denote the optimal solution of Problem 1 with a hat of "* $\wedge$ *" and denote the approximate optimal solution, obtained by the direct shooting method, with a tile of "* ~ *". If* $\tilde{\boldsymbol{u}}(t;\boldsymbol{p}) \to \hat{\boldsymbol{u}}$ *and* $\tilde{t}_f \to \hat{t}_f$ *, then we have the states* $\tilde{\boldsymbol{x}} \to \hat{\boldsymbol{x}}$ *and the costates* $\tilde{\boldsymbol{\lambda}} \to \hat{\boldsymbol{\lambda}}$ *, where "* $\to$ *" means "converge to".*

*Proof*. Let $\bar{t}_f = \max(\tilde{t}_f, \hat{t}_f)$ . With Lemma 7, it is easy to obtain that

$$\|\Delta\boldsymbol{x}\|_\infty \le c_{\boldsymbol{u}}^{\boldsymbol{x}}\|\Delta\boldsymbol{u}\|_\infty \tag{141}$$

where $\Delta\boldsymbol{x} = \tilde{\boldsymbol{x}} - \hat{\boldsymbol{x}}$ and $\Delta\boldsymbol{u} = \tilde{\boldsymbol{u}} - \hat{\boldsymbol{u}}$ .

$$c_{\boldsymbol{u}}^{\boldsymbol{x}} = \frac{L_{\boldsymbol{u}}^{\dot{\boldsymbol{x}}}}{L_{\boldsymbol{x}}^{\dot{\boldsymbol{x}}}}\left(e^{L_x(\bar{t}_f-t_0)}-1\right) \tag{142}$$

Here $L_{\boldsymbol{x}}^{\dot{\boldsymbol{x}}}$ and $L_{\boldsymbol{u}}^{\dot{\boldsymbol{x}}}$ are the Lipschitz constants with respect to $\boldsymbol{x}$ and $\boldsymbol{u}$ in Eq. (14).

For the error in the multiplier $\Delta\boldsymbol{\pi} = \tilde{\boldsymbol{\pi}} - \hat{\boldsymbol{\pi}}$ . Denote the matrix and the vector for the optimal solution, defined by Eq. (38) and Eq. (39) respectively, with $\hat{\boldsymbol{M}}_{\boldsymbol{\pi}}$ and $\hat{\boldsymbol{r}}_{\boldsymbol{\pi}}$ . Denote the quantities for the approximate optimal solution, defined by Eq. (36) and Eq. (37) respectively, with $\tilde{\boldsymbol{M}}_{1\boldsymbol{\pi}}$ and $\tilde{\boldsymbol{r}}_{1\boldsymbol{\pi}}$ . Let

$$\boldsymbol{f_u}^{\mathrm{T}}(t)\boldsymbol{\Phi}_o^{\mathrm{T}}(t_f,t)\boldsymbol{g}_{\boldsymbol{x}_f}{}^{\mathrm{T}} = \mathrm{Pro}_{\mathbb{S}}(\boldsymbol{f_u}^{\mathrm{T}}(t)\boldsymbol{\Phi}_o^{\mathrm{T}}(t_f,t)\boldsymbol{g}_{\boldsymbol{x}_f}{}^{\mathrm{T}}) + \mathrm{Pro}_{\mathbb{S}^\perp}(\boldsymbol{f_u}^{\mathrm{T}}(t)\boldsymbol{\Phi}_o^{\mathrm{T}}(t_f,t)\boldsymbol{g}_{\boldsymbol{x}_f}{}^{\mathrm{T}}) \tag{143}$$

$$\boldsymbol{p_u} = \mathrm{Pro}_{\mathbb{S}}(\boldsymbol{p_u}) + \mathrm{Pro}_{\mathbb{S}^\perp}(\boldsymbol{p_u}) \tag{144}$$

Then according to the orthogonality on the function projection (See Proposition 4), we have

$$\tilde{\boldsymbol{M}}_{1\boldsymbol{\pi}} = \hat{\boldsymbol{M}}_{\boldsymbol{\pi}} + \Delta\boldsymbol{M}_{\boldsymbol{\pi}} \tag{145}$$

$$\tilde{\boldsymbol{r}}_{1\boldsymbol{\pi}} = \hat{\boldsymbol{r}}_{\boldsymbol{\pi}} + \Delta\boldsymbol{r}_{\boldsymbol{\pi}} \tag{146}$$

where

$$\Delta\boldsymbol{M}_{\boldsymbol{\pi}} = -\int_{t_0}^{t_f}\mathrm{Pro}_{\mathbb{S}^\perp}(\boldsymbol{g}_{\boldsymbol{x}_f}\boldsymbol{\Phi}_o(t_f,t)\boldsymbol{f_u})\boldsymbol{K}\mathrm{Pro}_{\mathbb{S}^\perp}(\boldsymbol{f_u}^{\mathrm{T}}\boldsymbol{\Phi}_o^{\mathrm{T}}(t_f,t)\boldsymbol{g}_{\boldsymbol{x}_f}{}^{\mathrm{T}})\,\mathrm{d}t \tag{147}$$

$$\Delta\boldsymbol{r}_{\boldsymbol{\pi}} = -\int_{t_0}^{t_f}\mathrm{Pro}_{\mathbb{S}^\perp}(\boldsymbol{g}_{\boldsymbol{x}_f}\boldsymbol{\Phi}_o(t_f,t)\boldsymbol{f_u})\boldsymbol{K}\mathrm{Pro}_{\mathbb{S}^\perp}(\boldsymbol{p_u})\,\mathrm{d}t \tag{148}$$

From the relation of $(\hat{\boldsymbol{M}}_{\boldsymbol{\pi}}+\Delta\boldsymbol{M}_{\boldsymbol{\pi}})(\hat{\boldsymbol{\pi}}+\Delta\boldsymbol{\pi})=(\hat{\boldsymbol{r}}_{\boldsymbol{\pi}}+\Delta\boldsymbol{r}_{\boldsymbol{\pi}})$, we may derive that

$$\Delta\boldsymbol{\pi}=\tilde{\boldsymbol{M}}_{1\boldsymbol{\pi}}{}^{-1}\left(\Delta\boldsymbol{r}_{\boldsymbol{\pi}}-\Delta\boldsymbol{M}_{\boldsymbol{\pi}}\boldsymbol{\pi}\right) \tag{149}$$

Assume the vector functional $\Delta\boldsymbol{r}_{\boldsymbol{\pi}}-\Delta\boldsymbol{M}_{\boldsymbol{\pi}}\hat{\boldsymbol{\pi}}$ satisfies the Lipschitz condition and because the vector functional vanishes at the optimal solution (See Remark 7), then

$$\left\|\Delta\boldsymbol{r}_{\boldsymbol{\pi}}-\Delta\boldsymbol{M}_{\boldsymbol{\pi}}\hat{\boldsymbol{\pi}}\right\|\le L_{\boldsymbol{x}}^{F}\left\|\Delta\boldsymbol{x}\right\|_{\infty}+L_{\boldsymbol{u}}^{F}\left\|\Delta\boldsymbol{u}\right\|_{\infty}+L_{t_f}^{F}\left|\Delta t_f\right| \tag{150}$$

where $\Delta t_f=\tilde{t}_f-\hat{t}_f$. $L_{\boldsymbol{x}}^{F}$, $L_{\boldsymbol{u}}^{F}$, and $L_{t_f}^{F}$ are the Lipschitz constants. Combined with Eq. (141), it may be further obtained that

$$\left\|\Delta\boldsymbol{\pi}\right\|\le c_{\boldsymbol{u}}^{\boldsymbol{\pi}}\left\|\Delta\boldsymbol{u}\right\|_{\infty}+c_{t_f}^{\boldsymbol{\pi}}\left|\Delta t_f\right| \tag{151}$$

where

$$c_{\boldsymbol{u}}^{\boldsymbol{\pi}}=\frac{L_{\boldsymbol{x}}^{F}L_{\boldsymbol{u}}^{\dot{\boldsymbol{x}}}}{L_{\boldsymbol{x}}^{\dot{\boldsymbol{x}}}}\sigma_{\max}(\tilde{\boldsymbol{M}}_{1\boldsymbol{\pi}}{}^{-1})\left(e^{L_{\boldsymbol{x}}^{\dot{\boldsymbol{x}}}(\bar{t}_f-t_0)}-1\right)+\sigma_{\max}(\tilde{\boldsymbol{M}}_{1\boldsymbol{\pi}}{}^{-1})L_{\boldsymbol{u}}^{F} \tag{152}$$

$$c_{t_f}^{\boldsymbol{\pi}}=\sigma_{\max}(\tilde{\boldsymbol{M}}_{1\boldsymbol{\pi}}{}^{-1})L_{t_f}^{F} \tag{153}$$

and $\sigma_{\max}(\cdot)$ denotes the maximum singular value of the matrix.

Now we will consider the costate given by Eq. (117). Introducing an inverse time of $t_{inv}=-t+\bar{t}_f$, then the costate dynamic equation (118) and the terminal boundary condition (119) may be reformulated as

$$\frac{\mathrm{d}}{\mathrm{d}t_{inv}}\boldsymbol{\lambda}(t_{inv})=\boldsymbol{f}_{\boldsymbol{x}}{}^{\mathrm{T}}\boldsymbol{\lambda}(t_{inv})+L_{\boldsymbol{x}},\ \boldsymbol{\lambda}\big|_{t_{inv}=0}=\boldsymbol{\lambda}_0 \tag{154}$$

Similarly, we may obtain the following relation with Lemma 7, that is

$$\left\|\Delta\boldsymbol{\lambda}\right\|_{\infty}\le e^{\sigma_{\max}(\boldsymbol{f}_{\boldsymbol{x}})(\bar{t}_f-t_0)}\left\|\tilde{\boldsymbol{\lambda}}_0-\hat{\boldsymbol{\lambda}}_0\right\|+\frac{L_{\boldsymbol{x}}^{\dot{\boldsymbol{\lambda}}}}{\sigma_{\max}(\boldsymbol{f}_{\boldsymbol{x}})}\left(e^{\sigma_{\max}(\boldsymbol{f}_{\boldsymbol{x}})(\bar{t}_f-t_0)}-1\right)\left\|\Delta\boldsymbol{x}\right\|_{\infty}+\frac{L_{\boldsymbol{u}}^{\dot{\boldsymbol{\lambda}}}}{\sigma_{\max}(\boldsymbol{f}_{\boldsymbol{x}})}\left(e^{\sigma_{\max}(\boldsymbol{f}_{\boldsymbol{x}})(\bar{t}_f-t_0)}-1\right)\left\|\Delta\boldsymbol{u}\right\|_{\infty} \tag{155}$$

where $\Delta\boldsymbol{\lambda}=\tilde{\boldsymbol{\lambda}}-\hat{\boldsymbol{\lambda}}$. $L_{\boldsymbol{x}}^{\dot{\boldsymbol{\lambda}}}$ and $L_{\boldsymbol{u}}^{\dot{\boldsymbol{\lambda}}}$ are the Lipschitz constants with respect to the right function of Eq. (154). Regarding the initial value, there is

$$\begin{aligned}\left\|\tilde{\boldsymbol{\lambda}}_0-\hat{\boldsymbol{\lambda}}_0\right\|&\le\left\|(\varphi_{\boldsymbol{x}}(\tilde{t}_f)+\boldsymbol{g}_{\tilde{\boldsymbol{x}}_f}{}^{\mathrm{T}}\tilde{\boldsymbol{\pi}})-(\varphi_{\boldsymbol{x}}(\hat{t}_f)+\boldsymbol{g}_{\hat{\boldsymbol{x}}_f}{}^{\mathrm{T}}\hat{\boldsymbol{\pi}})\right\|+\left\|\dot{\boldsymbol{\lambda}}_{\mathrm{m}}\Delta t_f\right\|\\&\le\left\|(\varphi_{\boldsymbol{x}}(\tilde{t}_f)+\boldsymbol{g}_{\tilde{\boldsymbol{x}}_f}{}^{\mathrm{T}}\hat{\boldsymbol{\pi}})-(\varphi_{\boldsymbol{x}}(\hat{t}_f)+\boldsymbol{g}_{\hat{\boldsymbol{x}}_f}{}^{\mathrm{T}}\hat{\boldsymbol{\pi}})\right\|+\left\|\boldsymbol{g}_{\tilde{\boldsymbol{x}}_f}{}^{\mathrm{T}}\tilde{\boldsymbol{\pi}}-\boldsymbol{g}_{\tilde{\boldsymbol{x}}_f}{}^{\mathrm{T}}\hat{\boldsymbol{\pi}}\right\|+\left\|\dot{\boldsymbol{\lambda}}_{\mathrm{m}}\right\|\left|\Delta t_f\right|\end{aligned} \tag{156}$$

where $\dot{\boldsymbol{\lambda}}_m$ is the derivative upper limit that satisfies $\dot{\boldsymbol{\lambda}}_{\mathrm{m}i}\ge\left|\dot{\boldsymbol{\lambda}}_i\right|$ $(i=1,2,...,n)$ within the time interval of $\left[\underline{t}_f\ \ \bar{t}_f\right]$ and $\underline{t}_f=\min(\tilde{t}_f,\hat{t}_f)$. Assume that the function $\theta(\boldsymbol{x}_f,t_f)=\varphi_{\boldsymbol{x}}(t_f)+\boldsymbol{g}_{\boldsymbol{x}_f}{}^{\mathrm{T}}\hat{\boldsymbol{\pi}}$ satisfies

$$\left\|\theta(\tilde{\boldsymbol{x}}_f,\tilde{t}_f)-\theta(\hat{\boldsymbol{x}}_f,\hat{t}_f)\right\|\le L_{\boldsymbol{x}}^{\theta}\left\|\Delta\boldsymbol{x}_f\right\|+L_{t_f}^{\theta}\left|\Delta t_f\right| \tag{157}$$

where $\Delta\boldsymbol{x}_f=\tilde{\boldsymbol{x}}_f-\hat{\boldsymbol{x}}_f$. $L_{\boldsymbol{x}}^{\theta}$ and $L_{t_f}^{\theta}$ are the Lipschitz constants. Then we have

$$\begin{aligned}\left\|\tilde{\boldsymbol{\lambda}}_0-\hat{\boldsymbol{\lambda}}_0\right\|&\le L_{\boldsymbol{x}}^{\theta}\left\|\Delta\boldsymbol{x}_f\right\|+(L_{t_f}^{\theta}+\left\|\dot{\boldsymbol{\lambda}}_{\mathrm{m}}\right\|)\left|\Delta t_f\right|+\sigma_{\max}(\boldsymbol{g}_{\tilde{\boldsymbol{x}}_f}{}^{\mathrm{T}})\left\|\Delta\boldsymbol{\pi}\right\|\\&\le L_{\boldsymbol{x}}^{\theta}\left\|\Delta\boldsymbol{x}\right\|_{\infty}+(L_{t_f}^{\theta}+\left\|\dot{\boldsymbol{\lambda}}_{\mathrm{m}}\right\|)\left|\Delta t_f\right|+\sigma_{\max}(\boldsymbol{g}_{\tilde{\boldsymbol{x}}_f}{}^{\mathrm{T}})\left\|\Delta\boldsymbol{\pi}\right\|\end{aligned} \tag{158}$$

Therefore, Eq. (155) may be re-presented as

$$\|\Delta\boldsymbol{\lambda}\|_\infty \le \left( e^{\sigma_{\max}(\boldsymbol{f}_x)(\bar{t}_f - t_0)} L_{\boldsymbol{x}}^{\theta} + \frac{L_{\boldsymbol{x}}^{\dot{\boldsymbol{\lambda}}}}{\sigma_{\max}(\boldsymbol{f}_{\boldsymbol{x}})} (e^{\sigma_{\max}(\boldsymbol{f}_x)(\bar{t}_f - t_0)} - 1) \right) \|\Delta\boldsymbol{x}\|_\infty + e^{\sigma_{\max}(\boldsymbol{f}_x)(\bar{t}_f - t_0)} \sigma_{\max}(\boldsymbol{g}_{\tilde{\boldsymbol{x}}_f}{}^{\mathrm{T}}) \|\Delta\boldsymbol{\pi}\| \\ + \left( \frac{L_{\boldsymbol{u}}^{\dot{\boldsymbol{\lambda}}}}{\sigma_{\max}(\boldsymbol{f}_{\boldsymbol{x}})} \left( e^{\sigma_{\max}(\boldsymbol{f}_x)(\bar{t}_f - t_0)} - 1 \right) \right) \|\Delta\boldsymbol{u}\|_\infty + e^{\sigma_{\max}(\boldsymbol{f}_x)(\bar{t}_f - t_0)} (L_{t_f}^{\theta} + \|\dot{\boldsymbol{\lambda}}_{\mathrm{m}}\|) |\Delta t_f| \tag{159}$$

Upon Eq. (141) and Eq. (151), there is

$$\|\Delta\boldsymbol{\lambda}\|_\infty \le c_{\boldsymbol{u}}^{\boldsymbol{\lambda}} \|\Delta\boldsymbol{u}\|_\infty + c_{t_f}^{\boldsymbol{\lambda}} |\Delta t_f| \tag{160}$$

where

$$c_{\boldsymbol{u}}^{\boldsymbol{\lambda}} = \left( e^{\sigma_{\max}(\boldsymbol{f}_x)(\bar{t}_f - t_0)} L_{\boldsymbol{x}}^{\theta} + \frac{L_{\boldsymbol{x}}^{\dot{\boldsymbol{\lambda}}}}{\sigma_{\max}(\boldsymbol{f}_{\boldsymbol{x}})} (e^{\sigma_{\max}(\boldsymbol{f}_x)(\bar{t}_f - t_0)} - 1) \right) c_{\boldsymbol{u}}^{\boldsymbol{x}} + e^{\sigma_{\max}(\boldsymbol{f}_x)(\bar{t}_f - t_0)} \sigma_{\max}(\boldsymbol{g}_{\tilde{\boldsymbol{x}}_f}{}^{\mathrm{T}}) c_{\boldsymbol{u}}^{\boldsymbol{\pi}} + \left( \frac{L_{\boldsymbol{u}}^{\dot{\boldsymbol{\lambda}}}}{\sigma_{\max}(\boldsymbol{f}_{\boldsymbol{x}})} \left( e^{\sigma_{\max}(\boldsymbol{f}_x)(\bar{t}_f - t_0)} - 1 \right) \right) \tag{161}$$

$$c_{t_f}^{\boldsymbol{\lambda}} = e^{\sigma_{\max}(\boldsymbol{f}_x)(\bar{t}_f - t_0)} \sigma_{\max}(\boldsymbol{g}_{\tilde{\boldsymbol{x}}_f}{}^{\mathrm{T}}) c_{t_f}^{\boldsymbol{\pi}} + e^{\sigma_{\max}(\boldsymbol{f}_x)(\bar{t}_f - t_0)} (L_{t_f}^{\theta} + \|\dot{\boldsymbol{\lambda}}_{\mathrm{m}}\|) \tag{162}$$

Therefore, according to Eqs. (141), (151) and (160), the error to the optimal solution, on the states, the multipliers and the costates, will vanish when the control and the terminal time approaches the optimal. ■

**Remark 9:** Denote the optimal performance index with $\hat{J}$ and denote the performance index on the approximate optimal solution with $\tilde{J}$; upon Theorem 3, if $\tilde{\boldsymbol{u}}(t;\boldsymbol{p}) \to \hat{\boldsymbol{u}}$ and $\tilde{t}_f \to \hat{t}_f$, then it is obvious that $\tilde{J} \to \hat{J}$.

Combined with Proposition 1, it may be claimed that:

**Theorem 5**: *Given a set of parameterization as $\boldsymbol{u}^{(k)}(t) = \boldsymbol{u}(t;\boldsymbol{p}^{(k)})$, $k = 1,2,...,k_u$, where $k_u$ may be infinite or finite. The representation capacity of $\boldsymbol{u}^{(k)}(t)$ increases as $k$ increases and $\boldsymbol{u}^{(k_u)}(t)$ completely characterizes the optimal control. Assuming that there is only one optimal approximation to an target optimal control solution in the parameterized function space $\boldsymbol{u}(t) = \boldsymbol{u}(t;\boldsymbol{p}^{(k)}), \boldsymbol{p} \in \mathbb{R}^{s^{(k)}}$, then $\tilde{\boldsymbol{u}}(t;\boldsymbol{p}) \to \hat{\boldsymbol{u}}$, $\tilde{t}_f \to \hat{t}_f$, $\tilde{\boldsymbol{x}} \to \hat{\boldsymbol{x}}$, $\tilde{\boldsymbol{\lambda}} \to \hat{\boldsymbol{\lambda}}$ and $\tilde{J} \to \hat{J}$.*

## VI. Illustrative Examples

First a linear example taken from Xie [36] is solved.

**Example 1**: Consider the following dynamic system

$$\dot{\boldsymbol{x}} = \boldsymbol{A}\boldsymbol{x} + \boldsymbol{b}u$$

where $\boldsymbol{x} = \begin{bmatrix} x_1 \\ x_2 \end{bmatrix}$, $\boldsymbol{A} = \begin{bmatrix} 0 & 1 \\ 0 & 0 \end{bmatrix}$, and $\boldsymbol{b} = \begin{bmatrix} 0 \\ 1 \end{bmatrix}$. Find the solution that minimizes the performance index

$$J = \frac{1}{2}\int_{t_0}^{t_f} u^2 \mathrm{d}t$$

with the boundary conditions

$$\boldsymbol{x}(t_0) = \begin{bmatrix} 1 \\ 1 \end{bmatrix}, \ \boldsymbol{x}(t_f) = \begin{bmatrix} 0 \\ 0 \end{bmatrix}$$

where the initial time $t_0 = 0$ and the terminal time $t_f = 2$ are fixed.

In Ref. [29], this example was solved with the EPDE. Here we will solve it with the direct shooting method, in the following parameterization.

$$u = p_1 + p_2 t + p_3 t^2 + p_4 t^3$$

The concrete form of the parameterization evolution equations for the parameters $\boldsymbol{p} = \begin{bmatrix} p_1 & p_2 & p_3 & p_4 \end{bmatrix}^{\mathrm{T}}$ is

$$\frac{\mathrm{d}\boldsymbol{p}}{\mathrm{d}\tau} = -\boldsymbol{M}_{\boldsymbol{p}}^{-1} \left( \begin{bmatrix} 2 & 2 & 8/3 & 4 \\ 2 & 8/3 & 4 & 6.4 \\ 8/3 & 4 & 6.4 & 32/3 \\ 4 & 64 & 32/3 & 128/7 \end{bmatrix} \begin{bmatrix} p_1 \\ p_2 \\ p_3 \\ p_4 \end{bmatrix} + \begin{bmatrix} 2.0000 & 2.0000 \\ 1.3325 & 2.0000 \\ 1.3325 & 2.6675 \\ 1.5983 & 4.0025 \end{bmatrix} \boldsymbol{\pi} \right)$$

The one-dimensional gain matrix $K$ was set as $K = 0.1$, and then the positive-definite constant matrix $\boldsymbol{M}_{\boldsymbol{p}}$ and its inverse were

$$\boldsymbol{M}_{\boldsymbol{p}} = \begin{bmatrix} 20 & 20 & 80/3 & 40 \\ 20 & 80/3 & 40 & 64 \\ 80/3 & 40 & 64 & 320/3 \\ 40 & 64 & 320/3 & 1280/7 \end{bmatrix}, \ \boldsymbol{M}_{\boldsymbol{p}}^{-1} = \begin{bmatrix} 0.8 & -3 & 3 & -0.875 \\ -3 & 15 & -16.875 & 5.25 \\ 3 & -16.875 & 20.25 & -6.5625 \\ -0.875 & 5.25 & -6.5625 & 2.1875 \end{bmatrix}$$

The Lagrange multiplier $\boldsymbol{\pi}$ was

$$\boldsymbol{\pi} = -\begin{bmatrix} 0.2668 & 0.2000 \\ 0.2000 & 0.2000 \end{bmatrix}^{-1} \left( \begin{bmatrix} 2.0000 & 2.0000 \\ 1.3325 & 2.0000 \\ 1.3325 & 2.6675 \\ 1.5983 & 4.0025 \end{bmatrix}^{\mathrm{T}} \begin{bmatrix} p_1 \\ p_2 \\ p_3 \\ p_4 \end{bmatrix} - \boldsymbol{K}_g \boldsymbol{x}(t_f) \right)$$

with $\boldsymbol{K}_g = \begin{bmatrix} 0.1 & 0 \\ 0 & 0.1 \end{bmatrix}$. The state variables were computed by

$$\boldsymbol{x}(t,\tau) = e^{\boldsymbol{A}(t-t_0)} \begin{bmatrix} 1 \\ 1 \end{bmatrix} + \int_{t_0}^{t} e^{\boldsymbol{A}(t-s)} \boldsymbol{b} u(s,\tau;\boldsymbol{p}) \,\mathrm{d}s$$

and the costates indirectly determined was

$$\boldsymbol{\lambda}(t,\tau) = \left( e^{\boldsymbol{A}(t_f - t)} \right)^{\mathrm{T}} \boldsymbol{\pi}$$

The initial conditions were $\boldsymbol{p}(\tau)\big|_{\tau=0} = \boldsymbol{0}$, and thus a 4 dimensional IVP is produced, much less quantities to be determined than addressing the EPDE. To solve the IVP, the ODE integrator "ode45" in Matlab, with default relative error tolerance $1\times10^{-3}$ and default absolute error tolerance $1\times10^{-6}$, was employed. For comparison, the analytic solution of this example is presented as follows.

$$\begin{aligned} \hat{x}_1 &= 0.5t^3 - 1.75t^2 + t + 1 \\ \hat{x}_2 &= 1.5t^2 - 3.5t + 1 \\ \hat{\lambda}_1 &= 3 \\ \hat{\lambda}_2 &= -3t + 3.5 \\ \hat{u} &= 3t - 3.5 \\ \hat{\boldsymbol{\pi}} &= \begin{bmatrix} 3 & -2.5 \end{bmatrix}^{\mathrm{T}} \end{aligned}$$

It is shown that the optimal control is linear to the time, which is within the representation capacity of the parameterized control. Therefore, it is expected that the parameterized solution will achieve the optimal exactly, and this will be demonstrated in the following.

Figure 5 gives the dynamic curve for the parameters against the virtual time $\tau$ , and they converge to their equilibrium value . At $\tau$ =300s, the numerical solution is $\boldsymbol{p}(\tau)\big|_{\tau=300} = \begin{bmatrix}-3.5000 & 3.0000 & 0.0000 & 0.0000\end{bmatrix}^{\mathrm{T}}$ . Compared with the coefficients in $\hat{u}$ , it is found the optimal control solution is solved precisely. Figure 6 shows the evolution process of the corresponding parameterized control $u(t)$ towards the analytic. The numerical solution coincides with the analytic solution at $\tau$ =300s, as expected. In Fig. 7, the numerical solutions of the states during the computation are plotted. It is shown that they converge to the optimal solutions, and the infeasibilities in the terminal conditions are gradually removed. Figures 8 and 9 plot the evolution profiles of the implicit quantities in the direct shooting method, including the Lagrange multipliers and the costates. Consistent with the theoretical analysis, it is shown that they all tend to the optimal results. The Lagrange multipliers start with the value of $\boldsymbol{\pi}\big|_{\tau=0} = \begin{bmatrix} 2.9924 \\ -2.4924 \end{bmatrix}$ at $\tau$ =0s, and at $\tau$ =300s, the Lagrange multipliers are $\boldsymbol{\pi}\big|_{\tau=300} = \begin{bmatrix} 3.0000 \\ -2.5000 \end{bmatrix}$, which equal the analytic results exactly. The costate $\lambda_2$ , at the first glance, seems unchanged during the evolution. However, when investigating the close-ups, it is found that the analytic optimal solution is achieved from different solutions, even if their error is small.

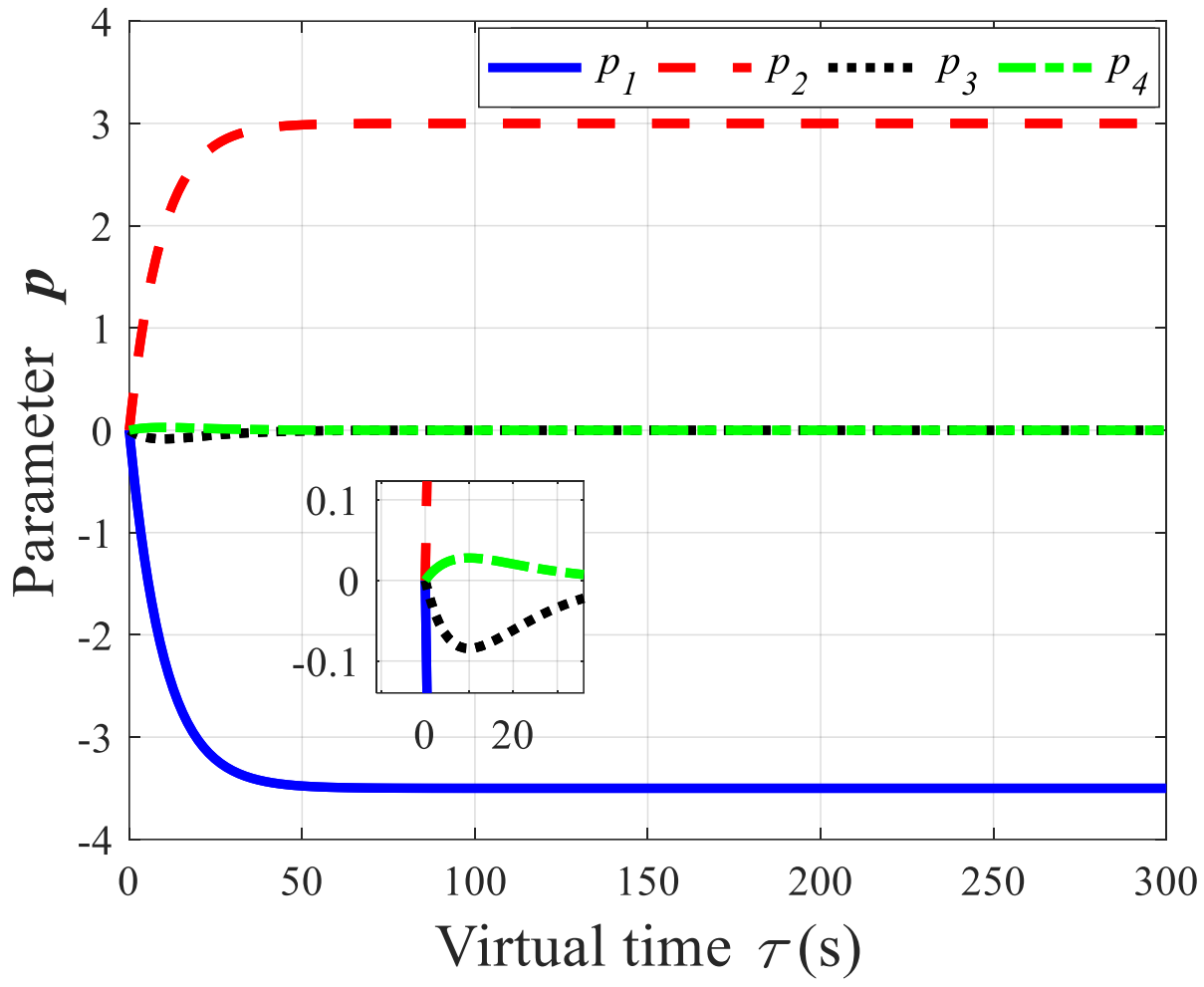

Fig. 5 The evolution profiles of the control parameters.

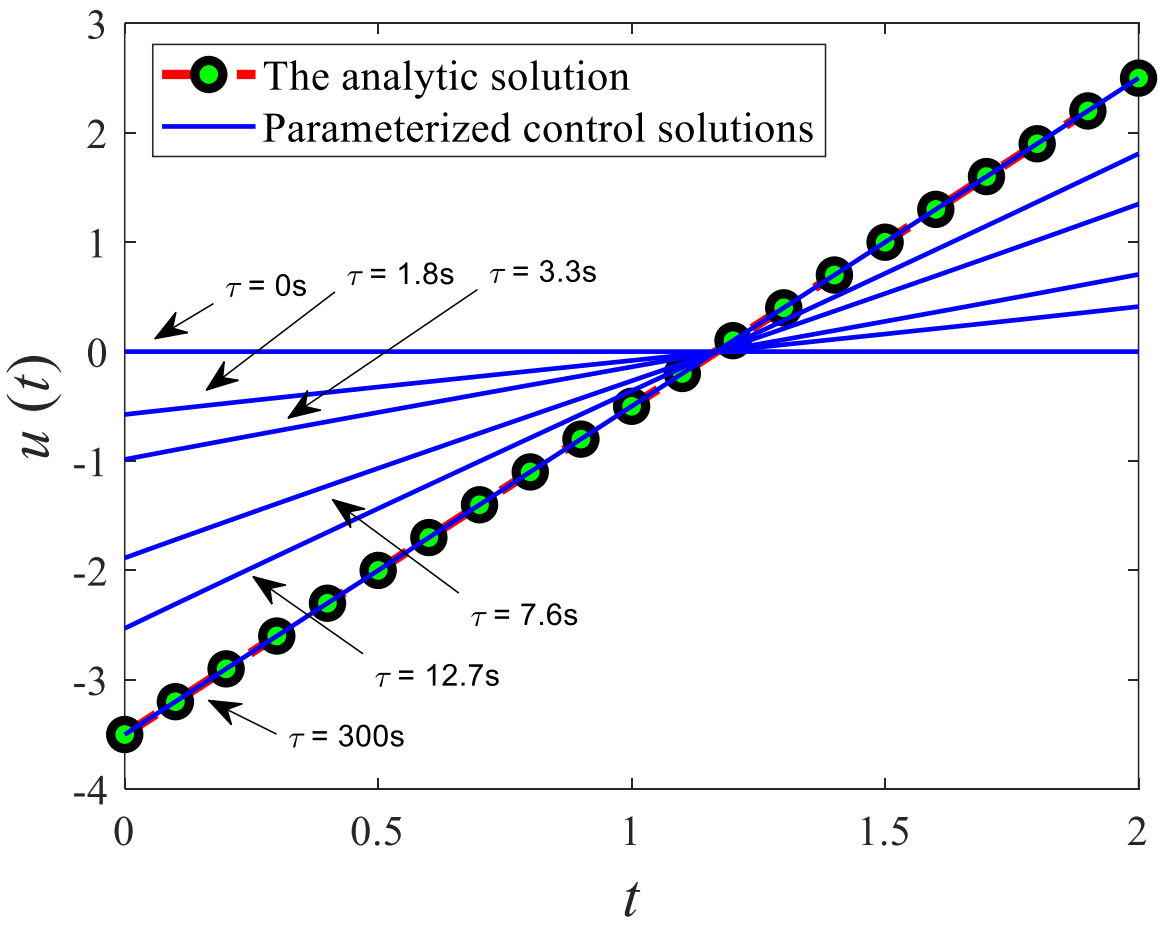

Fig. 6 The evolution of parameterized control solution $u$ to the analytic solution.

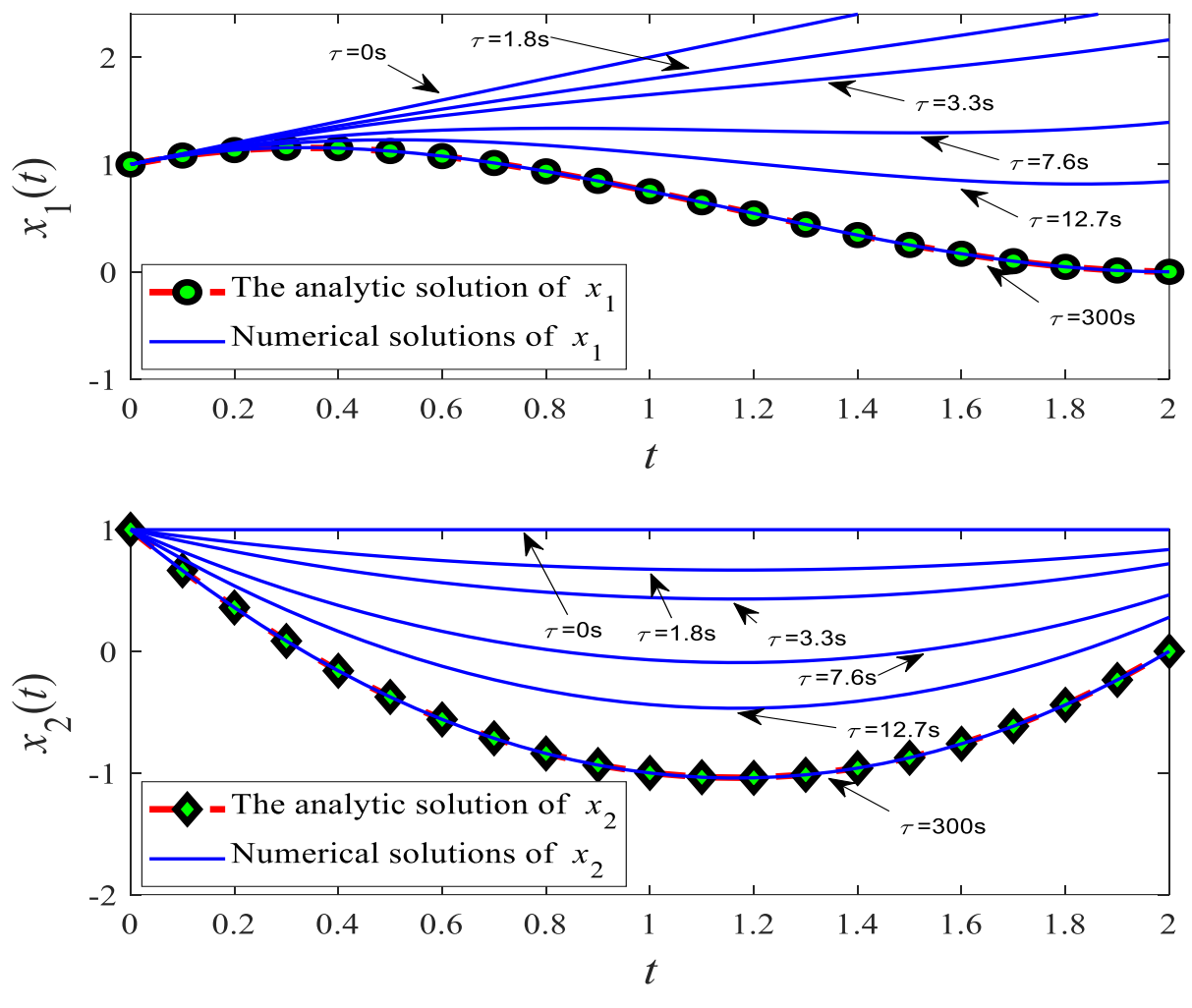

Fig. 7 The evolution of state solutions $x_1$ and $x_2$ to the analytic solution.

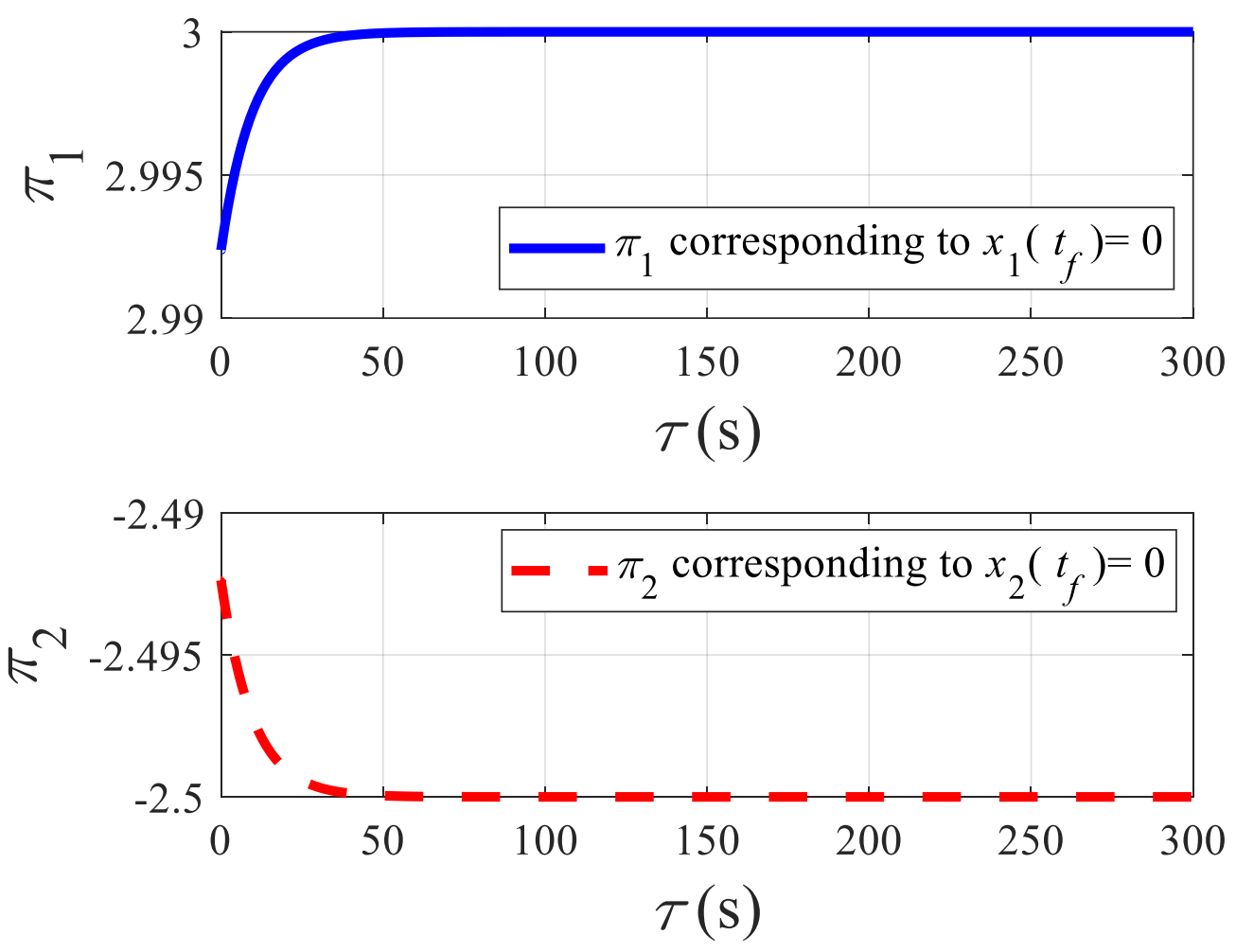

Fig. 8 The evolution profiles of Lagrange multipliers.

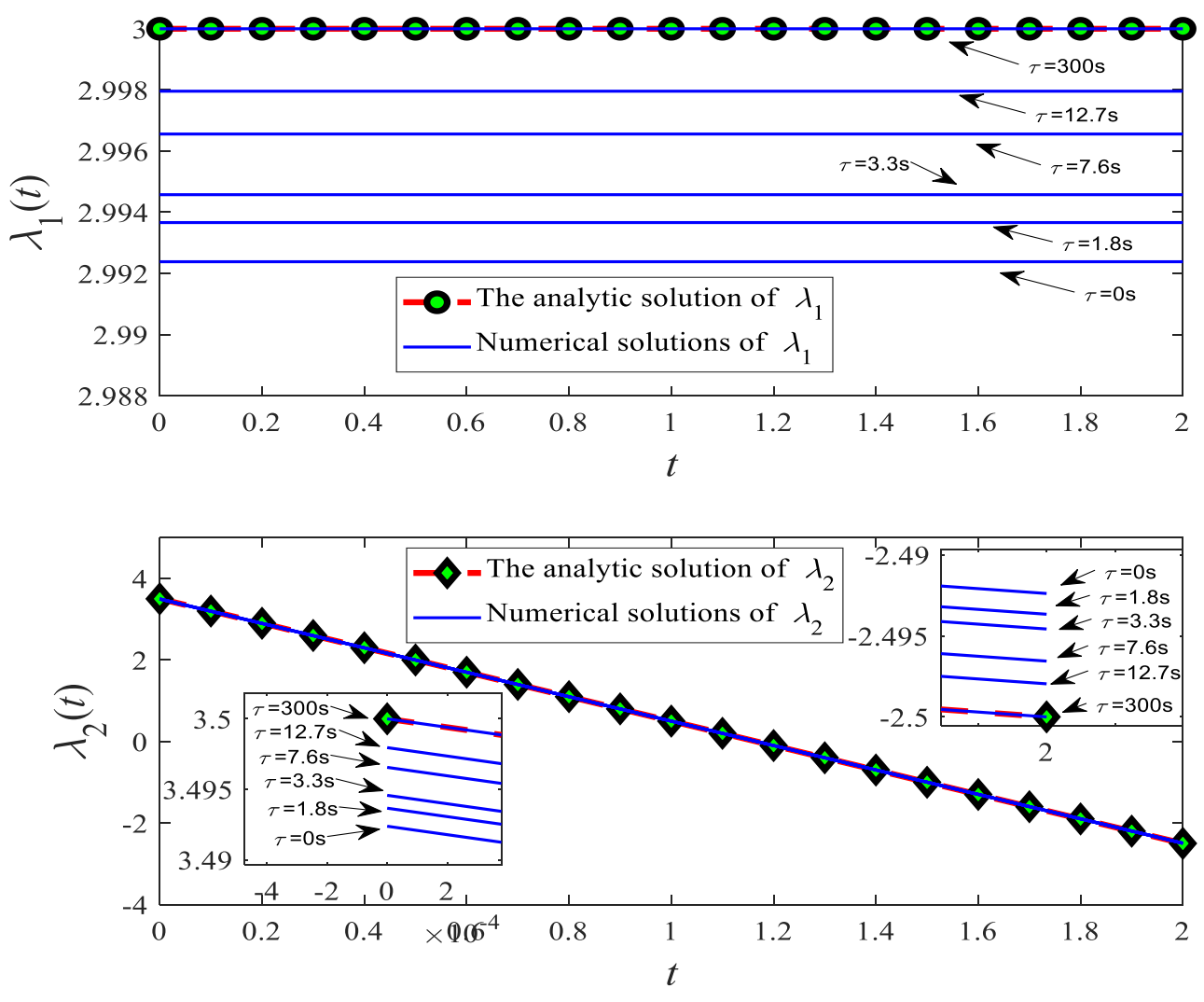

Fig. 9 The evolution of the costates $\lambda_1$ and $\lambda_2$ to the analytic solution.

In particular, we investigated solving the OCP with an arbitrary positive-definite matrix $\boldsymbol{M}_p$ in the parameterization evolution equation. In Sec. IV.C, the analysis shows that this gain matrix may be set more flexibly. From this view, we directly set it to be $\boldsymbol{M}_p = \begin{bmatrix} 0.1 & 0 & 0 & 0 \\ 0 & 0.1 & 0 & 0 \\ 0 & 0 & 0.1 & 0 \\ 0 & 0 & 0 & 0.1 \end{bmatrix}$ and then solve the IVP. Figure 10 gives the profiles of the parameters. Same to the previous, they also converge to the value of $\begin{bmatrix} -3.5000 & 3.0000 & 0.0000 & 0.0000 \end{bmatrix}^{\mathrm{T}}$, demonstrating the effectiveness in solving the OCP with an arbitrary constant positive-definite gain matrix. Therefore, the control profiles and the states all converge to the optimal solutions, in the manner similar to those illustrated in Figs. 6-7. Regarding the costates, figure 11 plots their evolution profiles. It is shown that they reach the optimal solution, which also implies the convergence of the multipliers. Different from Fig. 9, here the initial solutions are quite different and the convergence to optimal solution is more distinct.

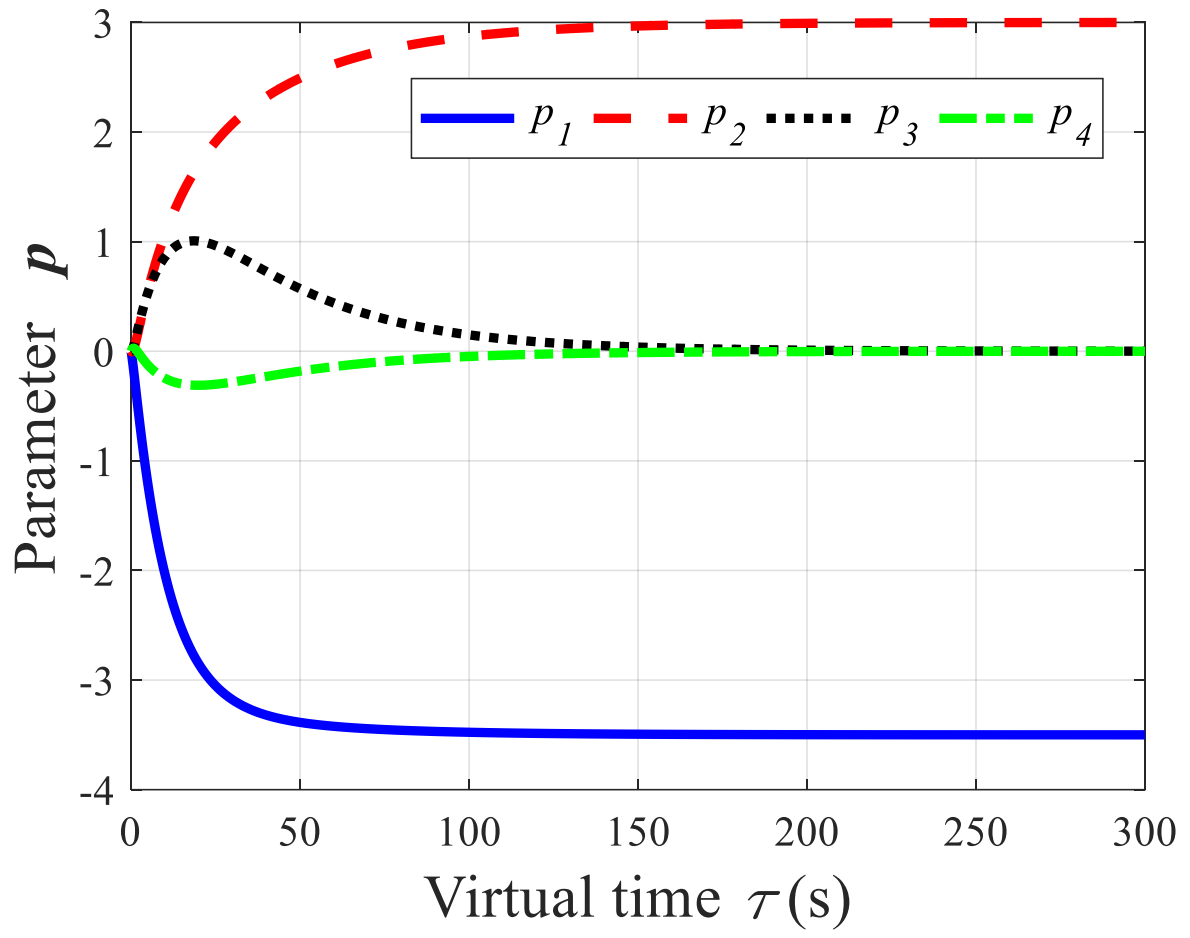

Fig. 10 The evolution profiles of the parameters with $\boldsymbol{M}_p$ a diagonal positive-definite matrix.

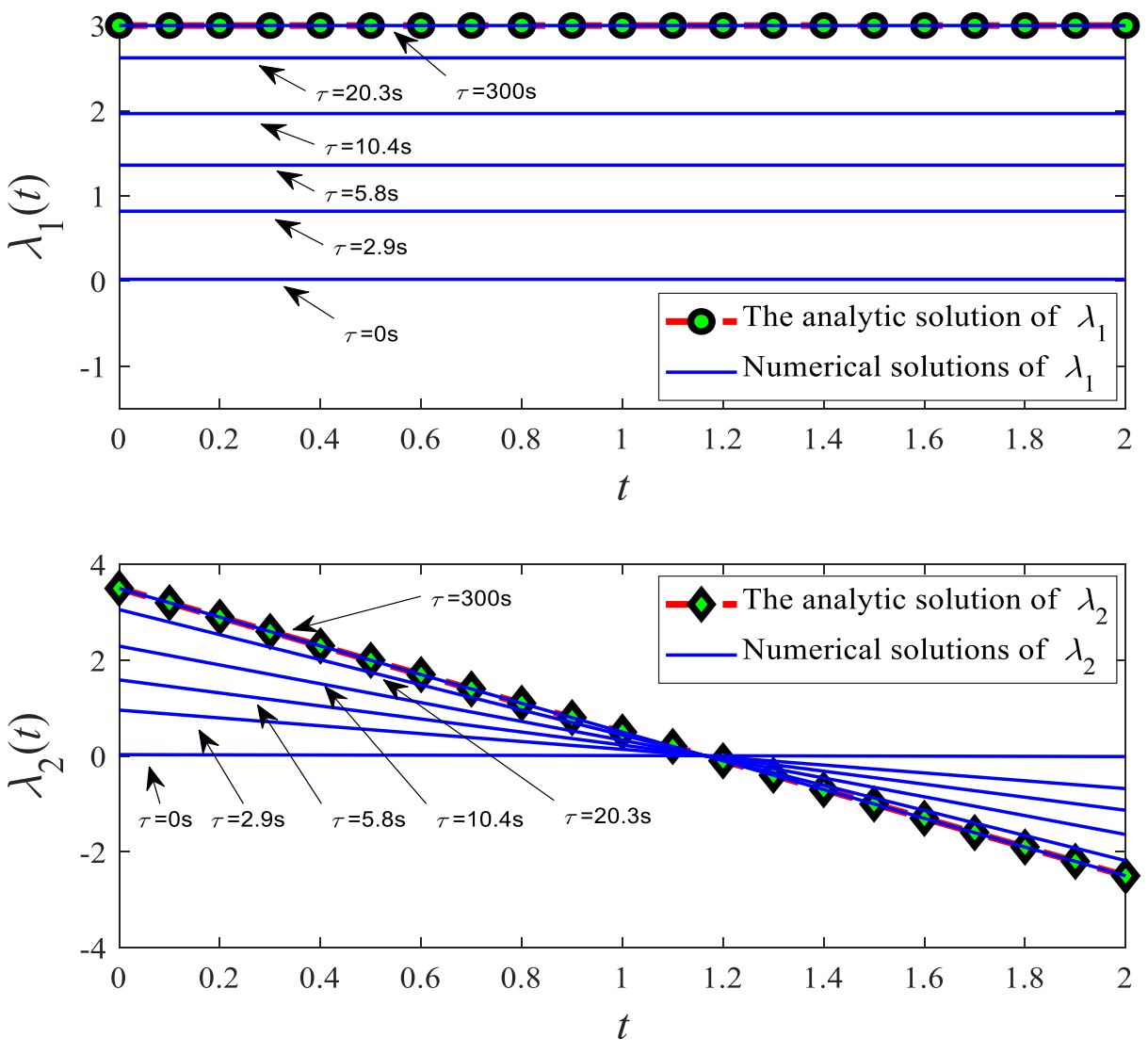

Fig. 11 The evolution profiles of the costates with $\boldsymbol{M}_p$ a diagonal positive-definite matrix.

Compared with solving the infinite-dimensional EPDE, one advantage with the direct shooting method is that the number of parameters may be greatly reduced during the computation. Thus the efficiency may be improved. In Ref. [29], with the semi-discrete method, the time horizon $[t_0, t_f]$ was discretized uniformly with 41 points. The resulting dynamic system has 41 states and it takes 1.31s to get the solution at $\tau$ =300s, slightly longer than the parameterized approach, which takes 0.92s.

Now we consider a nonlinear example with free terminal time $t_f$, the Brachistochrone problem [37], which describes the motion curve of the fastest descending.

**Example 2**: Consider the following dynamic system

$$\dot{\boldsymbol{x}} = \boldsymbol{f}(\boldsymbol{x}, u)$$

where $\boldsymbol{x} = \begin{bmatrix} x \\ y \\ V \end{bmatrix}$, $\boldsymbol{f} = \begin{bmatrix} V\sin(u) \\ -V\cos(u) \\ g\cos(u) \end{bmatrix}$, and $g = 10$ is the gravity constant. Find the solution that minimizes the performance index

$$J = t_f$$

with the boundary conditions

$$\begin{bmatrix} x \\ y \\ V \end{bmatrix}\Bigg|_{t_0=0} = \begin{bmatrix} 0 \\ 0 \\ 0 \end{bmatrix}, \begin{bmatrix} x \\ y \end{bmatrix}\Bigg|_{t_f} = \begin{bmatrix} 2 \\ -2 \end{bmatrix}$$

This example was solved with the EPDE in Ref. [29] before. Here it will be addressed via the direct shooting method. To solve this problem, we consider four cases of parameterization as follows.

1) Parameterization with the global polynomials

$$u = p_1 + p_2 t + p_3 t^2 + p_4 t^3 + p_5 t^4$$

where $\boldsymbol{p} = \begin{bmatrix} p_1 & p_2 & p_3 & p_4 & p_5 \end{bmatrix}^{\mathrm{T}}$ represents the coefficients of polynomials.

2) Parameterization with the Lagrange interpolation polynomials

$$\boldsymbol{u}(t) = \sum_{i=1}^{5} p_i \prod_{j=1, j\neq i}^{5} \frac{(t - t_j)}{(t_i - t_j)}$$

where $\boldsymbol{p} = \begin{bmatrix} p_1 & p_2 & p_3 & p_4 & p_5 \end{bmatrix}^{\mathrm{T}}$ represents the control value at time points of $t_i$, equally partitioned from $[t_0, t_f]$.

3) Parameterization with the linear interpolation between control values at equally partitioned time points, as Eq. (85) gives. The total number of time points is $N = 20$.

4) Parameterization with the step approximation given by Eq. (88), which implements the discrete control. The total number of control parameters is also $N = 20$.

The first and the second cases belong to the global parameterization, and they have the same representation capacity. The third and the forth parameterization employ a local type. For the resulting parameterization evolution equations of the four cases, the gain parameters were all set as $K = 0.1$, $k_{t_f} = 0.1$ and $\boldsymbol{K}_g = \begin{bmatrix} 0.1 & 0 \\ 0 & 0.1 \end{bmatrix}$. The initial conditions were all set to be $\boldsymbol{p}(\tau)\big|_{\tau=0} = \boldsymbol{0}$ and $t_f(\tau)\big|_{\tau=0}$=1s. We again employed "ode45" in Matlab to carry out the numerical integration for the resulting IVPs. In the integrator setting, the integration time horizon was 300s. The relative error tolerance and the absolute error tolerance were $1\times10^{-3}$ and $1\times10^{-6}$, respectively. For comparison, we computed the optimal solution with GPOPS-II [38], a Radau PS method based OCP solver.

The evolution profiles of the parameter solutions are plotted in Fig. 12, and the asymptotical approach to their equilibrium values are demonstrated. In Fig. 13, the terminal time profile against the virtual time $\tau$ is plotted. For the four cases, the results of $t_f$ decline rapidly at first and then approach the equilibrium values gradually. They are almost unchanged after $\tau$ = 50s. Table 1 presents the results of $t_f$ at $\tau$ = 300s. It is shown that the results from the four cases of parameterization are all same to the minimum decline time from GPOPS-II. However, when examining the close-up in Fig. 13, it is found that the results of Case 1-3 are indistinguishable from the optimal solution, while the result of Case 4 is slightly different from others. In particular, the parameter results at $\tau$ =300s for the first case are $\boldsymbol{p}(\tau)\big|_{\tau=300} = \begin{bmatrix}0.0000 & 1.4771 & 0.0000 & 0.0000 & 0.0000\end{bmatrix}^{\mathrm{T}}$, and for the second case are $\boldsymbol{p}(\tau)\big|_{\tau=300} = \begin{bmatrix}0.0000 & 0.3015 & 0.6030 & 0.9045 & 1.2060\end{bmatrix}^{\mathrm{T}}$. Substituting these results to their parameterization control expressions, it is found that the results are identical, which supports the view that the Form 1 and Form 2 parameterization are equivalent when they have the same representation capacity.

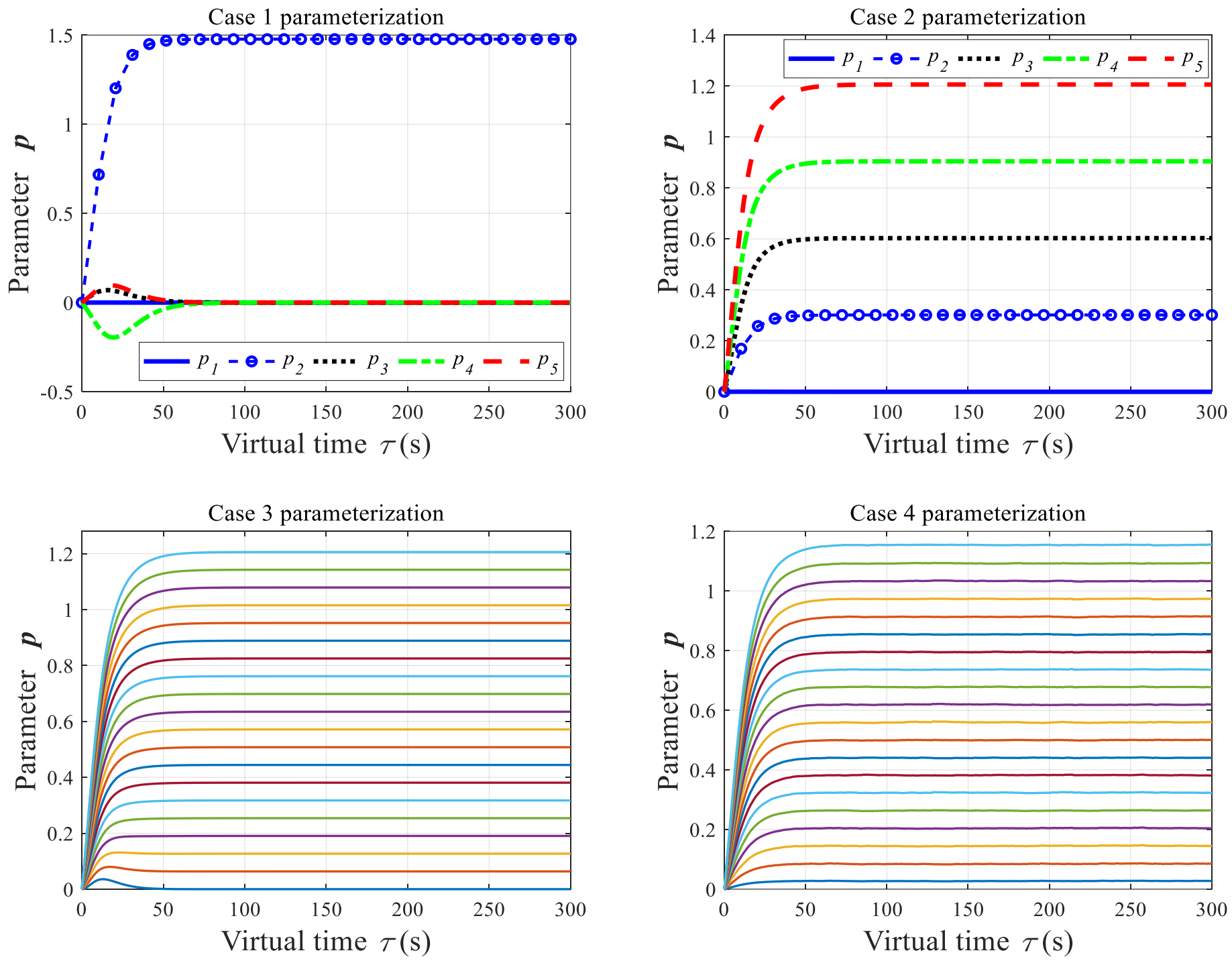

Fig. 12 The evolution profiles of the control parameters for the four parameterization cases.

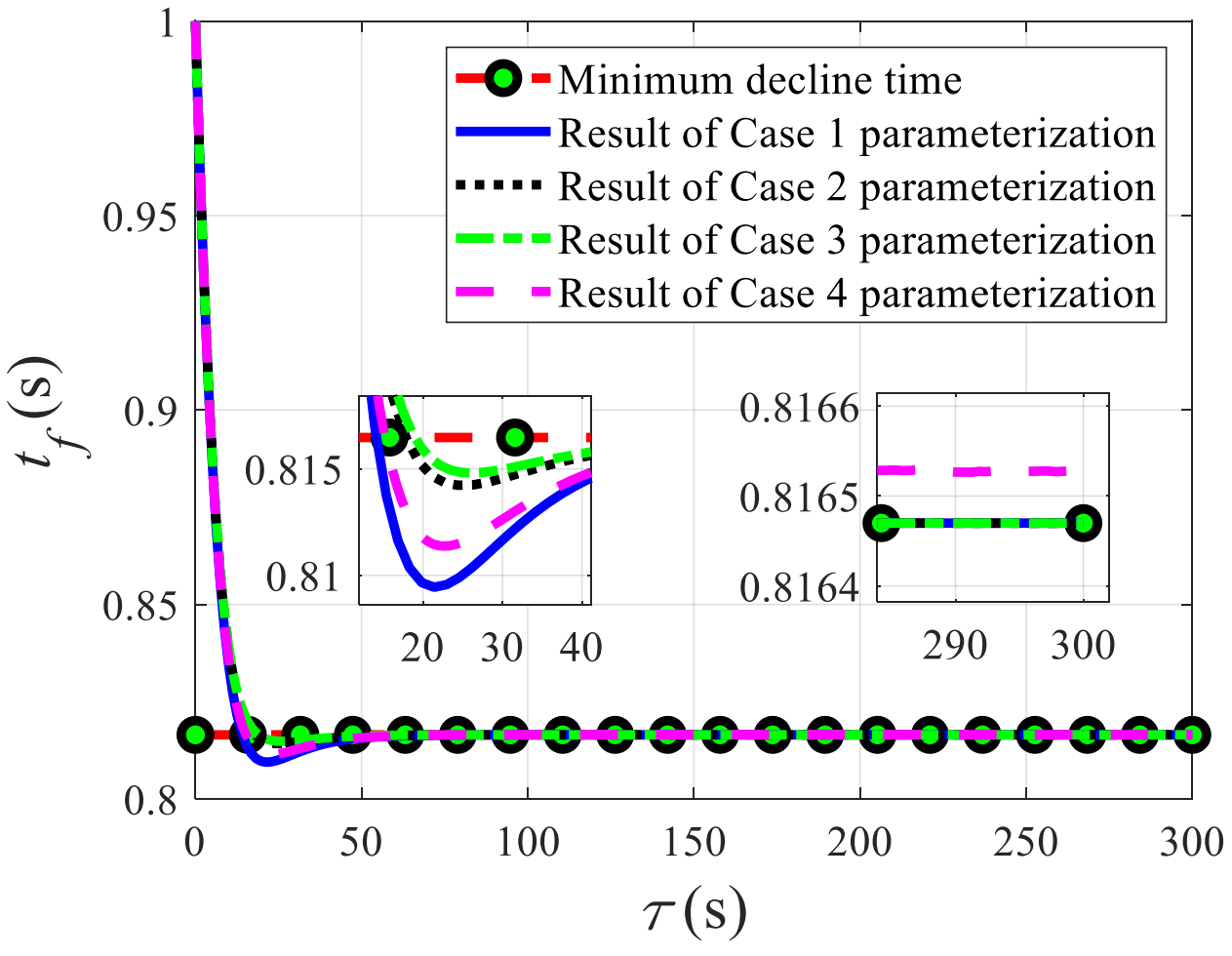

Fig. 13 The evolution profiles of $t_f$ to the minimum decline time.

Table 1 Comparison of results at $\tau = 300\text{s}$

| | Case 1 parameterization | Case 2 parameterization | Case 3 parameterization | Case 4 parameterization | Optimal solution from GPOPS-II |
|---|---|---|---|---|---|
| $t_f$ at $\tau$ = 300s | 0.8165 | 0.8165 | 0.8165 | 0.8165 | 0.8165 |
| $\boldsymbol{\pi}$ at $\tau$ = 300s | $[-0.1477\ \ 0.0564]^{T}$ | $[-0.1477\ \ 0.0564]^{T}$ | $[-0.1477\ \ 0.0564]^{T}$ | $[-0.1469\ \ 0.0589]^{T}$ | not given |

The parameterized controls, based on the discrete parameters, are presented in Fig. 14. Each subplot depicts the evolution of the parameterized control solution for one specific parameterization. It is shown that they all approach the optimal solution well. Since Cases 1-3 have very similar solutions, in the following only the solutions from Case 1 and Case 4 are compared. Figure 15 gives the evolution curves of states in the $xy$ coordinate plane for Case 1 and 4, showing that the numerical results starting from the vertical line approach the optimal solution over time.

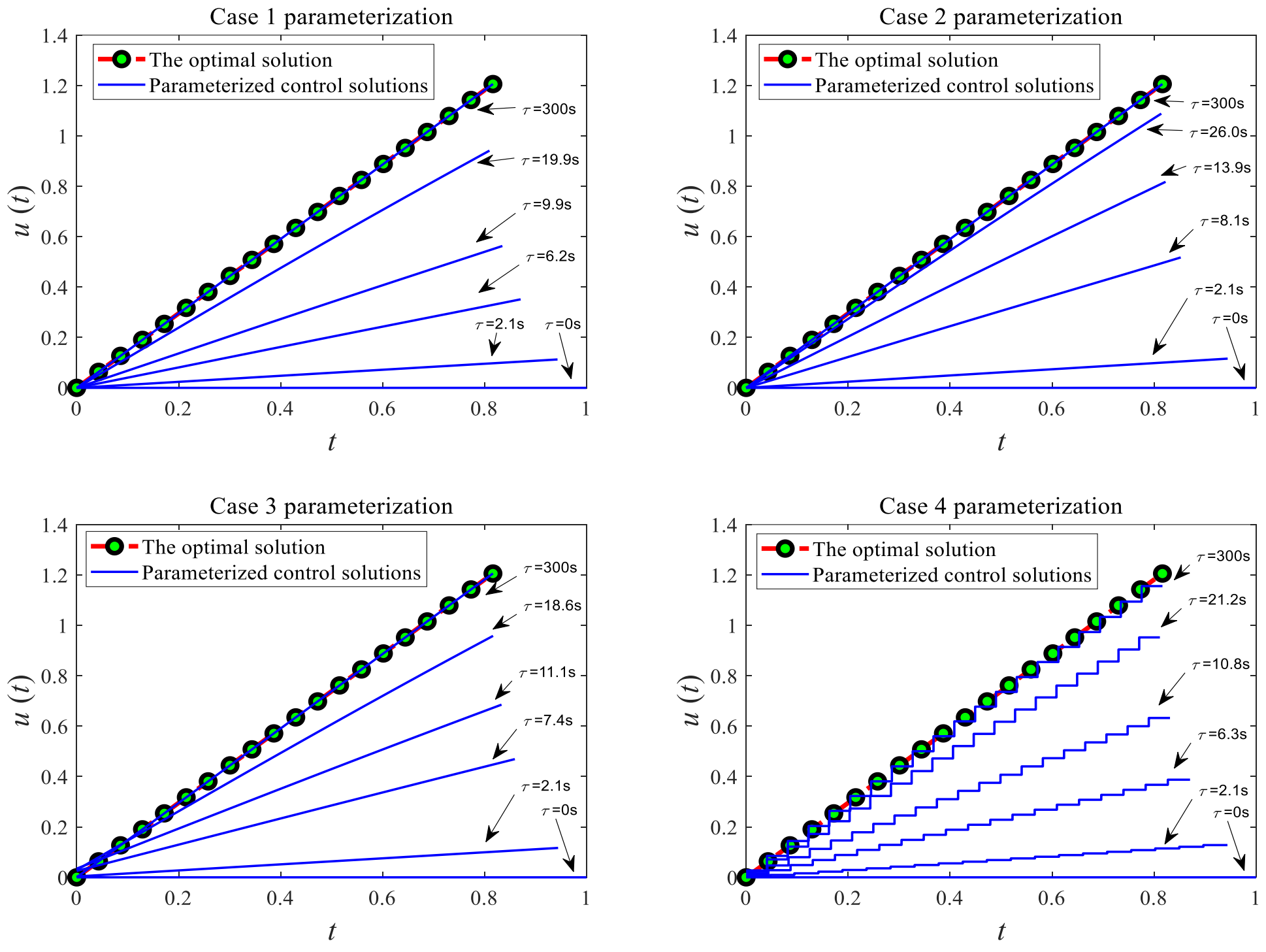

Fig. 14 The evolution of parameterized control solution $u$ for the four parameterization cases.

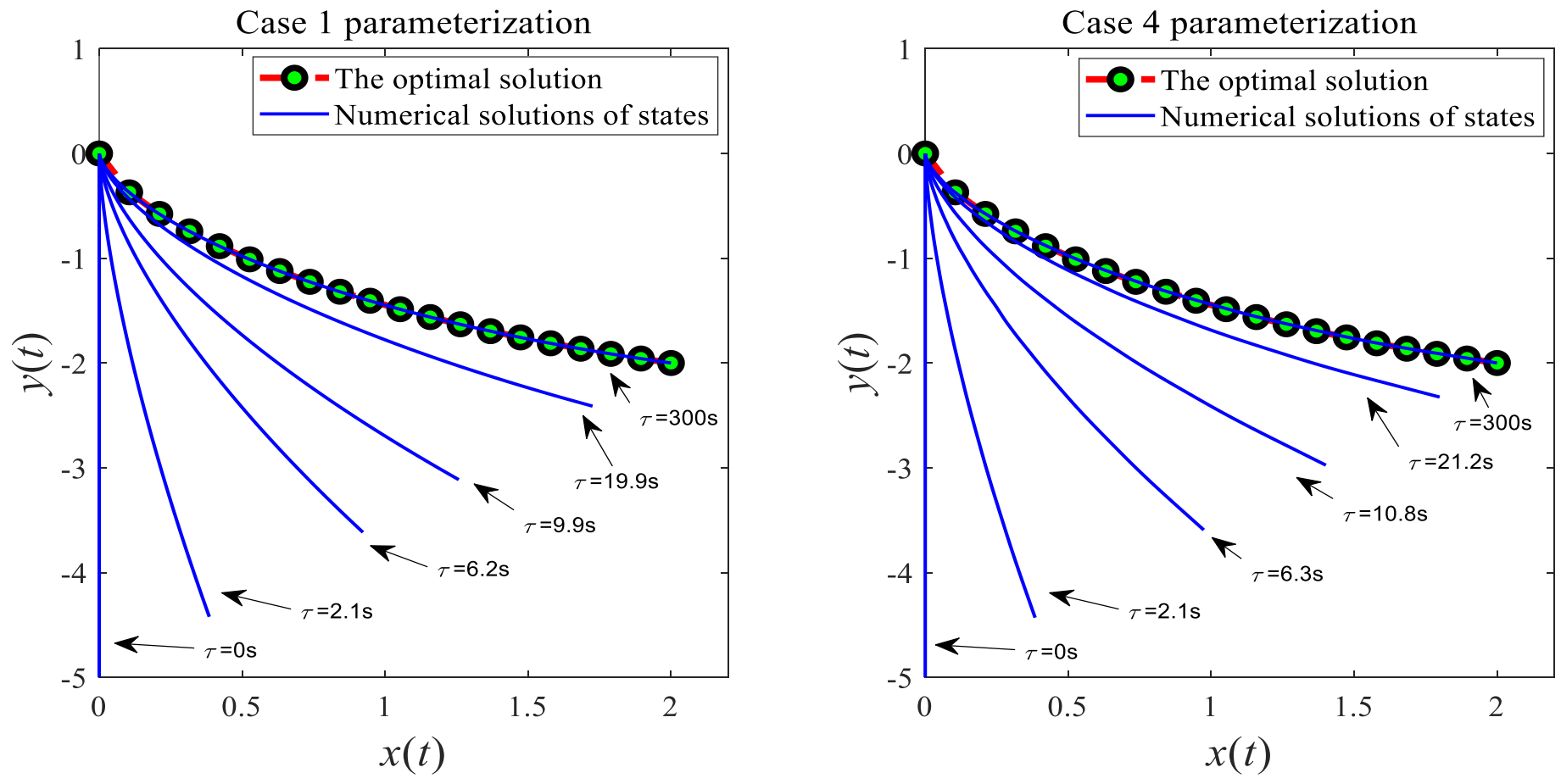

Fig. 15 The evolution of state solutions in the $xy$ coordinate plane for the Case 1 and Case 4 parameterization.

Based on the control and state solutions, the implicit costates and the Lagrange multipliers may be reconstructed. According to the previous analysis, as the control solution approaches the optimal solution, the costates and the multipliers will all tend to the optimal. Figure 16 gives the evolution profiles of the costates for the Case 1 parameterization. At $\tau$ =300s, they are identical with the optimal. In Fig. 17, the evolution profiles of the costates for the Case 4 parameterization are presented. However, there is some difference between the results and the optimal solution, especially in $\lambda_2$ . This phenomenon may arise from the error in the step approximation. In regards to the Lagrange multipliers, in Table 1 it may be found that the results from Cases 1-3 are the same, while the results of Case 4 is slightly different.

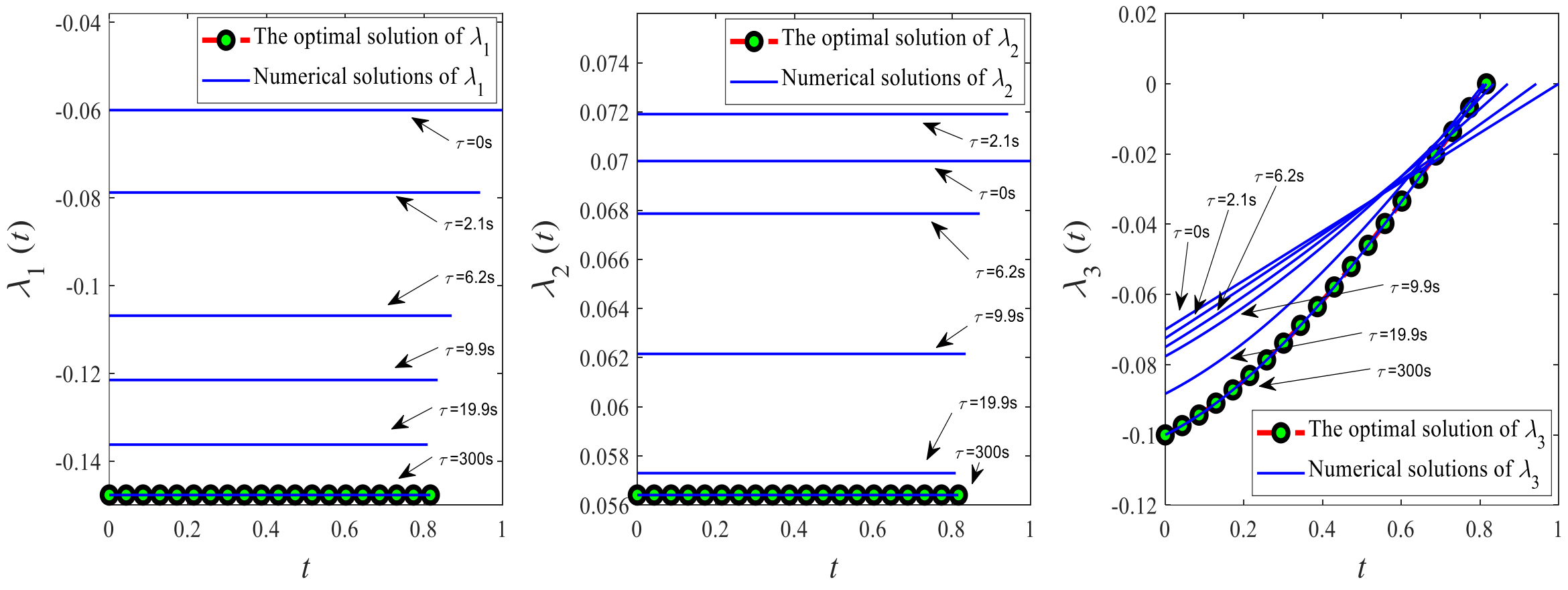

Fig. 16 The evolution of the costates $\lambda_1$, $\lambda_2$ and $\lambda_3$ for the Case 1 parameterization.

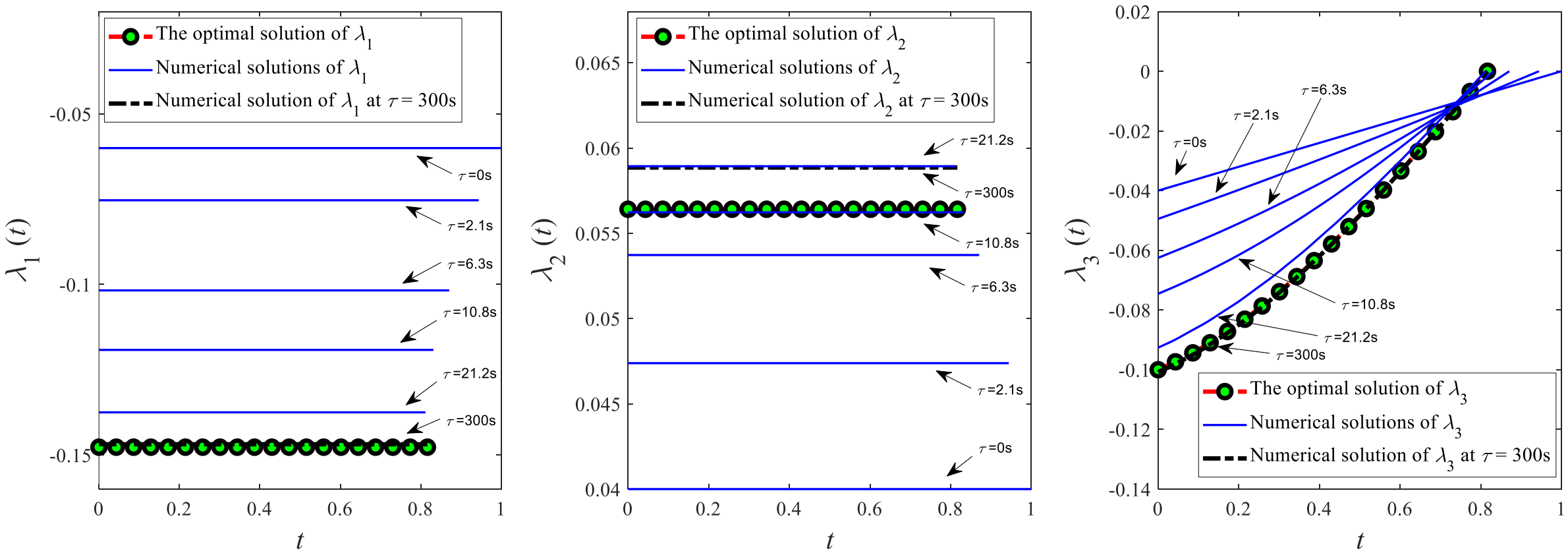

Fig. 17 The evolution of the costates $\lambda_1$, $\lambda_2$ and $\lambda_3$ for the Case 4 parameterization.

## VII. CONCLUSION

The direct shooting method is widely employed to solve the Optimal Control Problems (OCPs) due to its simplicity and ease of use. Compared with various other direct methods, this method has a small number of parameters to be determined. However, it has long been criticized for the lack of theoretical guarantee to achieve the optimality condition of the OCP, and is thought incapable to provide the information on the costates. This paper gives an enlightening endeavor to change this view. It is shown that the direct shooting method may produce the costates and its solution may also reach the optimal solution, upon the fundamental obtained from the Variation Evolving Method (VEM). With a reasonable parameterization, the control solution, the state solution and the costate solution may all approach the optimal solution. Thus, the direct shooting method may be regarded as a complete method in the view that the state and the costate solutions both converge.

In this paper, two forms of parameterization are investigated, and they are equivalent in results when endowed the same representation capacity. Both the global and the local parameterization are applicable within the proposed parameter evolution equations, which are consistent to the gradient dynamic equations, the continuous implementation of the discrete iterative gradient methods for the Nonlinear Programming (NLP) problems obtained by the direct shooting method. Illustrative examples are solved to demonstrate the effectiveness of the direct shooting method in seeking the optimal solutions, and it may have a higher efficiency than solving the Evolution Partial Differential Equation (EPDE) from the VEM since the scale of the transformed Initial-value Problem (IVP) is significantly reduced. Currently, the OCPs considered do not include the path constraints. It is more interesting to further consider the problems with complex path constraints, and the research will be carried out in the future.